\documentclass[11pt]{article}

\author{T.~Fiedler}

\usepackage{amsfonts, amssymb, theorem}

\title{Global knot theory in $F^2 \times \mathbb{R}$}
\newtheorem{propo}{Proposition}[section]
\newtheorem{lem}{Lemma}[section]
\newtheorem{theor}{Theorem}[section]
\newtheorem{defin}{Definition}[section]

\renewcommand{\Box}
{%
\mbox{}%
\nolinebreak%
\rule{2mm}{2mm}%
}

\newcommand{\stratapm}
{\makebox(18,10){%
\begin{picture}(16,8)(-10,4)%
\includegraphics{stratapm.pstex}%
\end{picture}%
\hspace{7mm}
}
}
\newcommand{\dr}
{\makebox(8,8)
{\begin{picture}(0,0)(3,5)%
\includegraphics{dr.pstex}%
\end{picture}
}
}
\newcommand{\br}
{\makebox(8,8)
{\begin{picture}(0,0)(3,5)%
\includegraphics{br.pstex}%
\end{picture}
}
}
\newcommand{\pmp}
{\,
\begin{picture}(0,0)(0,5)%
\includegraphics{pmp.pstex}%
\end{picture}
\hspace{7mm}
}
\newcommand{\trip}
{\,
\begin{picture}(0,0)(0,5)%
\includegraphics{trip.pstex}%
\end{picture}
\hspace{7mm}
}

\begin{document}
\maketitle
\tableofcontents
\begin{abstract}
\par We introduce a special class of knots, called global knots, in $F^2
\times \mathbb{R}$ and we construct new isotopy invariants, called
$T$-invariants, for global knots. 
\par Some $T$-invariants are of finite type but they cannot be extracted from
the generalized Kontsevitch integral (which is consequently not the universal
invariant of finite type for the restricted class of global knots).
\par We prove that $T$-invariants separate all global knots of a certain type. 
As a corollary, we prove the non-invertibility of some links in $S^3$ without
making use of the link group.
\end{abstract}
\section{Introduction and main results}
\hspace{7mm} Let $F^2$ be a compact oriented surface with or without boundary,
and let $v$ be a Morse-Smale vector field on $F^2$ which is transversal to the
boundary $\partial F^2$. (For us, a Morse-Smale vector field is a smooth 
vector field, having at most isolated non-degenerated singularities, and no 
limit cycles.)
We study oriented knots $K$ in the oriented manifold $F^2 \times \mathbb{R}$.
\par It turns out that there are naturally three types of knot theory, which we
call respectively local, global and general.

\subsection{Local knot theory}

\hspace{7mm} Let $F^2$ be the disk $D^2$ or the sphere $S^2$ and let $v$ be a
vector field which has only critical points of index +1. Alexander's theorem
says that each knot type (i.e. a knot up to ambient isotopy) has a
representative, also called $K$, such that the projection $K \hookrightarrow
F^2 \times \mathbb{R} \rightarrow F^2$ is transversal to $v$. Markov's theorem
says that two transversal representatives of the same knot type can be joined
by an almost transversal isotopy i.e. an isotopy through transversal knots,
besides for a finite number of values of the parameter where a Reidemeister
move of type $I$ occurs in a singularity of $v$ (such a move is usually called
Markov move). This type of knot theory is traditionally the most studied one
and we cannot add anything new here.

\subsection{Global knot theory}

\hspace{7mm} A knot type $K$ is called a {\em global knot\/} if there is a
vector field $v$ without critical points or only with critical points of index
-1, and there is a representative $K$ such that the projection $K
\hookrightarrow F^2 \times \mathbb{R} \rightarrow F^2$ is transversal to $v$. We call such $v$ {\em non-elliptic\/} vector fields. Notice that a knot
$K \hookrightarrow D^2 \times \mathbb{R} \hookrightarrow F^2 \times
\mathbb{R}$ can never be a global knot. This implies that there is no analogue
of Alexander's theorem here. However, there is an analogue of Markov's theorem
for global knots which is even better than Markov's theorem:
{\em If two representatives of the same global knot have projections 
transversal to the same non-elliptic vector field $v$, then there is an 
isotopy between them transversal to $v$.\/} 
We give a proof in a special case and we indicate the idea of the proof in the general case.
\par Global knots are the main object of our work. They have almost never been studied before, except for the very special case of closed braids.

\subsection{General knot theory}

\hspace{7mm} A knot type $K$ is called a {\em general knot\/} if for each
representative $K$ the following property is verified: if $v$ is a vector
field such that $K \hookrightarrow F^2 \times \mathbb{R} \rightarrow F^2$ is
transversal to $v$, then $v$ has critical points of both indices $+1$ and
$-1$. General knots do exist. For example: the Whitehead link seen as a knot
in the solid torus. Evidently, if a knot is neither local nor global, then it
is general. We prove that in the general setting there is no analogue of
Markov's theorem (in difference to the global setting). More precisely, we
give two representatives of some general knot in the solid torus with
projections transversal to the same vector field $v$, and we show that they
cannot be joined by any transversal isotopy. This indicates that the general
case is even much more complicated than the local and the global case.
\par The main achievement of our work is the construction of new isotopy
invariants, called {\em T-invariants\/} ($T$ means "transversal"), for global
knots. These invariants depend neither on the chosen non-elliptic vector field
$v$, nor on the chosen representative of the knot whose projection is
transversal to $v$. Hence, $T$-invariants are isotopy invariants in the usual
sense. $T$-invariants are defined as "Gauss diagram invariants". Consequently,
their calculation has polynomial complexity with respect to the number of
crossings of the knot diagrams. However, not all $T$-invariants are of finite
type in the sense of Vassiliev. Moreover, we show that even some
$T$-invariants of finite type cannot be extracted from the generalized
Kontsevitch integral (see [A-M-R]). This comes from the fact that
$T$-invariants are well defined only for global knots and not for all knots in
$F^2 \times \mathbb{R}$.
\par A knot $K \hookrightarrow F^2 \times \mathbb{R}$ is called a {\em solid
  torus knot\/} (or a closed braid) if it has a representative whose
  projection is contained in some annulus $S^1 \times I \hookrightarrow F^2$
  (of course, local knots are a special case of solid torus knots). We
  conjecture that $T$-invariants separate all global knots in $F^2 \times \mathbb{R}$ which are not solid torus knots.
\par We prove this conjecture in the following special case: Let $T^2$ be the
  torus. An oriented global knot $K \hookrightarrow T^2 \times\mathbb{R}$ is
  called $\mathbb{Z}/2\mathbb{Z}$-pure if for each crossing of $K$, each of
  the two oriented loops obtained by splitting the crossing is non-trivial in $H_1(T^2;\mathbb{Z})/\langle [K] \rangle \bigotimes \mathbb{Z}/2\mathbb{Z}\cong \mathbb{Z}/2\mathbb{Z}$
(We consider $K$ as a diagram in $F^2 \times \mathbb{R}$ over $F^2$ and we
  denote by $[K]$ the homology class represented by $K$. If $K$ is not a solid
  torus knot then $H_1(T^2;\mathbb{Z})/\langle [K] \rangle \bigotimes
  \mathbb{Z}/2\mathbb{Z}$ is automatically isomorphic to $\mathbb{Z}/2\mathbb{Z}$.)
\begin{theor}
{\it $T$-invariants separate all $\mathbb{Z}/2\mathbb{Z}$-pure global
knots in $T^2 \times \mathbb{R}$.\/}
\end{theor}
\par Let $flip: T^2 \times \mathbb{R} \to T^2 \times \mathbb{R}$ be the
hyper-elliptic involution on $T^2$ multiplied by the identity on the lines
$\mathbb{R}$. An oriented knot $K \hookrightarrow T^2 \times \mathbb{R}$ is
called {\em invertible\/} if it is ambient isotopic to
$flip(-K)=-flip(K)$. One easily shows that both the generalized HOMFLY-PT and
the generalized Kauffman polynomial can never distinguish $K$ from
$flip(-K)$. (The polynomials of the cables do not make the distinction
either.) On the other hand, we prove the non-invertibility of some global
knots in $T^2 \times \mathbb{R}$ using $T$-invariants of finite
type. Consequently, these invariants cannot be extracted from the above knot
polynomials. 
\par Knots in $T^2 \times \mathbb{R}$ are in 1-1 correspondence with ordered 3-component links in $S^3$ containing the Hopf link $H$ as a sublink. Using a $T$-invariant of degree 6 for the knot $K$, we show that the link $L=K \cup H \hookrightarrow S^3$ (see Fig.~1). is not invertible for any chosen orientation on it. Notice that this is the first proof of the non-invertibility of a link with a numerical invariant, which does not make any use of the link group $\pi_1(S^3\setminus $K$;*)$. (Compare with $L$ for approaches which use the link group.)
$$
\begin{picture}(0,0)%
\includegraphics{im1.pstex}%
\end{picture}%
\setlength{\unitlength}{4144sp}%
\begingroup\makeatletter\ifx\SetFigFont\undefined%
\gdef\SetFigFont#1#2#3#4#5{%
  \reset@font\fontsize{#1}{#2pt}%
  \fontfamily{#3}\fontseries{#4}\fontshape{#5}%
  \selectfont}%
\fi\endgroup%
\begin{picture}(4506,3244)(264,-2535)
\put(2601,-744){\makebox(0,0)[lb]{\smash{\SetFigFont{12}{14.4}{\rmdefault}{\mddefault}{\updefault}$H$}}}
\put(4770,-1204){\makebox(0,0)[lb]{\smash{\SetFigFont{12}{14.4}{\rmdefault}{\mddefault}{\updefault}$K$}}}
\put(2510,-2477){\makebox(0,0)[lb]{\smash{\SetFigFont{12}{14.4}{\rmdefault}{\mddefault}{\updefault}$L=K \cup H$}}}
\end{picture}
$$
\begin{center}
{\bf Fig. 1}
\end{center}
\vspace{1cm}

\par Let us consider the space of all diagrams of a given knot type $K
\hookrightarrow F^2 \times \mathbb{R}$. The discriminant is the subspace of
all non-generic diagrams. Each path in the space of diagrams which cuts the
discriminant only in strata of codimension 1 is a generic isotopy of knots
(see [F], sect.~1). 
\par The construction of $T$-invariants relies on the combination of two
different approaches. On one hand, there is the concept of $G$-pure knots and
$G$-pure isotopy. A $G$-pure isotopy is an isotopy which does not cut the
discriminant in strata with certain homological markings of the
crossings. $T$-invariants are, roughly speaking, Gauss diagram formulas which
are invariant under $G$-pure isotopy. The problem is now, that even if two
knots are isotopic, we cannot grant that there exists a $G$-pure isotopy
joining them. 
\par At this place, the concept of global knots is introduced. If two
isotopic knots are global, then there exists an isotopy through global knots
between them. Such an isotopy does not cut the strata of the discriminant depicted in Fig.~2.
$$
\begin{picture}(0,0)%
\includegraphics{ima2.pstex}%
\end{picture}%
\setlength{\unitlength}{4144sp}%
\begingroup\makeatletter\ifx\SetFigFont\undefined%
\gdef\SetFigFont#1#2#3#4#5{%
  \reset@font\fontsize{#1}{#2pt}%
  \fontfamily{#3}\fontseries{#4}\fontshape{#5}%
  \selectfont}%
\fi\endgroup%
\begin{picture}(4551,1049)(464,-728)
\end{picture}
$$
\begin{center}
{\bf Fig. 2}
\end{center}
\par We still do not know wether isotopic $G$-pure global knots can be joined
by a $G$-pure isotopy. But for an isotopy through global knots, we have enough
control over the "cycles of crossings" (compare with [F], sect.~4) to enable
us to show that $T$-invariants are actually isotopy invariants! 
\par Using $T$-invariants, we show in sect.~9 that there exist isotopic 
$G$-pure general knots such that
\begin{enumerate}
\item
there is no $G$-pure isotopy between them
\item
there is no isotopy transversal to the vector field $v$ between them
\end{enumerate}
\newpage
\section{Isotopy of global knots}

\hspace{7mm} We fix orientations on $F^2$ and $F^2 \times \mathbb{R}$. Let
$pr: F^2 \times \mathbb{R} \to F^2$ be the canonical projection, and let $v$
be a non-elliptic vector field. We fix the orientation on the global knot $K
\hookrightarrow F^2 \times \mathbb{R}$ in such a way that ($pr(K), v$) induce the given orientation on $F^2$.
\par The natural equivalence relation for global knots is: ambient isotopy
through global knots keeping the vector field $v$ fixed. However, we believe
that this equivalence relation coincides with the usual ambient isotopy.

{\sl Conjecture}. {\it Two global knots with respect to the same vector field
  $v$ are ambient isotopic if and only if they are ambient isotopic through global knots with respect to $v$.\/}

\begin{description}

\item[Remarks:]

\begin{enumerate}

\item

The conjecture is true in the case of closed braids: this is a consequence of
 Artin's theorem.

\item 

Below we give a proof in the case where $F^2=T^2$.

\item

We outline the strategy of the proof in the general case (and we will come
back to it in another paper): One can easily prove in a geometrical way that
braids are isotopic as tangles if and only if they are isotopics as braids
(see e.g. [P-S]). This proof can be generalized to all such isotopies of global
 knots which do not pass through the singularities of $v$. If an isotopy
 passes through a singularity of $v$, then the resulting knot is no longer
 transversal to $v$, except for the case where the isotopy passes again
 through the same singularity but in the opposite direction. This can be
 proven by constructing 2-disks with piecewise smooth boundary in the trace of
 the isotopy on $F^2$. The vector field $v$, having only critical points of
 index -1, can never be transversal to such a disk. Using this fact, one can remove successively these 2-disks in the isotopy.

\end{enumerate}

\end{description}

\par We consider now the case of the torus $T^2$. Let $pr_2: S^1_1 \times
S^1_2 \to S^1_2$ be the projection and let $v$ be the unit vector field
tangential to the fibres of $pr_2$ (we might perturbate $v$ slightly so that
it has no longer any closed orbits). Consequently, a knot $K \hookrightarrow
T^2 \times \mathbb{R}$ is global with respect to $v$ if and only if the restriction $pr_2: K \to S^1_2$ is a covering.
\begin{lem}
{\it Isotopic global knots with respect to $v$ in $T^2 \times \mathbb{R}$
are isotopic through global knots with respect to $v$.\/}
\end{lem}
\par Artin's theorem implies that isotopic closed braids are isotopic through
closed braids (see e.g. [M]). Let $T^2 \times \mathbb{R} \hookrightarrow S^3$
be the standard embedding and let $A_1$ resp. $A_2$ be the cores of the solid
tori $S^3 \setminus (T^2 \times \mathbb{R})$. The oriented link $A_1 \cup A_2$
is determined by the following conventions:
$lk(A_1, \{*\} \times S^1_2) = lk(A_2, S^1_1 \times \{*\}) = 1$
and $K$ is a global knot in $T^2 \times \mathbb{R}$ with respect to $v$ if and
only if $K \cup A_2$ is a closed braid for the natural fibering $S^3 \setminus
A_1 \to S^1_2$. Consequently, if $K$ and $K'$ are isotopic global knots then
$K \cup A_2$ and $K' \cup A_2$ are closed braids which are isotopic as
links. According to Artin's theorem, they are also isotopic as closed
braids. But $A_2$ is just the closure of a 1-string braid and therefore we may
assume that it remains fixed in the isotopy of closed braids. This implies
that the isotopy of $K$ in $S^3 \setminus (A_1 \cup A_2)$ is an isotopy
throughglobal knots $\Box$.
\par From now on, all the isotopies considered will be isotopies through global knots.
\par {\sl Basic observation.} {\it In an isotopy of global knots, the
  Reidemeister moves depicted in Fig.~3 can never occur.\/}

$$
\begin{picture}(0,0)%
\includegraphics{im2.pstex}%
\end{picture}%
\setlength{\unitlength}{4144sp}%
\begingroup\makeatletter\ifx\SetFigFont\undefined%
\gdef\SetFigFont#1#2#3#4#5{%
  \reset@font\fontsize{#1}{#2pt}%
  \fontfamily{#3}\fontseries{#4}\fontshape{#5}%
  \selectfont}%
\fi\endgroup%
\begin{picture}(5188,3464)(203,-2864)
\end{picture}
$$
\begin{center}
{\bf Fig. 3}
\end{center}
{\bf Proof}. Obviously, in the singularities there is always one branch of $K$
which has not the correct orientation (see the middle part of Fig.~3).
$\Box$

{\sl Question}. {\it Is there an analogue of Alexander's theorem (in the right
  sense) for %
  knots in $T^2 \times \mathbb{R}$?\/} More precisely, given a diagram $K$, is
there an algorithm which either constructs some non-elliptic vector field $v$
and a representative of $K$ transversal to $v$, or otherwise, which proves
that $K$ is not a global knot?

\newpage
\section{Construction of $T$-invariants for global knots.}

\par Let $K_0 \hookrightarrow F^2 \times \mathbb{R}$ be an oriented global
knot and let $K_t, t \in [0,1]$ be an isotopy of $K_0$ (through global knots). Let $G$ be a fixed quotient group of $H_1(F^2; \mathbb{Z})$ and let $[K]_G$ be the homology class in $G$ represented by
$K=K_0$.

\begin{defin} {\rm A global knot $K$ is {\it $G$-pure\/} if for each crossing
    $p$ of $K$ (with respect to the projection $pr$), $[K^+_p]_G \notin \{0,
    \pm [K]_G \}$. The isotopy $K_t$ is {\em $G$-pure\/} if $K_t$ is $G$-pure
    for each $t$. A knot type (called $K$ as the knot itself) is {\em
    $G$-pure\/} if it has a global representative which is $G$-pure.\/}
\end{defin}
\par We remind the definition of the oriented (global) knot $K^+_p$ in Fig.~4 (see also [F], sect.~ 0 and 1).
$$
\begin{picture}(0,0)%
\includegraphics{im3.pstex}%
\end{picture}%
\setlength{\unitlength}{4144sp}%
\begingroup\makeatletter\ifx\SetFigFont\undefined%
\gdef\SetFigFont#1#2#3#4#5{%
  \reset@font\fontsize{#1}{#2pt}%
  \fontfamily{#3}\fontseries{#4}\fontshape{#5}%
  \selectfont}%
\fi\endgroup%
\begin{picture}(3894,1698)(439,-1289)
\put(1120,-689){\makebox(0,0)[lb]{\smash{\SetFigFont{12}{14.4}{\rmdefault}{\mddefault}{\updefault}$p$}}}
\put(1147,-1232){\makebox(0,0)[lb]{\smash{\SetFigFont{12}{14.4}{\rmdefault}{\mddefault}{\updefault}$K$}}}
\put(4316,214){\makebox(0,0)[lb]{\smash{\SetFigFont{12}{14.4}{\rmdefault}{\mddefault}{\updefault}$K_p^-$}}}
\put(3516,-1231){\makebox(0,0)[lb]{\smash{\SetFigFont{12}{14.4}{\rmdefault}{\mddefault}{\updefault}$K_p^+$}}}
\end{picture}
$$
\begin{center}
{\bf Fig. 4}
\end{center}
\par The {\it Gauss diagram\/} of $K$ is the abstract chord diagram of
$pr(K)$, where each chord $p$ is oriented from the undercross to the overcross
of the crossing $p$, and it is marked by the writhe $w(p)$ and by the homology
class $[K^+_p]_G \in G$. A {\it configuration\/} is a given (abstract) chord
diagram with given orientations and homological markings (in $G$) of the
chords, but without writhes. A {\it Gauss sum of degree $d$ for the knot
  $K$\/} (also called a {\it Gauss diagram formula\/}) is a sum which is
defined in the following way: Let $D$ be a given configuration of $d$ chords. We consider the integer

\begin{displaymath}
\sum_{D} {\mbox{function (writhes of the crossings of $K$ corresponding to the chords of $D$)}}
\end{displaymath}
where $D$ runs over all the subdiagrams of the Gauss diagram of $K$. The
function is called the {\em weight function\/} (see also [F], sect.~0 and 1).

\begin{defin}
{\rm A {\it $G$-pure configuration of degree $m$\/} is a given configuration
  $D$ of $m$ oriented chords $(p_1, \dots, p_m)$ with corresponding markings
  $(a_1, \dots, a_m)$ where each $a_i \in G \setminus \{0, \pm [K]_G \}$,
  verifying the following 2 conditions:}
\begin{enumerate}
\item
{\rm Let $D$ be represented as a subdiagram $D_0$ of a Gauss diagram of any
  $G$-pure knot $K_0$, and let $K_t, t \in [0,1]$ be any $G$-pure isotopy of
  $K_0$ without Reidemeister moves of type $II$ involving one of the crossings
  $p_i$. Then, $D$ is preserved as a subdiagram of $K_t$ i.e. there is a
  continuous family $D_t$ of subdiagrams of the Gauss diagrams of $K_t$ and
  almost each $D_t$ represents $D$.\/} 
\item
{\rm If $D$ contains a fragment of the type depicted in Fig.~5, then $a_i
\not= a_j$.\/}
\end{enumerate}
\end{defin}
$$
\begin{picture}(0,0)%
\includegraphics{manu20.pstex}%
\end{picture}%
\setlength{\unitlength}{4144sp}%
\begingroup\makeatletter\ifx\SetFigFont\undefined%
\gdef\SetFigFont#1#2#3#4#5{%
  \reset@font\fontsize{#1}{#2pt}%
  \fontfamily{#3}\fontseries{#4}\fontshape{#5}%
  \selectfont}%
\fi\endgroup%
\begin{picture}(2302,2421)(645,-1707)
\put(1131,519){\makebox(0,0)[lb]{\smash{\SetFigFont{12}{14.4}{\rmdefault}{\mddefault}{\updefault}$a_i$}}}
\put(2406,519){\makebox(0,0)[lb]{\smash{\SetFigFont{12}{14.4}{\rmdefault}{\mddefault}{\updefault}$a_j$}}}
\put(1411,-116){\makebox(0,0)[rb]{\smash{\SetFigFont{12}{14.4}{\rmdefault}{\mddefault}{\updefault}$p_i$}}}
\put(2251,-16){\makebox(0,0)[lb]{\smash{\SetFigFont{12}{14.4}{\rmdefault}{\mddefault}{\updefault}$p_j$}}}
\end{picture}
$$
\begin{center}
{\bf Fig. 5}
\end{center}
{\it Exemple. $2.1$}
$G:= \mathbb{Z}$ and $[K]_G = 0$. Let $a \in \mathbb{Z} \setminus 0$ be
fixed. There is exactly one $\mathbb{Z}$-pure configuration of degree 1 which
involves the class $a$ (see Fig.~6). There are exactly three $\mathbb{Z}$-pure
configurations of degree 2 which involve the class $a$ (see Fig.~7).
$$
\begin{picture}(0,0)%
\includegraphics{manu21.pstex}%
\end{picture}%
\setlength{\unitlength}{4144sp}%
\begingroup\makeatletter\ifx\SetFigFont\undefined%
\gdef\SetFigFont#1#2#3#4#5{%
  \reset@font\fontsize{#1}{#2pt}%
  \fontfamily{#3}\fontseries{#4}\fontshape{#5}%
  \selectfont}%
\fi\endgroup%
\begin{picture}(1164,1316)(739,-617)
\put(1400,-110){\makebox(0,0)[lb]{\smash{\SetFigFont{12}{14.4}{\rmdefault}{\mddefault}{\updefault}$p_1$}}}
\put(1361,534){\makebox(0,0)[lb]{\smash{\SetFigFont{12}{14.4}{\rmdefault}{\mddefault}{\updefault}$a$}}}
\end{picture}
$$
\begin{center}
{\bf Fig. 6}
\end{center}
\newpage
$$
\begin{picture}(0,0)%
\includegraphics{manu21bis.pstex}%
\end{picture}%
\setlength{\unitlength}{4144sp}%
\begingroup\makeatletter\ifx\SetFigFont\undefined%
\gdef\SetFigFont#1#2#3#4#5{%
  \reset@font\fontsize{#1}{#2pt}%
  \fontfamily{#3}\fontseries{#4}\fontshape{#5}%
  \selectfont}%
\fi\endgroup%
\begin{picture}(4948,1420)(465,-796)
\put(821,-156){\makebox(0,0)[lb]{\smash{\SetFigFont{12}{14.4}{\rmdefault}{\mddefault}{\updefault}$p_1$}}}
\put(2676,-37){\makebox(0,0)[lb]{\smash{\SetFigFont{12}{14.4}{\rmdefault}{\mddefault}{\updefault}$p_1$}}}
\put(796,429){\makebox(0,0)[lb]{\smash{\SetFigFont{12}{14.4}{\rmdefault}{\mddefault}{\updefault}$a$}}}
\put(1351,-796){\makebox(0,0)[lb]{\smash{\SetFigFont{12}{14.4}{\rmdefault}{\mddefault}{\updefault}$-a$}}}
\put(3236,-116){\makebox(0,0)[lb]{\smash{\SetFigFont{12}{14.4}{\rmdefault}{\mddefault}{\updefault}$p_2$}}}
\put(2496,-751){\makebox(0,0)[lb]{\smash{\SetFigFont{12}{14.4}{\rmdefault}{\mddefault}{\updefault}$a$}}}
\put(3196,424){\makebox(0,0)[lb]{\smash{\SetFigFont{12}{14.4}{\rmdefault}{\mddefault}{\updefault}$-a$}}}
\put(4603,-35){\makebox(0,0)[lb]{\smash{\SetFigFont{12}{14.4}{\rmdefault}{\mddefault}{\updefault}$p_1$}}}
\put(4436,439){\makebox(0,0)[lb]{\smash{\SetFigFont{12}{14.4}{\rmdefault}{\mddefault}{\updefault}$a$}}}
\put(5071,459){\makebox(0,0)[lb]{\smash{\SetFigFont{12}{14.4}{\rmdefault}{\mddefault}{\updefault}$a$}}}
\put(5176,-91){\makebox(0,0)[lb]{\smash{\SetFigFont{12}{14.4}{\rmdefault}{\mddefault}{\updefault}$p_2$}}}
\put(1376,-161){\makebox(0,0)[lb]{\smash{\SetFigFont{12}{14.4}{\rmdefault}{\mddefault}{\updefault}$p_2$}}}
\put(1896,-611){\makebox(0,0)[lb]{\smash{\SetFigFont{14}{16.8}{\rmdefault}{\mddefault}{\updefault},}}}
\put(3856,-621){\makebox(0,0)[lb]{\smash{\SetFigFont{14}{16.8}{\rmdefault}{\mddefault}{\updefault},}}}
\end{picture}
$$\nolinebreak
\begin{center}
{\bf Fig. 7}
\end{center}

\par This follows immediately from the fact that a $G$-pure isotopy does not
intersect the following strata in the discriminant: $a_{\dr}^{+(-)}(0|a,0)$,
$a_{\br}^{+(-)}(0|a,-a)$, $a_{\dr}^{+(-)}(a|a,0)$, $a_{\br}^{+(-)}(a|0,a)$,  (See [F], sect.~4.11).
\par In fact, in [F] we replaced all triple points of type \trip  by triple
points of type \pmp. A closer look to this replacement shows that a
Reidemeister move of type $III$ corresponding to \trip is $G$-pure and a
configuration $D$ is invariant under this move if and only if it is invariant
under the corresponding move of type \pmp. We illustrate this in Fig.~8.
$$
\begin{picture}(0,0)%
\includegraphics{im4.pstex}%
\end{picture}%
\setlength{\unitlength}{4144sp}%
\begingroup\makeatletter\ifx\SetFigFont\undefined%
\gdef\SetFigFont#1#2#3#4#5{%
  \reset@font\fontsize{#1}{#2pt}%
  \fontfamily{#3}\fontseries{#4}\fontshape{#5}%
  \selectfont}%
\fi\endgroup%
\begin{picture}(3584,2021)(465,-1505)
\put(1141,-146){\makebox(0,0)[lb]{\smash{\SetFigFont{12}{14.4}{\rmdefault}{\mddefault}{\updefault}$a$}}}
\put(1326,-1505){\makebox(0,0)[lb]{\smash{\SetFigFont{12}{14.4}{\rmdefault}{\mddefault}{\updefault}$a$}}}
\put(1681,-1210){\makebox(0,0)[lb]{\smash{\SetFigFont{12}{14.4}{\rmdefault}{\mddefault}{\updefault}$a$}}}
\put(3538,-1366){\makebox(0,0)[lb]{\smash{\SetFigFont{12}{14.4}{\rmdefault}{\mddefault}{\updefault}$a$}}}
\put(3506,145){\makebox(0,0)[lb]{\smash{\SetFigFont{12}{14.4}{\rmdefault}{\mddefault}{\updefault}$a$}}}
\put(2281,-670){\makebox(0,0)[lb]{\smash{\SetFigFont{12}{14.4}{\rmdefault}{\mddefault}{\updefault}$a$}}}
\put(3171,-290){\makebox(0,0)[b]{\smash{\SetFigFont{12}{14.4}{\rmdefault}{\mddefault}{\updefault}$a$}}}
\put(1956,-860){\makebox(0,0)[lb]{\smash{\SetFigFont{12}{14.4}{\rmdefault}{\mddefault}{\updefault}$a$}}}
\end{picture}
$$
\begin{center}
{\bf Fig. 8}
\end{center}
\par The homological markings involved in \trip are exactly the same as those
involved in the corresponding \pmp. Notice that the $G$-pure configuration
depicted in the left-hand part of Fig.~9, which appears in the replacement,
disappears again. 
\par Consequently, the chords $p_1$ and $p_2$ of the above configurations
cannot become crossed (as shown in the right-hand part of Fig.~9) in a
$G$-pure isotopy.

$$
\begin{picture}(0,0)%
\includegraphics{ima9.pstex}%
\end{picture}%
\setlength{\unitlength}{4144sp}%
\begingroup\makeatletter\ifx\SetFigFont\undefined%
\gdef\SetFigFont#1#2#3#4#5{%
  \reset@font\fontsize{#1}{#2pt}%
  \fontfamily{#3}\fontseries{#4}\fontshape{#5}%
  \selectfont}%
\fi\endgroup%
\begin{picture}(3609,1282)(469,-878)
\put(656,219){\makebox(0,0)[lb]{\smash{\SetFigFont{12}{14.4}{\rmdefault}{\mddefault}{\updefault}$a$}}}
\put(1291,239){\makebox(0,0)[lb]{\smash{\SetFigFont{12}{14.4}{\rmdefault}{\mddefault}{\updefault}$a$}}}
\put(2116,-711){\makebox(0,0)[lb]{\smash{\SetFigFont{14}{16.8}{\rmdefault}{\mddefault}{\updefault},}}}
\put(3696,-183){\makebox(0,0)[lb]{\smash{\SetFigFont{12}{14.4}{\rmdefault}{\mddefault}{\updefault}$p_2$}}}
\put(3113,-178){\makebox(0,0)[lb]{\smash{\SetFigFont{12}{14.4}{\rmdefault}{\mddefault}{\updefault}$p_1$}}}
\end{picture}
$$
\begin{center}
{\bf Fig.9}
\end{center}
{\it Example $3.2$}
$G := \mathbb{Z}/2\mathbb{Z}$ and $[K]_G=0$. Each chord is marked by the
non-trivial element in $\mathbb{Z}/2\mathbb{Z}$. In this case, each
configuration which does not contain a fragment as depicted in Fig.~10 
$$
\begin{picture}(0,0)%
\includegraphics{ima10.pstex}%
\end{picture}%
\setlength{\unitlength}{4144sp}%
\begingroup\makeatletter\ifx\SetFigFont\undefined%
\gdef\SetFigFont#1#2#3#4#5{%
  \reset@font\fontsize{#1}{#2pt}%
  \fontfamily{#3}\fontseries{#4}\fontshape{#5}%
  \selectfont}%
\fi\endgroup%
\begin{picture}(2357,2242)(610,-1707)
\put(1116,-136){\makebox(0,0)[lb]{\smash{\SetFigFont{12}{14.4}{\rmdefault}{\mddefault}{\updefault}$p_i$}}}
\put(2201,-51){\makebox(0,0)[lb]{\smash{\SetFigFont{12}{14.4}{\rmdefault}{\mddefault}{\updefault}$p_j$}}}
\end{picture}
$$
\begin{center}
{\bf Fig. 10}
\end{center}
is a $\mathbb{Z}/2\mathbb{Z}$-pure configuration.
\par Let $D$ be a fixed $G$-pure configuration of degree $m$.

\begin{defin}
 {\rm Let $\mathcal{D}_i, i \in I$ be a finite set of configurations, each of
 them with at most $m+n$ oriented chords $(p_1, \dots, p_m, q^{(i)}_1, \dots,
 q^{(i)}_{n_i})$. Let $n = \max_{i}(n_i)$. All chords have markings in $G
 \setminus \{0, \pm[K]_G \}$. The given (unordered) chords $(p_1, \dots, p_m)$
 form the given $G$-pure subconfiguration $D$ in each $\mathcal{D}_i$. Let $f_i, i \in I$ be functions
\begin{displaymath}
f_i: \underbrace{\mathbb{Z}/2\mathbb{Z} \times \cdots \times \mathbb{Z}/
2\mathbb{Z}}_{n_i} \rightarrow \mathbb{Z}
\end{displaymath}
The linear combination of Gauss diagram formulas
\begin{displaymath}
c(D) := \sum_{\mathcal{D}_i}{f_i(w(q^{(i)}_1), \dots, w(q^{(i)}_{n_i}))}
\end{displaymath}
is called a {\it $G$-pure class of $D$ of degree at most $n$ } if the
following condition holds: 
\par Let $D$ be represented as a subdiagram $D_0$ of a Gauss diagram of any
$G$-pure knot $K_0$ and let $K_t, t \in [0,1]$ be any $G$-pure isotopy of
$K_0$ without Reidemeister moves of type $II$ involving one of the
$p_i$. (Hence, there exists a continuous family $D_t$ as in Def.~3.2.). Let

\begin{displaymath}
c(D_t) := \sum_{i} \sum_{\Delta_i}{f_i(w(q^({i})_1, \dots,
w_(q^({i})_{n_i}))}
\end{displaymath}
where $\Delta_i$ runs through all subdiagrams which represent $\mathcal{D}_i$
in the Gauss diagram of $K_t$, and which contain $D_t$ as the given
subconfiguration $D$. The integer $c(D_t)$ is the same for all $t \in [0,1]$
such that the projection $pr: K_t \to F^2$ is generic. (Here, $w(q^{(i)}_j)$
is the writhe of the crossing $q^{(i)}_j$.)}
\end{defin}
{\bf Remarks}:
\begin{enumerate}
\item
In the calculation of $c(D_t)$, subdiagrams which coincide up to different
numerations of the chords are identified, 
$$
\begin{picture}(0,0)%
\includegraphics{ima11.pstex}%
\end{picture}%
\setlength{\unitlength}{4144sp}%
\begingroup\makeatletter\ifx\SetFigFont\undefined%
\gdef\SetFigFont#1#2#3#4#5{%
  \reset@font\fontsize{#1}{#2pt}%
  \fontfamily{#3}\fontseries{#4}\fontshape{#5}%
  \selectfont}%
\fi\endgroup%
\begin{picture}(5531,2405)(388,-2164)
\put(4295,-935){\makebox(0,0)[lb]{\smash{\SetFigFont{12}{14.4}{\rmdefault}{\mddefault}{\updefault}$p$}}}
\put(4775,-510){\makebox(0,0)[lb]{\smash{\SetFigFont{12}{14.4}{\rmdefault}{\mddefault}{\updefault}$q_2$}}}
\put(4810,-1510){\makebox(0,0)[lb]{\smash{\SetFigFont{12}{14.4}{\rmdefault}{\mddefault}{\updefault}$q_1$}}}
\put(1215,-950){\makebox(0,0)[lb]{\smash{\SetFigFont{12}{14.4}{\rmdefault}{\mddefault}{\updefault}$p$}}}
\put(1808,-1298){\makebox(0,0)[lb]{\smash{\SetFigFont{12}{14.4}{\rmdefault}{\mddefault}{\updefault}$q_2$}}}
\put(1707,-505){\makebox(0,0)[lb]{\smash{\SetFigFont{12}{14.4}{\rmdefault}{\mddefault}{\updefault}$q_1$}}}
\put(3076,-1066){\makebox(0,0)[lb]{\smash{\SetFigFont{12}{14.4}{\rmdefault}{\mddefault}{\updefault}$=$}}}
\end{picture}
$$
\begin{center}
{\bf Fig. 11}
\end{center}
so that they bring only one term intothe sum (see Fig~11).
\item
We say that the class $c(D)$ is of degree $n$ if for each class $c'(D)$ which uses no more than $n-1$ chords, $c(D_t) - c'(D_t)$ is not constant.
\end{enumerate}

{\it Example $3.3$}
$G := \mathbb{Z}$, $[K]_G=0$, $a \in \mathbb{Z} \setminus 0$, $D$,
$\mathcal{D}_1$, $\mathcal{D}_2$ are as in Fig.~12, $f_1 \equiv f_2 = w(q_1)w(q_2)$.
$$
\begin{picture}(0,0)%
\includegraphics{ima12.pstex}%
\end{picture}%
\setlength{\unitlength}{4144sp}%
\begingroup\makeatletter\ifx\SetFigFont\undefined%
\gdef\SetFigFont#1#2#3#4#5{%
  \reset@font\fontsize{#1}{#2pt}%
  \fontfamily{#3}\fontseries{#4}\fontshape{#5}%
  \selectfont}%
\fi\endgroup%
\begin{picture}(5712,1630)(423,-1232)
\put(5609,-1232){\makebox(0,0)[lb]{\smash{\SetFigFont{12}{14.4}{\rmdefault}{\mddefault}{\updefault}$-a$}}}
\put(5136,148){\makebox(0,0)[lb]{\smash{\SetFigFont{12}{14.4}{\rmdefault}{\mddefault}{\updefault}$a$}}}
\put(5654,148){\makebox(0,0)[lb]{\smash{\SetFigFont{12}{14.4}{\rmdefault}{\mddefault}{\updefault}$a$}}}
\put(5525,-827){\makebox(0,0)[lb]{\smash{\SetFigFont{12}{14.4}{\rmdefault}{\mddefault}{\updefault}$q_2$}}}
\put(5554,-245){\makebox(0,0)[lb]{\smash{\SetFigFont{12}{14.4}{\rmdefault}{\mddefault}{\updefault}$q_1$}}}
\put(5302,-584){\makebox(0,0)[lb]{\smash{\SetFigFont{12}{14.4}{\rmdefault}{\mddefault}{\updefault}$p$}}}
\put(423,-479){\makebox(0,0)[lb]{\smash{\SetFigFont{12}{14.4}{\rmdefault}{\mddefault}{\updefault}$D=$}}}
\put(1531,-509){\makebox(0,0)[lb]{\smash{\SetFigFont{12}{14.4}{\rmdefault}{\mddefault}{\updefault}$p$}}}
\put(1536,233){\makebox(0,0)[lb]{\smash{\SetFigFont{12}{14.4}{\rmdefault}{\mddefault}{\updefault}$a$}}}
\put(2146,-544){\makebox(0,0)[lb]{\smash{\SetFigFont{12}{14.4}{\rmdefault}{\mddefault}{\updefault},}}}
\put(4114,-473){\makebox(0,0)[lb]{\smash{\SetFigFont{9}{10.8}{\rmdefault}{\mddefault}{\updefault},}}}
\put(2333,-534){\makebox(0,0)[lb]{\smash{\SetFigFont{12}{14.4}{\rmdefault}{\mddefault}{\updefault}$\mathcal{D}_1=$}}}
\put(3029,171){\makebox(0,0)[lb]{\smash{\SetFigFont{12}{14.4}{\rmdefault}{\mddefault}{\updefault}$a$}}}
\put(3630,-1161){\makebox(0,0)[lb]{\smash{\SetFigFont{12}{14.4}{\rmdefault}{\mddefault}{\updefault}$-a$}}}
\put(3171,-533){\makebox(0,0)[lb]{\smash{\SetFigFont{12}{14.4}{\rmdefault}{\mddefault}{\updefault}$p$}}}
\put(3667,128){\makebox(0,0)[lb]{\smash{\SetFigFont{12}{14.4}{\rmdefault}{\mddefault}{\updefault}$a$}}}
\put(3601,-768){\makebox(0,0)[lb]{\smash{\SetFigFont{12}{14.4}{\rmdefault}{\mddefault}{\updefault}$q_2$}}}
\put(3603,-267){\makebox(0,0)[lb]{\smash{\SetFigFont{12}{14.4}{\rmdefault}{\mddefault}{\updefault}$q_1$}}}
\put(4359,-550){\makebox(0,0)[lb]{\smash{\SetFigFont{12}{14.4}{\rmdefault}{\mddefault}{\updefault}$\mathcal{D}_2=$}}}
\end{picture}
$$
\begin{center}
{\bf Fig. 12}
\end{center}

\begin{displaymath}
c(D) = \sum_{\mathcal{D}_1}{w(q_1)w(q_2)} + \sum_{\mathcal{D}_2}{w(q_1)w(q_2)}
\end{displaymath}
is a class of $D$ of degree 2. Indeed, in any $G$-pure isotopy of a knot $K$
such that the crossing $p$ does not disappear, we observe the following: the
chord $p$ can get crossed with none of the $q_i, i=1, 2$. The chords $q_1$ and
$q_2$ can get crossed together by passing e.g. a stratum of the type
$a_{\dr}^+(a|-a,2a)$. But we count them now in $\mathcal{D}_2$ instead of
$\mathcal{D}_1$. Notice that the move depicted in Fig.~13 is again not
possible. 
$$
\begin{picture}(0,0)%
\includegraphics{ima13.pstex}%
\end{picture}%
\setlength{\unitlength}{4144sp}%
\begingroup\makeatletter\ifx\SetFigFont\undefined%
\gdef\SetFigFont#1#2#3#4#5{%
  \reset@font\fontsize{#1}{#2pt}%
  \fontfamily{#3}\fontseries{#4}\fontshape{#5}%
  \selectfont}%
\fi\endgroup%
\begin{picture}(3468,3319)(439,-2841)
\put(878,-1016){\makebox(0,0)[lb]{\smash{\SetFigFont{12}{14.4}{\rmdefault}{\mddefault}{\updefault}$p$}}}
\put(1200,-1589){\makebox(0,0)[lb]{\smash{\SetFigFont{12}{14.4}{\rmdefault}{\mddefault}{\updefault}$-a$}}}
\put(727,-209){\makebox(0,0)[lb]{\smash{\SetFigFont{12}{14.4}{\rmdefault}{\mddefault}{\updefault}$a$}}}
\put(1245,-209){\makebox(0,0)[lb]{\smash{\SetFigFont{12}{14.4}{\rmdefault}{\mddefault}{\updefault}$a$}}}
\put(3261,313){\makebox(0,0)[lb]{\smash{\SetFigFont{12}{14.4}{\rmdefault}{\mddefault}{\updefault}$a$}}}
\put(2754,-2144){\makebox(0,0)[lb]{\smash{\SetFigFont{12}{14.4}{\rmdefault}{\mddefault}{\updefault}$p$}}}
\put(2908,276){\makebox(0,0)[lb]{\smash{\SetFigFont{12}{14.4}{\rmdefault}{\mddefault}{\updefault}$a$}}}
\put(3533,-1034){\makebox(0,0)[lb]{\smash{\SetFigFont{12}{14.4}{\rmdefault}{\mddefault}{\updefault}$-a$}}}
\put(3226,-2841){\makebox(0,0)[lb]{\smash{\SetFigFont{12}{14.4}{\rmdefault}{\mddefault}{\updefault}$-a$}}}
\put(3606,-1499){\makebox(0,0)[lb]{\smash{\SetFigFont{12}{14.4}{\rmdefault}{\mddefault}{\updefault}$a$}}}
\put(2916,-1402){\makebox(0,0)[lb]{\smash{\SetFigFont{12}{14.4}{\rmdefault}{\mddefault}{\updefault}$a$}}}
\put(2072,-1711){\makebox(0,0)[lb]{\smash{\SetFigFont{12}{14.4}{\rmdefault}{\mddefault}{\updefault}or}}}
\put(3253,-247){\makebox(0,0)[lb]{\smash{\SetFigFont{12}{14.4}{\rmdefault}{\mddefault}{\updefault}$q_1$}}}
\put(3563,-484){\makebox(0,0)[lb]{\smash{\SetFigFont{12}{14.4}{\rmdefault}{\mddefault}{\updefault}$q_2$}}}
\put(2663,-286){\makebox(0,0)[lb]{\smash{\SetFigFont{12}{14.4}{\rmdefault}{\mddefault}{\updefault}$p$}}}
\put(1140,-572){\makebox(0,0)[lb]{\smash{\SetFigFont{12}{14.4}{\rmdefault}{\mddefault}{\updefault}$q_1$}}}
\put(1186,-1096){\makebox(0,0)[lb]{\smash{\SetFigFont{12}{14.4}{\rmdefault}{\mddefault}{\updefault}$q_2$}}}
\put(3330,-1789){\makebox(0,0)[lb]{\smash{\SetFigFont{12}{14.4}{\rmdefault}{\mddefault}{\updefault}$q_2$}}}
\put(3645,-2055){\makebox(0,0)[lb]{\smash{\SetFigFont{12}{14.4}{\rmdefault}{\mddefault}{\updefault}$q_1$}}}
\end{picture}
$$
\begin{center}
{\bf Fig. 13}
\end{center}
\par Let e.g. $(\mathcal{D}_1)_t$ be a subdiagram which represents
$\mathcal{D}_1$ for a knot $K_t$. If in the isotopy e.g. $(q_1)_t$ disappears
(see the left-hand part of Fig.~14), then there is a crossing $(q'_1)_t$ such
that the diagram depicted in the right-hand part of Fig.~14 represents also
$\mathcal{D}_1$. But for the writhes $w((q_1)_t)=-w((q'_1)_t)$, and hence for
$f_1=f_2=w(q_1)w(q_2)$, the contributions of $(q_1)_t$ and $(q'_1)_t$ in
$c(D_t)$ cancel out. This shows that $c(D)$ is a class, and calculating
examples one easily establishes that it is of degree 3.
$$
\begin{picture}(0,0)%
\includegraphics{ima14.pstex}%
\end{picture}%
\setlength{\unitlength}{4144sp}%
\begingroup\makeatletter\ifx\SetFigFont\undefined%
\gdef\SetFigFont#1#2#3#4#5{%
  \reset@font\fontsize{#1}{#2pt}%
  \fontfamily{#3}\fontseries{#4}\fontshape{#5}%
  \selectfont}%
\fi\endgroup%
\begin{picture}(5042,1497)(464,-1280)
\put(4510, 52){\makebox(0,0)[lb]{\smash{\SetFigFont{12}{14.4}{\rmdefault}{\mddefault}{\updefault}$a$}}}
\put(5111,-1280){\makebox(0,0)[lb]{\smash{\SetFigFont{12}{14.4}{\rmdefault}{\mddefault}{\updefault}$-a$}}}
\put(5148,  9){\makebox(0,0)[lb]{\smash{\SetFigFont{12}{14.4}{\rmdefault}{\mddefault}{\updefault}$a$}}}
\put(4352,-615){\makebox(0,0)[lb]{\smash{\SetFigFont{12}{14.4}{\rmdefault}{\mddefault}{\updefault}$p_t$}}}
\put(4777,-401){\makebox(0,0)[lb]{\smash{\SetFigFont{12}{14.4}{\rmdefault}{\mddefault}{\updefault}$(q_1')_t$}}}
\put(4835,-858){\makebox(0,0)[lb]{\smash{\SetFigFont{12}{14.4}{\rmdefault}{\mddefault}{\updefault}$(q_2)_t$}}}
\put(1192,-679){\makebox(0,0)[lb]{\smash{\SetFigFont{12}{14.4}{\rmdefault}{\mddefault}{\updefault}$(q_1')_t$}}}
\put(554,-671){\makebox(0,0)[lb]{\smash{\SetFigFont{12}{14.4}{\rmdefault}{\mddefault}{\updefault}$(q_1)_t$}}}
\put(2043,-346){\makebox(0,0)[lb]{\smash{\SetFigFont{12}{14.4}{\rmdefault}{\mddefault}{\updefault}$t$}}}
\end{picture}
$$
\begin{center}
{\bf Fig. 14}
\end{center}
\par Let $K \hookrightarrow F^2 \times \mathbb{R}$ be a $G$-pure global knot
and let $D$ be a $G$-pure configuration of degree $m$. Let $c_i(D), i \in \{1,
\dots, k\}$ be a finite collection of $G$-pure classes of $D$ having degree
$n_i$ respectively, and let $c_i, i \in \{1, \dots, k\}$ be fixed integers. Let $n = \max_{i}{n_i}$.

\begin{defin}
{\rm The {\it $T$-invariant\/} $T_K(D; c_1(D)=c_1, \dots c_k(D)=c_k)$ for $G$-pure global knots is defined as:
\begin{displaymath}
\sum_{D}{w(p_1) \cdots w(p_m)}
\end{displaymath}
Where we sum over all $D$ occuring as subdiagrams of the Gauss diagram of $K$,
such that $c_1(D)=c_1, \dots c_k(D)=c_k$. Here, $w(p_i)$ are the writhes of
the crossings of $K$ which correspond to the chords of the subdiagram
representing $D$. If $n=0$ and $m \not= 0$, then $T_K(D; \emptyset)$ is defined as
\begin{displaymath}
\sum_{D}{w(p_1) \cdots w(p_m)}
\end{displaymath}
Where we sum over all $D \subset$ Gauss diagram of $K$. If $m=0$ and $n \not=
0$ then $T_K(\emptyset; c(\emptyset))$ is defined as
\begin{displaymath}
c(\emptyset) = \sum_{i} \sum_{\Delta_i}{f_i(w(q^{(i)}_1,
\dots, q^{(i)}_{n_i}))}
\end{displaymath}
Where $\Delta_i$ runs through all subdiagrams which represent $\mathcal{D}_i$
in the Gauss diagram of $K$.If $n= m= 0$, i.e. $D= \mathcal{D}_i= \emptyset$,
then $T_K(\emptyset, \emptyset)$ is defined as the {\it free regular homotopy
  class\/} of $pr(K) \subset F^2$. (This is the universal invariant of degree
0.)}
\end{defin}

\begin{defin}
{\rm The set $\{ c_1(D), \dots, c_k(D) \}$ is called a {\it multi-class\/} of $D$.}
\end{defin}
{\bf Remark}.
\par If there is no risk of confusion, we will denote shortly by $T_K$ the
invariant $T_K(D; c_1(D)= c_1, \dots, c_k(D)= c_k)$.
\par We are now ready to formulate our main result.
\begin{theor}
{\it Let $K_0$ and $K_1$ be $G$-pure global knots which are ambient isotopic. 
Then, for each $T$-invariant for $G$-pure global knots}
\end{theor}
\begin{displaymath}
T_{K_0}(D; c_1(D)= c_1, \dots, c_k(D)= c_k)= T_{K_1}(D; c_1(D)= c_1, \dots, c_k(D)= c_k)
\end{displaymath}

{\bf Remark}.
\par The isotopy in the theorem need not to be $G$-pure!

{\bf Proof}.
The formal proof is very complicated. We just outline the main steps and let
the verification of the details for the reader. It follows from the definition
of a class $c(D)$ that $T_K$ is invariant for all such $G$-pure isotopies in
which no crossing $p_i$ disappears for a subdiagram (of the Gauss diagram of
$K$) representing $D$. Assume now that $D$ is represented by $\{ p_1, \dots,
p_i, \dots p_m \}$ and that the move depicted in Fig.~15 occurs in the isotopy:
$$
\begin{picture}(0,0)%
\includegraphics{ima15.pstex}%
\end{picture}%
\setlength{\unitlength}{4144sp}%
\begingroup\makeatletter\ifx\SetFigFont\undefined%
\gdef\SetFigFont#1#2#3#4#5{%
  \reset@font\fontsize{#1}{#2pt}%
  \fontfamily{#3}\fontseries{#4}\fontshape{#5}%
  \selectfont}%
\fi\endgroup%
\begin{picture}(4599,650)(441,-271)
\put(1668,-271){\makebox(0,0)[lb]{\smash{\SetFigFont{12}{14.4}{\rmdefault}{\mddefault}{\updefault}$K_t$}}}
\put(576,-136){\makebox(0,0)[lb]{\smash{\SetFigFont{12}{14.4}{\rmdefault}{\mddefault}{\updefault}$p_i$}}}
\put(1281,-136){\makebox(0,0)[rb]{\smash{\SetFigFont{12}{14.4}{\rmdefault}{\mddefault}{\updefault}$p_i'$}}}
\end{picture}
$$
\begin{center}
{\bf Fig. 15}
\end{center}
(Remember that \stratapm does not occur in a transversal isotopy.) Then the crossing $p'_i$, with which $p_i$ disappears, verifies:
\begin{enumerate}
\item
$w(p_i)= -w(p'_i)$
\item
$p'_i \notin \{p_1, \dots, p_i, \dots, p_m \}$. (This follows from 2. in Definition.~3.2.)
\item
$\{ p_1, \dots, p'_i, \dots, p_m \}$ represents $D$ ($p_i$ is replaced by $p'_i$).
\item
$c(\{ p_1, \dots, p_i, \dots, p_m \} = c(\{ p_1, \dots, p'_i, \dots p_m \}$, (in the right-hand term, $p_i$ is replaced by $p'_i$).
\end{enumerate}
Consequently, the contributions of $p_i$ and $p'_i$ in $T_{K_t}$
cancel out. This implies that $T_K$ is invariant for $G$-pure
isotopies.
\par Assume now that $K_t, t \in [0,1]$ is an isotopy which is not necessarily
$G$-pure (but $K_0$ and $K_1$ are $G$-pure!) The (value of) a class $c(D)$ can
change only if two of the crossings among $\{ p_1, \dots, p_m; q^{(i)}_1,
\dots, q^{(i)}_n \}$ of a configuration $\mathcal{D}_i$ are involved in a
Reidemeister move of type $III$, such that the third crossing involved has the
homological marking 0 or $\pm[K]_G$. Performing the isotopy $K_t$, let us
watch the traces on $F^2$ of the crossings with markings in $\{0, [K]_G,
-[K]_G \}$. These traces form immersed circles, called {\em Rudolph
  diagrams\/} (see [F], sect.~4.11). To each such circle, we associate a
family of disks $D_t$ (which are immersed in $F^2$), exactly as in the proof
of Theorem.~4.3 in [F]. This is possible because in the isotopy $K_t$, there
are no Reidemeister moves of type $II$ with opposite directions of the
tangencies. Such moves could destroy the disks $D_t$, as shown in Fig.~16.
$$
\begin{picture}(0,0)%
\includegraphics{im5.pstex}%
\end{picture}%
\setlength{\unitlength}{4144sp}%
\begingroup\makeatletter\ifx\SetFigFont\undefined%
\gdef\SetFigFont#1#2#3#4#5{%
  \reset@font\fontsize{#1}{#2pt}%
  \fontfamily{#3}\fontseries{#4}\fontshape{#5}%
  \selectfont}%
\fi\endgroup%
\begin{picture}(5707,3415)(386,-2645)
\put(5563,-2012){\makebox(0,0)[lb]{\smash{\SetFigFont{17}{20.4}{\rmdefault}{\mddefault}{\updefault}$?$}}}
\put(1651,-307){\makebox(0,0)[lb]{\smash{\SetFigFont{12}{14.4}{\rmdefault}{\mddefault}{\updefault}$K_t$}}}
\put(5552,317){\makebox(0,0)[lb]{\smash{\SetFigFont{12}{14.4}{\rmdefault}{\mddefault}{\updefault}$D_t$}}}
\put(5177,605){\makebox(0,0)[lb]{\smash{\SetFigFont{12}{14.4}{\rmdefault}{\mddefault}{\updefault}$0$}}}
\put(5267, 71){\makebox(0,0)[lb]{\smash{\SetFigFont{12}{14.4}{\rmdefault}{\mddefault}{\updefault}$0$}}}
\put(1131,-1477){\makebox(0,0)[lb]{\smash{\SetFigFont{12}{14.4}{\rmdefault}{\mddefault}{\updefault}$0$}}}
\put(1590,-1855){\makebox(0,0)[lb]{\smash{\SetFigFont{12}{14.4}{\rmdefault}{\mddefault}{\updefault}$D_t$}}}
\put(771,-2335){\makebox(0,0)[lb]{\smash{\SetFigFont{12}{14.4}{\rmdefault}{\mddefault}{\updefault}$0$}}}
\put(3787,-1944){\makebox(0,0)[lb]{\smash{\SetFigFont{12}{14.4}{\rmdefault}{\mddefault}{\updefault}$D_t$}}}
\put(2932,-2460){\makebox(0,0)[lb]{\smash{\SetFigFont{12}{14.4}{\rmdefault}{\mddefault}{\updefault}$0$}}}
\put(4814,-2496){\makebox(0,0)[lb]{\smash{\SetFigFont{12}{14.4}{\rmdefault}{\mddefault}{\updefault}$0$}}}
\end{picture}
$$
\begin{center}
{\bf Fig. 16}
\end{center}
When Reidemeister moves of type $II$, with equal directions of the tangencies,
occur in the isotopy $K_t$, one gets the usual surgeries of the disks, as shown in Fig.~17 and 18.
$$
\begin{picture}(0,0)%
\includegraphics{im6.pstex}%
\end{picture}%
\setlength{\unitlength}{4144sp}%
\begingroup\makeatletter\ifx\SetFigFont\undefined%
\gdef\SetFigFont#1#2#3#4#5{%
  \reset@font\fontsize{#1}{#2pt}%
  \fontfamily{#3}\fontseries{#4}\fontshape{#5}%
  \selectfont}%
\fi\endgroup%
\begin{picture}(4605,2831)(43,-2111)
\put(1889,174){\makebox(0,0)[lb]{\smash{\SetFigFont{12}{14.4}{\rmdefault}{\mddefault}{\updefault}$K_t$}}}
\put(3600,-521){\makebox(0,0)[lb]{\smash{\SetFigFont{12}{14.4}{\rmdefault}{\mddefault}{\updefault}$D_t$}}}
\put(2855,140){\makebox(0,0)[lb]{\smash{\SetFigFont{12}{14.4}{\rmdefault}{\mddefault}{\updefault}$0$}}}
\put(4415,118){\makebox(0,0)[lb]{\smash{\SetFigFont{12}{14.4}{\rmdefault}{\mddefault}{\updefault}$0$}}}
\put(1080,-1580){\makebox(0,0)[lb]{\smash{\SetFigFont{12}{14.4}{\rmdefault}{\mddefault}{\updefault}$0$}}}
\put(2775,-1555){\makebox(0,0)[lb]{\smash{\SetFigFont{12}{14.4}{\rmdefault}{\mddefault}{\updefault}$0$}}}
\put(1870,-2111){\makebox(0,0)[lb]{\smash{\SetFigFont{12}{14.4}{\rmdefault}{\mddefault}{\updefault}$D_t$}}}
\put(475,470){\makebox(0,0)[lb]{\smash{\SetFigFont{12}{14.4}{\rmdefault}{\mddefault}{\updefault}$0$}}}
\put(160,468){\makebox(0,0)[lb]{\smash{\SetFigFont{12}{14.4}{\rmdefault}{\mddefault}{\updefault}$0$}}}
\put(1169,495){\makebox(0,0)[lb]{\smash{\SetFigFont{12}{14.4}{\rmdefault}{\mddefault}{\updefault}$0$}}}
\put(1524,475){\makebox(0,0)[lb]{\smash{\SetFigFont{12}{14.4}{\rmdefault}{\mddefault}{\updefault}$0$}}}
\put(769,-416){\makebox(0,0)[lb]{\smash{\SetFigFont{12}{14.4}{\rmdefault}{\mddefault}{\updefault}$D_t$}}}
\put(135,151){\makebox(0,0)[lb]{\smash{\SetFigFont{12}{14.4}{\rmdefault}{\mddefault}{\updefault}(or}}}
\put( 43,-83){\makebox(0,0)[lb]{\smash{\SetFigFont{12}{14.4}{\rmdefault}{\mddefault}{\updefault}$\pm [K]_G$)}}}
\end{picture}
$$
\begin{center}
{\bf Fig. 17}
\end{center}

$$
\begin{picture}(0,0)%
\includegraphics{im7.pstex}%
\end{picture}%
\setlength{\unitlength}{4144sp}%
\begingroup\makeatletter\ifx\SetFigFont\undefined%
\gdef\SetFigFont#1#2#3#4#5{%
  \reset@font\fontsize{#1}{#2pt}%
  \fontfamily{#3}\fontseries{#4}\fontshape{#5}%
  \selectfont}%
\fi\endgroup%
\begin{picture}(3826,1277)(453,-861)
\put(1331,-236){\makebox(0,0)[lb]{\smash{\SetFigFont{12}{14.4}{\rmdefault}{\mddefault}{\updefault}$K_t$}}}
\put(3526,-321){\makebox(0,0)[lb]{\smash{\SetFigFont{12}{14.4}{\rmdefault}{\mddefault}{\updefault}$0$}}}
\put(4051,-311){\makebox(0,0)[lb]{\smash{\SetFigFont{12}{14.4}{\rmdefault}{\mddefault}{\updefault}$0$}}}
\put(3111,-861){\makebox(0,0)[lb]{\smash{\SetFigFont{12}{14.4}{\rmdefault}{\mddefault}{\updefault}$D_t$}}}
\end{picture}
$$
\begin{center}
{\bf Fig. 18}
\end{center}
Remember that in a transversal isotopy, there are no Reidemeister moves of
type $I$. Each of the families of disks $D_t$ starts and ends with the empty
set because the phenomen explained in Fig.~101 in [F] can still not
appear. Each disk in $D_t$ has exactly two vertices corresponding both to
crossings of type $0$ or $[K]_G$ or $-[K]_G$ which appeared in the underlying circle of Rudolph's diagram.
\par We distinguish now two cases:
\begin{itemize}
\item
{\em The simple case:\/} Let $r_1, r_2 \in \{ p_1, \dots, p_m; q^{(i)}_1,
\dots, q^{(i)}_n \}$ be two crossings of a configuration $\mathcal{D}_i$ used
in the definition of $T_K$. We assume that $r_1$ and $r_2$ pass together for
the first time $t$ in the isotopy $K_t$ over a vertex of a disk in $D_t$, and
neither $r_1$ nor $r_2$ passes over a vertex of {\em another\/} disk in $D_t$.
\par Then, similar arguments to those used in the proof of Theorem.~4.3 in [F]
show that one of the following three possibilities is realized:
\begin{enumerate}
\item
$r_1$ and $r_2$ will pass a second time together over the same vertex, but in opposite directions
\item
$r_1$ and $r_2$ will pass together also over the second vertex of the disk
\item
There exists a third crossing $r'_1$ which has appeared together with $r_1$ by
a Reidemeister move of type $II$ (hence, $w(r_1)= -w(r'_1)$, and $r'_1$ can
replace $r_1$ in the configuration $\mathcal{D}_i$). Moreover, $r'_1$ will
pass together with $r_2$ over a vertex of the disk (compare with the proof of Theorem.~4.3 in [F]). (See Fig.~19.)
\end{enumerate}
$$
\begin{picture}(0,0)%
\includegraphics{im8.pstex}%
\end{picture}%
\setlength{\unitlength}{4144sp}%
\begingroup\makeatletter\ifx\SetFigFont\undefined%
\gdef\SetFigFont#1#2#3#4#5{%
  \reset@font\fontsize{#1}{#2pt}%
  \fontfamily{#3}\fontseries{#4}\fontshape{#5}%
  \selectfont}%
\fi\endgroup%
\begin{picture}(5186,4907)(324,-4596)
\put(885, 34){\makebox(0,0)[lb]{\smash{\SetFigFont{12}{14.4}{\rmdefault}{\mddefault}{\updefault}$r_2$}}}
\put(1055,-564){\makebox(0,0)[lb]{\smash{\SetFigFont{12}{14.4}{\rmdefault}{\mddefault}{\updefault}$0$}}}
\put(1925,-528){\makebox(0,0)[lb]{\smash{\SetFigFont{12}{14.4}{\rmdefault}{\mddefault}{\updefault}$0$}}}
\put(1630,-993){\makebox(0,0)[lb]{\smash{\SetFigFont{12}{14.4}{\rmdefault}{\mddefault}{\updefault}$D_t$}}}
\put(3417,-3749){\makebox(0,0)[lb]{\smash{\SetFigFont{12}{14.4}{\rmdefault}{\mddefault}{\updefault}$0$}}}
\put(4599,-3762){\makebox(0,0)[lb]{\smash{\SetFigFont{12}{14.4}{\rmdefault}{\mddefault}{\updefault}$0$}}}
\put(3362,-819){\makebox(0,0)[lb]{\smash{\SetFigFont{12}{14.4}{\rmdefault}{\mddefault}{\updefault}$0$}}}
\put(4507,-764){\makebox(0,0)[lb]{\smash{\SetFigFont{12}{14.4}{\rmdefault}{\mddefault}{\updefault}$0$}}}
\put(3676,-179){\makebox(0,0)[lb]{\smash{\SetFigFont{12}{14.4}{\rmdefault}{\mddefault}{\updefault}$r_1$}}}
\put(4952,-3884){\makebox(0,0)[lb]{\smash{\SetFigFont{12}{14.4}{\rmdefault}{\mddefault}{\updefault}$r_2$}}}
\put(4935,-4178){\makebox(0,0)[lb]{\smash{\SetFigFont{12}{14.4}{\rmdefault}{\mddefault}{\updefault}$r_1'$}}}
\put(3242,-29){\makebox(0,0)[lb]{\smash{\SetFigFont{12}{14.4}{\rmdefault}{\mddefault}{\updefault}$r_2$}}}
\put(4095,-83){\makebox(0,0)[lb]{\smash{\SetFigFont{12}{14.4}{\rmdefault}{\mddefault}{\updefault}$r_1'$}}}
\put(3376,-4469){\makebox(0,0)[lb]{\smash{\SetFigFont{12}{14.4}{\rmdefault}{\mddefault}{\updefault}$r_1$}}}
\put(627,-2384){\makebox(0,0)[lb]{\smash{\SetFigFont{12}{14.4}{\rmdefault}{\mddefault}{\updefault}$0$}}}
\put(1455,-2093){\makebox(0,0)[lb]{\smash{\SetFigFont{12}{14.4}{\rmdefault}{\mddefault}{\updefault}$r_1'$}}}
\put(601,-2939){\makebox(0,0)[lb]{\smash{\SetFigFont{12}{14.4}{\rmdefault}{\mddefault}{\updefault}$r_1$}}}
\put(962,-3014){\makebox(0,0)[lb]{\smash{\SetFigFont{12}{14.4}{\rmdefault}{\mddefault}{\updefault}$r_2$}}}
\put(3512,-1454){\makebox(0,0)[lb]{\smash{\SetFigFont{12}{14.4}{\rmdefault}{\mddefault}{\updefault}$r_2$}}}
\put(3552,-1769){\makebox(0,0)[lb]{\smash{\SetFigFont{12}{14.4}{\rmdefault}{\mddefault}{\updefault}$0$}}}
\put(4802,-1701){\makebox(0,0)[lb]{\smash{\SetFigFont{12}{14.4}{\rmdefault}{\mddefault}{\updefault}$0$}}}
\put(4522,-2616){\makebox(0,0)[lb]{\smash{\SetFigFont{12}{14.4}{\rmdefault}{\mddefault}{\updefault}$D_t$}}}
\put(4672,-4596){\makebox(0,0)[lb]{\smash{\SetFigFont{12}{14.4}{\rmdefault}{\mddefault}{\updefault}$D_t$}}}
\put(1912,-2799){\makebox(0,0)[lb]{\smash{\SetFigFont{12}{14.4}{\rmdefault}{\mddefault}{\updefault}$0$}}}
\put(2713,-2811){\makebox(0,0)[lb]{\smash{\SetFigFont{12}{14.4}{\rmdefault}{\mddefault}{\updefault}or}}}
\put(2717,-2552){\makebox(0,0)[lb]{\smash{\SetFigFont{12}{14.4}{\rmdefault}{\mddefault}{\updefault}either}}}
\end{picture}
$$
\begin{center}
{\bf Fig. 19}
\end{center}
\par Obviously, $T_{K_t}$ does not change for the first two possibilities. If
in the third possibility, $r_1$ and $r_2$ enter together into some
configuration $\mathcal{D}_i$, then $r_2$ and $r'_1$ (instead of $r_1$) enter
also in this configuration $\mathcal{D}_i$. There are two cases to distinguish:
\begin{itemize}
\item
A) $r_1$ and hence $r'_1$ belong to $\{ p_1, \dots, p_m \} \subset
\mathcal{D}_i$. But then, the contributions of $r_1$ and $r'_1$ cancel out in $T_{K_t}$ because the weight function is $w(p_1) \cdots w(p_m)$.
\item
B) $r_1$ (and hence $r'_1$) belongs to $\{ q^{(i)}_1, \dots, q^{(i)}_n \}
\subset \mathcal{D}_i$. But then the definition of a class $c(D)$ implies that
the contributions of $r_1$ and $r'_1$ cancel out in $c(D)$, and hence $T_{K_t}$ is again invariant.
\end{itemize}
\par {\it The general case:\/} During the isotopy $K_t$, $r_1$ and $r_2$ pass
together over a vertex of a disk in $D_t$. Let $t_0$ be the smallest value of
$t$ for which this occurs. At a time $t_1> t_0$, $r_1$ passes over the vertex
of another disk in $D_t$. Notice, that $r_1$ and $r_2$ cannot pass together
over the vertex of another disk in $D_t$. We illustrate the general case with
an example, shown in Fig.~20.
$$
\begin{picture}(0,0)%
\includegraphics{im9.pstex}%
\end{picture}%
\setlength{\unitlength}{4144sp}%
\begingroup\makeatletter\ifx\SetFigFont\undefined%
\gdef\SetFigFont#1#2#3#4#5{%
  \reset@font\fontsize{#1}{#2pt}%
  \fontfamily{#3}\fontseries{#4}\fontshape{#5}%
  \selectfont}%
\fi\endgroup%
\begin{picture}(4487,4329)(446,-3872)
\put(1358,262){\makebox(0,0)[lb]{\smash{\SetFigFont{12}{14.4}{\rmdefault}{\mddefault}{\updefault}$r_3$}}}
\put(766,-286){\makebox(0,0)[lb]{\smash{\SetFigFont{12}{14.4}{\rmdefault}{\mddefault}{\updefault}$r_2$}}}
\put(833,-1073){\makebox(0,0)[lb]{\smash{\SetFigFont{12}{14.4}{\rmdefault}{\mddefault}{\updefault}$r_1$}}}
\put(2416,-76){\makebox(0,0)[lb]{\smash{\SetFigFont{12}{14.4}{\rmdefault}{\mddefault}{\updefault}$[K]_G$}}}
\put(1276,-811){\makebox(0,0)[lb]{\smash{\SetFigFont{12}{14.4}{\rmdefault}{\mddefault}{\updefault}$0$}}}
\put(2813,-346){\makebox(0,0)[lb]{\smash{\SetFigFont{12}{14.4}{\rmdefault}{\mddefault}{\updefault}$0$}}}
\put(3496,-1253){\makebox(0,0)[lb]{\smash{\SetFigFont{12}{14.4}{\rmdefault}{\mddefault}{\updefault}$[K]_G$}}}
\put(2318,-1201){\makebox(0,0)[lb]{\smash{\SetFigFont{12}{14.4}{\rmdefault}{\mddefault}{\updefault}$S$}}}
\put(3421,-1988){\makebox(0,0)[lb]{\smash{\SetFigFont{12}{14.4}{\rmdefault}{\mddefault}{\updefault}$r_2$}}}
\put(3983,-2768){\makebox(0,0)[lb]{\smash{\SetFigFont{12}{14.4}{\rmdefault}{\mddefault}{\updefault}$r_3$}}}
\put(3916,-3458){\makebox(0,0)[lb]{\smash{\SetFigFont{12}{14.4}{\rmdefault}{\mddefault}{\updefault}$r_1$}}}
\put(3091,-2805){\makebox(0,0)[lb]{\smash{\SetFigFont{12}{14.4}{\rmdefault}{\mddefault}{\updefault}$[K]_G$}}}
\put(2176,-2355){\makebox(0,0)[lb]{\smash{\SetFigFont{12}{14.4}{\rmdefault}{\mddefault}{\updefault}$[K]_G$}}}
\put(2685,-2551){\makebox(0,0)[lb]{\smash{\SetFigFont{12}{14.4}{\rmdefault}{\mddefault}{\updefault}$0$}}}
\put(1171,-3031){\makebox(0,0)[lb]{\smash{\SetFigFont{12}{14.4}{\rmdefault}{\mddefault}{\updefault}$0$}}}
\put(2221,-3458){\makebox(0,0)[lb]{\smash{\SetFigFont{12}{14.4}{\rmdefault}{\mddefault}{\updefault}$S$}}}
\put(4594,-811){\makebox(0,0)[lb]{\smash{\SetFigFont{12}{14.4}{\rmdefault}{\mddefault}{\updefault}$K_t$}}}
\put(1426,-1576){\makebox(0,0)[lb]{\smash{\SetFigFont{12}{14.4}{\rmdefault}{\mddefault}{\updefault}$(D_1)_t$}}}
\put(2948,-1793){\makebox(0,0)[lb]{\smash{\SetFigFont{12}{14.4}{\rmdefault}{\mddefault}{\updefault}$(D_2)_t$}}}
\end{picture}
$$
\begin{center}
{\bf Fig. 20}
\end{center}
We can assume that the crossing $s$ in Fig.~20 is not of type $0$ or $\pm
[K]_G$ because otherwise, the crossings $r_2$ and $r_3$ on the left-hand side
would already have passed together over the vertex of a disk in $D_t$.
\par We observe that $r_1$ with $r_2$ on the left-hand side of Fig.~20 form
exactly the same configurations as $r_1$ with $r_3$ on the right-hand side of
Fig.~20. But the mutual configuration of $r_2$ and $r_3$ has changed (by passing $s$).
\par Notice that we cannot eliminate the disks $(D_1)_t$, $(D_2)_t$ by an
isotopy, which is supported in a 3-ball containing just the fragment of $K_t$
drawn in Fig.~20. Nevertheless, we can replace the local piece of the isotopy
$K_t$ shown in Fig.~20 by the local piece of a $G$-pure isotopy $K'_t$ shown
in Fig.~21. 
$$
\begin{picture}(0,0)%
\includegraphics{im10.pstex}%
\end{picture}%
\setlength{\unitlength}{4144sp}%
\begingroup\makeatletter\ifx\SetFigFont\undefined%
\gdef\SetFigFont#1#2#3#4#5{%
  \reset@font\fontsize{#1}{#2pt}%
  \fontfamily{#3}\fontseries{#4}\fontshape{#5}%
  \selectfont}%
\fi\endgroup%
\begin{picture}(4727,4555)(206,-4098)
\put(4594,-811){\makebox(0,0)[lb]{\smash{\SetFigFont{12}{14.4}{\rmdefault}{\mddefault}{\updefault}$K_t'$}}}
\put(2243,-3234){\makebox(0,0)[lb]{\smash{\SetFigFont{12}{14.4}{\rmdefault}{\mddefault}{\updefault}$S$}}}
\put(3218,-2716){\makebox(0,0)[lb]{\smash{\SetFigFont{12}{14.4}{\rmdefault}{\mddefault}{\updefault}$r_2$}}}
\put(2782,-3834){\makebox(0,0)[lb]{\smash{\SetFigFont{12}{14.4}{\rmdefault}{\mddefault}{\updefault}$r_1$}}}
\put(2986,-3317){\makebox(0,0)[lb]{\smash{\SetFigFont{12}{14.4}{\rmdefault}{\mddefault}{\updefault}$r_3$}}}
\put(1358,262){\makebox(0,0)[lb]{\smash{\SetFigFont{12}{14.4}{\rmdefault}{\mddefault}{\updefault}$r_3$}}}
\put(2198,-841){\makebox(0,0)[lb]{\smash{\SetFigFont{12}{14.4}{\rmdefault}{\mddefault}{\updefault}$S$}}}
\put(765,-1118){\makebox(0,0)[lb]{\smash{\SetFigFont{12}{14.4}{\rmdefault}{\mddefault}{\updefault}$r_1$}}}
\put(893,-841){\makebox(0,0)[lb]{\smash{\SetFigFont{12}{14.4}{\rmdefault}{\mddefault}{\updefault}$r_2$}}}
\end{picture}
$$
\begin{center}
{\bf Fig. 21}
\end{center}
The crossings $r_1, r_2, r_3$ contribute to $T_{K_t}$ on the
left-hand-side (resp. right-hand side) of Fig.~20 exactly the same way as they
contribute to $T_{K'_t}$ in the left-hand side (resp. right-hand side) of
Fig.~21. But the isotopy in Fig.~21 is $G$-pure and this implies, as already
proven, that $T_{K'_t}$ is invariant. Consequently, $T_{K_t}$ for the isotopy
$K_t$ in fig.~20 is invariant too.
\par Clearly, these arguments can be generalized for the case of more than two
disks $D_t$ and several crossings instead of only $s$ in the local piece of the isotopy $K_t$.
$\Box$
\end{itemize}
\begin{defin}
{\rm Let $n= \max_i{n_i}$. The {\it degree of $T_K$ as Gauss diagram invariant\/} is equal to } $m+ n$.
\end{defin}
This definition is justified by the observation that the complexity of the
calculation of $T_K(D; c_1(D)= c_1, \dots, c_k(D)= c_k)$ for knots $K$ is a
polynomial of degree $m+ n$ in the number of crossings of $K$. Hence, the
invariant $T_K$ is calculated with the same (order of) complexity as a Vassiliev invariant of degree $m+n$.
\par {\bf Remarks}:
\begin{enumerate}
\item
Examples 1 and 3 show that the numbers of different $G$-pure configurations
$D$ of degree $m$, and of classes $c_i(D)$ of degree $n$ is in general not
finite (the configurations depend on the parameter $a \in \mathbb{Z} \setminus 0$).
\item
We show in an example in sect.~9 that multi-classes are usefull. As a matter
of fact, $T_K(D; c_1(D)= c_1, c_2(D)= c_2)$ contains sometimes more
information than $T_K(D; c_1(D)= c_1)$ and $T_K(D; c_2(D)= c_2)$ together.
\end{enumerate}
\begin{lem}
{\it If $m= 0$ or $n= 0$, then the invariant $T_K$ is of finite type in the sense of Vassiliev-Gussarov.}
\end{lem}
We omit the proof of this lemma: it is a straightforward generalization of
Oestlund's proof that the Gauss diagram invariants of Polyak-Viro for knots in
$\mathbb{R}^3$ are of finite type ([P-V], see also [F], sect.~2). 
\par Gussarov has proven that each Vassiliev invariant for knots in 
$\mathbb{R}^3$ can be represented as a Gauss diagram invariant of 
Polyak-Viro (see [G-P-V]). We do not know wether or not this is still true 
for finite type invariants of global knots. But evidently, each Gauss diagram 
invariant of finite type which does not use the homological markings 
$\{0, \pm [K]_G \}$ in $G$ is a $T$-invariant for $m= 0$ (i.e. $D= \emptyset$).
\begin{lem}
{\it Let $m> 0$ and $n>0$ and assume that $c(D)$ (from Definition~3.3) is not
  a Gauss diagram identity, i.e. there exists a knot $K_t$ and a crossing $q$
  of $K_t$ such that switching the crossing $q$ changes the value $c(D_t)$
  (see [F], sect.~4). Then the invariant $T_K$ is not of finite type.\/}
\end{lem}
We omit the general proof but show this in an example for a class of degree
1. It is clear that one can find similar examples for any class of degree at least 1.
{\it Example $3.4$}
We fix the system of generators $\{ \alpha, \beta \}$ of $H_1(T^2;
\mathbb{Z})$ as shown in Fig.~22. 
$$
\begin{picture}(0,0)%
\includegraphics{im11.pstex}%
\end{picture}%
\setlength{\unitlength}{4144sp}%
\begingroup\makeatletter\ifx\SetFigFont\undefined%
\gdef\SetFigFont#1#2#3#4#5{%
  \reset@font\fontsize{#1}{#2pt}%
  \fontfamily{#3}\fontseries{#4}\fontshape{#5}%
  \selectfont}%
\fi\endgroup%
\begin{picture}(4005,2040)(961,-2229)
\put(4966,-1059){\makebox(0,0)[lb]{\smash{\SetFigFont{12}{14.4}{\rmdefault}{\mddefault}{\updefault}$T^2$}}}
\put(4081,-1666){\makebox(0,0)[lb]{\smash{\SetFigFont{12}{14.4}{\rmdefault}{\mddefault}{\updefault}$\alpha$}}}
\put(3743,-2109){\makebox(0,0)[lb]{\smash{\SetFigFont{12}{14.4}{\rmdefault}{\mddefault}{\updefault}$\beta$}}}
\end{picture}
$$
\begin{center}
{\bf Fig. 22}
\end{center}
Let $f: T^2 \to S^1$ be a submersion such that the fibers represent
$\beta$. The vector field $v$ is defined as the unit tangent vector field to
the fibers of $f$. Let us consider the family of knots $K_n, n \in \mathbb{N}$
shown in Fig.~23.
$$
\begin{picture}(0,0)%
\includegraphics{im12.pstex}%
\end{picture}%
\setlength{\unitlength}{4144sp}%
\begingroup\makeatletter\ifx\SetFigFont\undefined%
\gdef\SetFigFont#1#2#3#4#5{%
  \reset@font\fontsize{#1}{#2pt}%
  \fontfamily{#3}\fontseries{#4}\fontshape{#5}%
  \selectfont}%
\fi\endgroup%
\begin{picture}(4942,2711)(476,-2327)
\put(1486,-1600){\makebox(0,0)[lb]{\smash{\SetFigFont{12}{14.4}{\rmdefault}{\mddefault}{\updefault}$+$}}}
\put(1876,-1637){\makebox(0,0)[lb]{\smash{\SetFigFont{12}{14.4}{\rmdefault}{\mddefault}{\updefault}$+$}}}
\put(1448,-1885){\makebox(0,0)[lb]{\smash{\SetFigFont{12}{14.4}{\rmdefault}{\mddefault}{\updefault}$p_1$}}}
\put(1928,-1975){\makebox(0,0)[lb]{\smash{\SetFigFont{12}{14.4}{\rmdefault}{\mddefault}{\updefault}$p_2$}}}
\put(5093,-1098){\makebox(0,0)[lb]{\smash{\SetFigFont{12}{14.4}{\rmdefault}{\mddefault}{\updefault}$K_n$}}}
\put(2558,-1143){\makebox(0,0)[lb]{\smash{\SetFigFont{12}{14.4}{\rmdefault}{\mddefault}{\updefault}$2n$}}}
\put(2138,-1720){\makebox(0,0)[lb]{\smash{\SetFigFont{12}{14.4}{\rmdefault}{\mddefault}{\updefault}$-$}}}
\put(2326,-1712){\makebox(0,0)[lb]{\smash{\SetFigFont{12}{14.4}{\rmdefault}{\mddefault}{\updefault}$-$}}}
\put(2963,-1607){\makebox(0,0)[lb]{\smash{\SetFigFont{12}{14.4}{\rmdefault}{\mddefault}{\updefault}$-$}}}
\end{picture}
$$
\begin{center}
{\bf Fig. 23}
\end{center}
$K_n \hookrightarrow T^2 \times \mathbb{R}$ is a global knot with repect to
$v$ and $[K_n]= 3\alpha+ \beta \in H_1(T^2)$. We take as group $G$:
\begin{displaymath}
(H_1(T^2)/\langle [K_n] \rangle) \otimes \mathbb{Z}/2\mathbb{Z}
\cong \mathbb{Z}/2\mathbb{Z}
\end{displaymath}
Each knot $K_n$ is $G$-pure, i.e. for each crossing $p$ the loop $pr(K^+_p)$
(as well as the loop $pr(K^-_p)$) represents the non-trivial element in $G
\cong \mathbb{Z}/2\mathbb{Z}$. (Thus, each chord is marked by the same element
in $G$ and we do not write the marking.)

\par Let us consider the unique $\mathbb{Z}/2\mathbb{Z}$-pure configuration
$D$ of degree 1:
$$
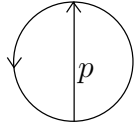
\begin{picture}(0,0)%
\includegraphics{pag46.pstex}%
\end{picture}%
\setlength{\unitlength}{4144sp}%
\begingroup\makeatletter\ifx\SetFigFont\undefined%
\gdef\SetFigFont#1#2#3#4#5{%
  \reset@font\fontsize{#1}{#2pt}%
  \fontfamily{#3}\fontseries{#4}\fontshape{#5}%
  \selectfont}%
\fi\endgroup%
\begin{picture}(772,744)(494,-258)
\put(931, 39){\makebox(0,0)[lb]{\smash{\SetFigFont{12}{14.4}{\rmdefault}{\mddefault}{\updefault}$p$}}}
\end{picture}
$$

Let us consider the class of $D$ of degree 1 defined by

\begin{displaymath}
c(D)= \sum_{\begin{picture}(0,0)%
\includegraphics{page47.pstex}%
\end{picture}%
\setlength{\unitlength}{4144sp}%
\begingroup\makeatletter\ifx\SetFigFont\undefined%
\gdef\SetFigFont#1#2#3#4#5{%
  \reset@font\fontsize{#1}{#2pt}%
  \fontfamily{#3}\fontseries{#4}\fontshape{#5}%
  \selectfont}%
\fi\endgroup%
\begin{picture}(631,578)(564,-175)
\put(1021, 59){\makebox(0,0)[lb]{\smash{\SetFigFont{12}{14.4}{\rmdefault}{\mddefault}{\updefault}$q$}}}
\put(841, 64){\makebox(0,0)[lb]{\smash{\SetFigFont{12}{14.4}{\rmdefault}{\mddefault}{\updefault}$p$}}}
\end{picture}
}{w(q)}
\end{displaymath}

(One easily verifies that $c(D)$ is indeed a class of $D$, because

the chords $p$ and $q$ cannot get crossed in a $\mathbb{Z}/2 \mathbb{Z}$-pure isotopy, and if a new couple of crossings $q_i, i= 1, 2$

appears by a Reidemeister move of type $II$, then $c(D_t)$ remains

unchanged.)

\par $K_n$ has the Gauss diagram depicted in Fig.~24.
$$
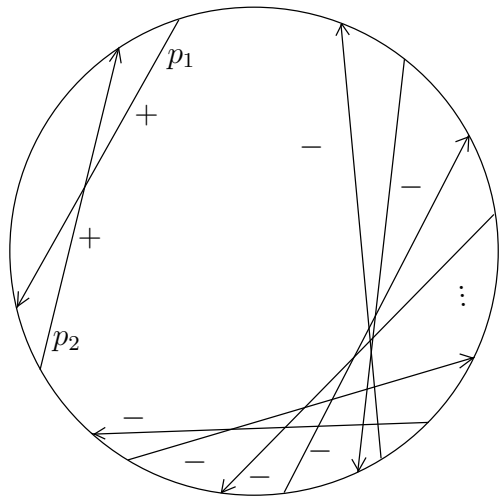
\begin{picture}(0,0)%
\includegraphics{ima24.pstex}%
\end{picture}%
\setlength{\unitlength}{4144sp}%
\begingroup\makeatletter\ifx\SetFigFont\undefined%
\gdef\SetFigFont#1#2#3#4#5{%
  \reset@font\fontsize{#1}{#2pt}%
  \fontfamily{#3}\fontseries{#4}\fontshape{#5}%
  \selectfont}%
\fi\endgroup%
\begin{picture}(2938,2938)(697,-2120)
\put(1661,489){\makebox(0,0)[lb]{\smash{\SetFigFont{12}{14.4}{\rmdefault}{\mddefault}{\updefault}$p_1$}}}
\put(1466,129){\makebox(0,0)[lb]{\smash{\SetFigFont{12}{14.4}{\rmdefault}{\mddefault}{\updefault}$+$}}}
\put(1126,-611){\makebox(0,0)[lb]{\smash{\SetFigFont{12}{14.4}{\rmdefault}{\mddefault}{\updefault}$+$}}}
\put(976,-1201){\makebox(0,0)[lb]{\smash{\SetFigFont{12}{14.4}{\rmdefault}{\mddefault}{\updefault}$p_2$}}}
\put(2451,-66){\makebox(0,0)[lb]{\smash{\SetFigFont{12}{14.4}{\rmdefault}{\mddefault}{\updefault}$-$}}}
\put(3046,-306){\makebox(0,0)[lb]{\smash{\SetFigFont{12}{14.4}{\rmdefault}{\mddefault}{\updefault}$-$}}}
\put(2501,-1886){\makebox(0,0)[lb]{\smash{\SetFigFont{12}{14.4}{\rmdefault}{\mddefault}{\updefault}$-$}}}
\put(2141,-2041){\makebox(0,0)[lb]{\smash{\SetFigFont{12}{14.4}{\rmdefault}{\mddefault}{\updefault}$-$}}}
\put(1751,-1951){\makebox(0,0)[lb]{\smash{\SetFigFont{12}{14.4}{\rmdefault}{\mddefault}{\updefault}$-$}}}
\put(1386,-1686){\makebox(0,0)[lb]{\smash{\SetFigFont{12}{14.4}{\rmdefault}{\mddefault}{\updefault}$-$}}}
\put(3396,-866){\makebox(0,0)[lb]{\smash{\SetFigFont{12}{14.4}{\rmdefault}{\mddefault}{\updefault}.}}}
\put(3401,-916){\makebox(0,0)[lb]{\smash{\SetFigFont{12}{14.4}{\rmdefault}{\mddefault}{\updefault}.}}}
\put(3406,-971){\makebox(0,0)[lb]{\smash{\SetFigFont{12}{14.4}{\rmdefault}{\mddefault}{\updefault}.}}}
\end{picture}
$$
\begin{center}
{\bf Fig. 24}
\end{center}

Consequently, $c(p_1)= -n$, $c(p_2)= 0$ and $c(p_i)= 0$ or $1$ for

each crossing $p_i$ of $K_n \setminus \{p_1, p_2 \}$. Therefore,

$T_{K_n}(D; c(D)=-n)= +1$ and $T_{K_n}(D; c(D)=r)= 0$  

for all $r \notin \{ -n, 0, 1 \}$.

Let $K'_n$ be any knot obtained from $K_n$ by changing any of the

$2n$ crossings of $K_n \setminus \{ p_1, p_2 \}$. Then $c(p_1)$

becomes strictly bigger than $-n$ and thus, $T_{K'_n}(D; c(D)=-n)= 0$

This shows that no linear combination with coefficients $\pm 1$ of

$T_{K_n}(D; c(D)=-n)$ with $T_{K'_n}(D; c(D)=-n)$ can ever be equal to

$0$. Moreover, replacing the fragment in $K_n$ depicted in the left-hand

side of Fig.~25 by the one depicted in the right-hand side of Fig.~25

does not change the knot $K_n$. Switching any of the $s$ crossings

$\{ q_1, \dots, q_s \}$ makes $c(p_1)$ strictly bigger than $-n$.

Therefore, $T_K(D; c(D)=-n)$ is of degree at least $2n+ s+ 1$ and hence

it is not of finite type. $\Box$

$$
\begin{picture}(0,0)%
\includegraphics{ima25.pstex}%
\end{picture}%
\setlength{\unitlength}{4144sp}%
\begingroup\makeatletter\ifx\SetFigFont\undefined%
\gdef\SetFigFont#1#2#3#4#5{%
  \reset@font\fontsize{#1}{#2pt}%
  \fontfamily{#3}\fontseries{#4}\fontshape{#5}%
  \selectfont}%
\fi\endgroup%
\begin{picture}(5439,2777)(409,-2386)
\put(1731,-781){\makebox(0,0)[lb]{\smash{\SetFigFont{12}{14.4}{\rmdefault}{\mddefault}{\updefault}$2n$}}}
\put(1703,-2386){\makebox(0,0)[lb]{\smash{\SetFigFont{12}{14.4}{\rmdefault}{\mddefault}{\updefault}$2n$}}}
\put(4426,-2311){\makebox(0,0)[lb]{\smash{\SetFigFont{12}{14.4}{\rmdefault}{\mddefault}{\updefault}$2s$}}}
\put(5386,-1088){\makebox(0,0)[lb]{\smash{\SetFigFont{12}{14.4}{\rmdefault}{\mddefault}{\updefault}$q_s$}}}
\put(3781,-1058){\makebox(0,0)[lb]{\smash{\SetFigFont{12}{14.4}{\rmdefault}{\mddefault}{\updefault}$q_1$}}}
\put(3901,-61){\makebox(0,0)[lb]{\smash{\SetFigFont{12}{14.4}{\rmdefault}{\mddefault}{\updefault}by}}}
\end{picture}
$$
\begin{center}
{\bf Fig. 25}
\end{center}

Clearly, $T$-invariants which are not of finite type cannot be extracted

from the generalized Kontsevitch integral (which is the universal invariant

of finite type). Moreover, if we restrict ourselves to $G$-pure global

knots, then the Kontsevitch integral is no longer even the universal

invariant of finite type for these knots.

{\it Exemple $2.5$}
Let $a \in G$ and $a \notin \{0, \pm [K]_G \}$.
Let 
\begin{center}
$\makebox(50,50){D=} \makebox(50,50){%
\begin{picture}(0,0)%
\includegraphics{pag50.pstex}%
\end{picture}%
\setlength{\unitlength}{4144sp}%
\begingroup\makeatletter\ifx\SetFigFont\undefined%
\gdef\SetFigFont#1#2#3#4#5{%
  \reset@font\fontsize{#1}{#2pt}%
  \fontfamily{#3}\fontseries{#4}\fontshape{#5}%
  \selectfont}%
\fi\endgroup%
\begin{picture}(631,890)(564,-336)
\put(758,419){\makebox(0,0)[lb]{\smash{\SetFigFont{10}{12.0}{\rmdefault}{\mddefault}{\updefault}$a$}}}
\put(810,139){\makebox(0,0)[lb]{\smash{\SetFigFont{10}{12.0}{\rmdefault}{\mddefault}{\updefault}$p_1$}}}
\put(999,-31){\makebox(0,0)[lb]{\smash{\SetFigFont{10}{12.0}{\rmdefault}{\mddefault}{\updefault}$p_2$}}}
\put(936,-336){\makebox(0,0)[lb]{\smash{\SetFigFont{10}{12.0}{\rmdefault}{\mddefault}{\updefault}$[K]_G-a$}}}
\end{picture}
}$ \makebox(30,50){ , } $\makebox(50,50){ n=0 }$
\end{center}

$T_K(D; \emptyset)$ is a Gauss diagram invariant of degree 2 for
$G$-pure global knots $K \hookrightarrow T^2 \times \mathbb{R}$.
One easily verifies that $T_K(D; \emptyset)$ is of degree 2 as a

finite type invariant too. If in a (not $G$-pure) isotopy $K_t$,

the knot $K$ crosses exactly once a stratum of type e.g. $a^{\mbox{}}_{\dr}(0|a,-a)$

(see Fig.~26)

$$
\begin{picture}(0,0)%
\includegraphics{ima26.pstex}%
\end{picture}%
\setlength{\unitlength}{4144sp}%
\begingroup\makeatletter\ifx\SetFigFont\undefined%
\gdef\SetFigFont#1#2#3#4#5{%
  \reset@font\fontsize{#1}{#2pt}%
  \fontfamily{#3}\fontseries{#4}\fontshape{#5}%
  \selectfont}%
\fi\endgroup%
\begin{picture}(5095,1884)(345,-1225)
\put(1166,-894){\makebox(0,0)[lb]{\smash{\SetFigFont{12}{14.4}{\rmdefault}{\mddefault}{\updefault}$-a$}}}
\put(1550,-916){\makebox(0,0)[lb]{\smash{\SetFigFont{12}{14.4}{\rmdefault}{\mddefault}{\updefault}$a$}}}
\put(3665,-129){\makebox(0,0)[lb]{\smash{\SetFigFont{12}{14.4}{\rmdefault}{\mddefault}{\updefault}$a$}}}
\put(4871,-122){\makebox(0,0)[lb]{\smash{\SetFigFont{12}{14.4}{\rmdefault}{\mddefault}{\updefault}$-a$}}}
\put(4327,-1085){\makebox(0,0)[lb]{\smash{\SetFigFont{12}{14.4}{\rmdefault}{\mddefault}{\updefault}$0$}}}
\put(1680,-252){\makebox(0,0)[lb]{\smash{\SetFigFont{12}{14.4}{\rmdefault}{\mddefault}{\updefault}$0$}}}
\end{picture}
$$
\begin{center}
{\bf Fig. 26}
\end{center}
 then $T_{K_t}(D; \emptyset)$ changes by $\pm 1$.

Consequently, $T_K(D; \emptyset)$ is not invariant for {\em all\/}

isotopies of $K$ and therefore cannot be extracted from the

generalized Kontsevitch integral.

\par{\bf Remarks}:

\begin{enumerate}

\item

The above invariant $T_K(D; \emptyset)$ could have been equally

considered as an invariant $T_K(\emptyset; c(D))$ with

$$
c(D)= \sum_{\begin{picture}(0,0)%
\includegraphics{pag50.pstex}%
\end{picture}%
\setlength{\unitlength}{4144sp}%
\begingroup\makeatletter\ifx\SetFigFont\undefined%
\gdef\SetFigFont#1#2#3#4#5{%
  \reset@font\fontsize{#1}{#2pt}%
  \fontfamily{#3}\fontseries{#4}\fontshape{#5}%
  \selectfont}%
\fi\endgroup%
\begin{picture}(631,890)(564,-336)
\put(758,419){\makebox(0,0)[lb]{\smash{\SetFigFont{10}{12.0}{\rmdefault}{\mddefault}{\updefault}$a$}}}
\put(810,139){\makebox(0,0)[lb]{\smash{\SetFigFont{10}{12.0}{\rmdefault}{\mddefault}{\updefault}$p_1$}}}
\put(999,-31){\makebox(0,0)[lb]{\smash{\SetFigFont{10}{12.0}{\rmdefault}{\mddefault}{\updefault}$p_2$}}}
\put(936,-336){\makebox(0,0)[lb]{\smash{\SetFigFont{10}{12.0}{\rmdefault}{\mddefault}{\updefault}$[K]_G-a$}}}
\end{picture}
}{w(q_1)w(q_2)}
$$

\item

Let $n= 0$ and let 
\begin{center}
\makebox(50,50){$D=$} \makebox(50,50){%
\begin{picture}(0,0)%
\includegraphics{pag52GH.pstex}%
\end{picture}%
\setlength{\unitlength}{4144sp}%
\begingroup\makeatletter\ifx\SetFigFont\undefined%
\gdef\SetFigFont#1#2#3#4#5{%
  \reset@font\fontsize{#1}{#2pt}%
  \fontfamily{#3}\fontseries{#4}\fontshape{#5}%
  \selectfont}%
\fi\endgroup%
\begin{picture}(631,956)(564,-336)
\put(723,-336){\makebox(0,0)[lb]{\smash{\SetFigFont{10}{12.0}{\rmdefault}{\mddefault}{\updefault}$a$}}}
\put(1003,455){\makebox(0,0)[lb]{\smash{\SetFigFont{10}{12.0}{\rmdefault}{\mddefault}{\updefault}$[K]_G-a$}}}
\end{picture}
}
\makebox(50,50){ or } 
\makebox(50,50){$D=$}\makebox(50,50){%
\begin{picture}(0,0)%
\includegraphics{pag52DH.pstex}%
\end{picture}%
\setlength{\unitlength}{4144sp}%
\begingroup\makeatletter\ifx\SetFigFont\undefined%
\gdef\SetFigFont#1#2#3#4#5{%
  \reset@font\fontsize{#1}{#2pt}%
  \fontfamily{#3}\fontseries{#4}\fontshape{#5}%
  \selectfont}%
\fi\endgroup%
\begin{picture}(631,710)(564,-175)
\put(1025,400){\makebox(0,0)[lb]{\smash{\SetFigFont{10}{12.0}{\rmdefault}{\mddefault}{\updefault}$a$}}}
\put(720,400){\makebox(0,0)[lb]{\smash{\SetFigFont{10}{12.0}{\rmdefault}{\mddefault}{\updefault}$a$}}}
\end{picture}}\makebox(50,50){ .}
\end{center}
The corresponding $T$-invariants $T_K(D; \emptyset)$ are the only other
invariants of degree 2 which are of finite type and which cannot be extracted
from the generalized Kontsevitch integral for knots in $F^2 \times \mathbb{R}$.
\end{enumerate}
\par More generally, let us consider $T$-invariants of finite type under
$G$-pure isotopy from the point of view of the works [V], [K], [BN]. In the
case of knots in $\mathbb{R}^3$, the famous theorem of Kontsevitch states that
each $\mathbb{C}$-valued functional on the $\mathbb{C}$-vector space of
(unmarked) chord diagrams can be integrated to a knot invariant of finite type
if it verifies the 1-$T$ and 4-$T$ relations. (It is an easy matter to see
that these relations are necessarily verified by each invariant of finite type).
\par Let $\mathcal{K}$ be a fixed free homotopy class of an oriented loop in
$F^2$ and let $[K] \in H_1(F^2; \mathbb{Z})$ be the corresponding homology
class. Let $M_{\mathcal{K}}$ be the space of all possibly singular knot
diagrams $K \hookrightarrow F^2 \times \mathbb{R}$ such that $pr(K)$
represents $\mathcal{K}$. Let $G$ be a fixed quotient group of $H_1(F^2;
\mathbb{Z})/\langle [K] \rangle$. Finally, let $M^G_{\mathcal{K}}
\hookrightarrow M_{\mathcal{K}}$ be the subspace of all possibly singular
$G$-pure diagrams. Here, the marking in $G$ of a double point of $K$ is given
by the corresponding marking of the positive resolution (see Fig.~27)
$$
\begin{picture}(0,0)%
\includegraphics{ima27.pstex}%
\end{picture}%
\setlength{\unitlength}{4144sp}%
\begingroup\makeatletter\ifx\SetFigFont\undefined%
\gdef\SetFigFont#1#2#3#4#5{%
  \reset@font\fontsize{#1}{#2pt}%
  \fontfamily{#3}\fontseries{#4}\fontshape{#5}%
  \selectfont}%
\fi\endgroup%
\begin{picture}(3781,1104)(626,-2770)
\put(1231,-2098){\makebox(0,0)[lb]{\smash{\SetFigFont{12}{14.4}{\rmdefault}{\mddefault}{\updefault}$a$}}}
\put(1268,-2585){\makebox(0,0)[lb]{\smash{\SetFigFont{12}{14.4}{\rmdefault}{\mddefault}{\updefault}$p$}}}
\put(3953,-2593){\makebox(0,0)[lb]{\smash{\SetFigFont{12}{14.4}{\rmdefault}{\mddefault}{\updefault}$p$}}}
\put(3983,-2075){\makebox(0,0)[lb]{\smash{\SetFigFont{12}{14.4}{\rmdefault}{\mddefault}{\updefault}$a=[K_p^+] \in G$}}}
\put(1628,-1910){\makebox(0,0)[lb]{\smash{\SetFigFont{12}{14.4}{\rmdefault}{\mddefault}{\updefault}$K$}}}
\end{picture}
$$
\begin{center}
{\bf Fig. 27}
\end{center}
(Near each crossing, we write the corresponding marking in $G$.)
Evidently, $M_{\mathcal{K}}$ is connected. We do not know wether $M^G_
{\mathcal{K}}$ is always connected or not. However, if we take out the set
$\Sigma$ of all singular diagrams in $M^G_{\mathcal{K}}$, then ther are in
general different components of $M^G_{\mathcal{K}} \setminus \Sigma$ which
represent the same knot type in $F^2 \times \mathbb{R}$ (see
sect.~9). Theorem~1 claims that for global knots, $T$-invariants do not depend
on the chosen component of $M^G_{\mathcal{K}} \setminus \Sigma$ for a given
knot type in $F^2 \times \mathbb{R}$. But let us forget global knots for one
moment, and let us consider $G$-pure knots only up to $G$-pure isotopy.
\par Let us have a look at the analogue of the above mentioned relations. Let
$\mathcal{A}_G$ be the $\mathbb{C}$-vector space generated by all chord
diagrams with homological markings in $G$ of the chords, and with the
homotopical marking $\mathcal{K}$ for the whole circle. Finally, let
$\mathcal{A}^0_G \hookrightarrow \mathcal{A}_G$ be the subspace generated by
the $G$-pure chord diagrams (i.e. there are no chords with marking $0 \in
G$). Let $I: \mathcal{A}_G \to \mathbb{C}$ be a functional. We want to integrate it to a knot invariant.
\begin{itemize}
\item {\em $1-T$ relation:\/} This relation is obtained by going in
  $M_{\mathcal{K}}$ around a diagram which has a double point (as a singular knot) in a cusp of the projection into $F^2$. See Fig.~28. 
$$
\makebox(10,80){I}
\makebox(10,80){\Bigg( }
\makebox(60,80){%
\begin{picture}(0,0)%
\includegraphics{ima28G.pstex}%
\end{picture}%
\setlength{\unitlength}{4144sp}%
\begingroup\makeatletter\ifx\SetFigFont\undefined%
\gdef\SetFigFont#1#2#3#4#5{%
  \reset@font\fontsize{#1}{#2pt}%
  \fontfamily{#3}\fontseries{#4}\fontshape{#5}%
  \selectfont}%
\fi\endgroup%
\begin{picture}(512,391)(444,-1707)
\put(571,-1707){\makebox(0,0)[lb]{\smash{\SetFigFont{12}{14.4}{\rmdefault}{\mddefault}{\updefault}$0$}}}
\end{picture}
}
\makebox(10,80){\Bigg) }
\makebox(10,80){-}
\makebox(10,80){I}
\makebox(10,80){\Bigg( }
\makebox(60,80){\begin{picture}(0,0)%
\includegraphics{ima28D.pstex}%
\end{picture}%
\setlength{\unitlength}{4144sp}%
\begingroup\makeatletter\ifx\SetFigFont\undefined%
\gdef\SetFigFont#1#2#3#4#5{%
  \reset@font\fontsize{#1}{#2pt}%
  \fontfamily{#3}\fontseries{#4}\fontshape{#5}%
  \selectfont}%
\fi\endgroup%
\begin{picture}(512,431)(444,-1742)
\put(571,-1742){\makebox(0,0)[lb]{\smash{\SetFigFont{12}{14.4}{\rmdefault}{\mddefault}{\updefault}$0$}}}
\end{picture}
}
\makebox(10,80){\Bigg) }
\makebox(10,80){=}
\makebox(10,80){0}
$$
\begin{center}
{\bf Fig. 28}
\end{center}
\item {\em $2-T$ relation:\/} We have this additional relation because,

instead of considering only embeddings $S^1 \hookrightarrow F^2 \times \mathbb{R}$, we consider diagrams, i.e. embeddings together with the

projection onto $F^2$. Then, the relation is obtained by going in

$M_{\mathcal{K}}$ around a diagram which has two double points in an

autotangency of the projection (see Fig.~29).

$$
\makebox(10,80){I}
\makebox(10,80){\Bigg( }
\makebox(60,80)
{%
\begin{picture}(0,0)%
\includegraphics{ima29G.pstex}%
\end{picture}%
\setlength{\unitlength}{4144sp}%
\begingroup\makeatletter\ifx\SetFigFont\undefined%
\gdef\SetFigFont#1#2#3#4#5{%
  \reset@font\fontsize{#1}{#2pt}%
  \fontfamily{#3}\fontseries{#4}\fontshape{#5}%
  \selectfont}%
\fi\endgroup%
\begin{picture}(549,329)(439,-2082)
\put(511,-2082){\makebox(0,0)[lb]{\smash{\SetFigFont{12}{14.4}{\rmdefault}{\mddefault}{\updefault}$a$}}}
\end{picture}
}\makebox(10,80){ \Bigg) }
\makebox(10,80){-}
\makebox(10,80){I}
\makebox(10,80){ \Bigg( }
\makebox(60,80)
{%
\begin{picture}(0,0)%
\includegraphics{ima29D.pstex}%
\end{picture}%
\setlength{\unitlength}{4144sp}%
\begingroup\makeatletter\ifx\SetFigFont\undefined%
\gdef\SetFigFont#1#2#3#4#5{%
  \reset@font\fontsize{#1}{#2pt}%
  \fontfamily{#3}\fontseries{#4}\fontshape{#5}%
  \selectfont}%
\fi\endgroup%
\begin{picture}(549,204)(439,-1517)
\end{picture}}
\makebox(10,80){ \Bigg) }
\makebox(10,80){=}
\makebox(10,80){0}
$$
\begin{center}
{\bf Fig. 29}
\end{center}

(The markings in all other crossings or double points of these

fragments are determined by the single marking $a \in G$.)

There are two cases to consider, which are shown in Fig.~30 and 31.
$$
\makebox(10,80){\Bigg(}
\makebox(10,80){I}
\makebox(10,80){\Bigg(}
\makebox(60,80){\begin{picture}(0,0)%
\includegraphics{ima301.pstex}%
\end{picture}%
\setlength{\unitlength}{4144sp}%
\begingroup\makeatletter\ifx\SetFigFont\undefined%
\gdef\SetFigFont#1#2#3#4#5{%
  \reset@font\fontsize{#1}{#2pt}%
  \fontfamily{#3}\fontseries{#4}\fontshape{#5}%
  \selectfont}%
\fi\endgroup%
\begin{picture}(553,357)(443,232)
\put(499,232){\makebox(0,0)[lb]{\smash{\SetFigFont{12}{14.4}{\rmdefault}{\mddefault}{\updefault}$a$}}}
\put(833,246){\makebox(0,0)[lb]{\smash{\SetFigFont{12}{14.4}{\rmdefault}{\mddefault}{\updefault}$a$}}}
\end{picture}
}
\makebox(10,80){\Bigg)}
\makebox(10,80){-}
\makebox(10,80){I}
\makebox(10,80){\Bigg(}
\makebox(60,80){\begin{picture}(0,0)%
\includegraphics{ima302.pstex}%
\end{picture}%
\setlength{\unitlength}{4144sp}%
\begingroup\makeatletter\ifx\SetFigFont\undefined%
\gdef\SetFigFont#1#2#3#4#5{%
  \reset@font\fontsize{#1}{#2pt}%
  \fontfamily{#3}\fontseries{#4}\fontshape{#5}%
  \selectfont}%
\fi\endgroup%
\begin{picture}(662,390)(431,-218)
\put(516,-207){\makebox(0,0)[lb]{\smash{\SetFigFont{12}{14.4}{\rmdefault}{\mddefault}{\updefault}$-a$}}}
\put(898,-218){\makebox(0,0)[lb]{\smash{\SetFigFont{12}{14.4}{\rmdefault}{\mddefault}{\updefault}$a$}}}
\end{picture}
}
\makebox(10,80){\Bigg)}
\makebox(10,80){\Bigg)}
\makebox(10,80){-}
$$
$$
\makebox(10,80){\Bigg(}
\makebox(10,80){I}
\makebox(10,80){\Bigg(}
\makebox(60,80){\begin{picture}(0,0)%
\includegraphics{ima303.pstex}%
\end{picture}%
\setlength{\unitlength}{4144sp}%
\begingroup\makeatletter\ifx\SetFigFont\undefined%
\gdef\SetFigFont#1#2#3#4#5{%
  \reset@font\fontsize{#1}{#2pt}%
  \fontfamily{#3}\fontseries{#4}\fontshape{#5}%
  \selectfont}%
\fi\endgroup%
\begin{picture}(553,357)(443,232)
\put(499,232){\makebox(0,0)[lb]{\smash{\SetFigFont{12}{14.4}{\rmdefault}{\mddefault}{\updefault}$-a$}}}
\put(833,246){\makebox(0,0)[lb]{\smash{\SetFigFont{12}{14.4}{\rmdefault}{\mddefault}{\updefault}$-a$}}}
\end{picture}
}
\makebox(10,80){\Bigg)}
\makebox(10,80){-}
\makebox(10,80){I}
\makebox(10,80){\Bigg(}
\makebox(60,80){\begin{picture}(0,0)%
\includegraphics{ima302.pstex}%
\end{picture}%
\setlength{\unitlength}{4144sp}%
\begingroup\makeatletter\ifx\SetFigFont\undefined%
\gdef\SetFigFont#1#2#3#4#5{%
  \reset@font\fontsize{#1}{#2pt}%
  \fontfamily{#3}\fontseries{#4}\fontshape{#5}%
  \selectfont}%
\fi\endgroup%
\begin{picture}(662,390)(431,-218)
\put(516,-207){\makebox(0,0)[lb]{\smash{\SetFigFont{12}{14.4}{\rmdefault}{\mddefault}{\updefault}$-a$}}}
\put(898,-218){\makebox(0,0)[lb]{\smash{\SetFigFont{12}{14.4}{\rmdefault}{\mddefault}{\updefault}$a$}}}
\end{picture}
}
\makebox(10,80){\Bigg)}
\makebox(10,80){\Bigg)}
\makebox(10,80){=}
\makebox(10,80){0}
$$
\begin{center}
{\bf Fig. 30}
\end{center}

$$
\makebox(10,80){\Bigg(}
\makebox(10,80){I}
\makebox(10,80){\Bigg(}
\makebox(60,80){\begin{picture}(0,0)%
\includegraphics{ima311.pstex}%
\end{picture}%
\setlength{\unitlength}{4144sp}%
\begingroup\makeatletter\ifx\SetFigFont\undefined%
\gdef\SetFigFont#1#2#3#4#5{%
  \reset@font\fontsize{#1}{#2pt}%
  \fontfamily{#3}\fontseries{#4}\fontshape{#5}%
  \selectfont}%
\fi\endgroup%
\begin{picture}(662,409)(431,-237)
\put(887,-218){\makebox(0,0)[lb]{\smash{\SetFigFont{12}{14.4}{\rmdefault}{\mddefault}{\updefault}$-a$}}}
\put(512,-237){\makebox(0,0)[lb]{\smash{\SetFigFont{12}{14.4}{\rmdefault}{\mddefault}{\updefault}$a$}}}
\end{picture}
}
\makebox(10,80){\Bigg)}
\makebox(10,80){-}
\makebox(10,80){I}
\makebox(10,80){\Bigg(}
\makebox(60,80){\begin{picture}(0,0)%
\includegraphics{ima312.pstex}%
\end{picture}%
\setlength{\unitlength}{4144sp}%
\begingroup\makeatletter\ifx\SetFigFont\undefined%
\gdef\SetFigFont#1#2#3#4#5{%
  \reset@font\fontsize{#1}{#2pt}%
  \fontfamily{#3}\fontseries{#4}\fontshape{#5}%
  \selectfont}%
\fi\endgroup%
\begin{picture}(553,338)(435,212)
\put(507,213){\makebox(0,0)[lb]{\smash{\SetFigFont{12}{14.4}{\rmdefault}{\mddefault}{\updefault}$-a$}}}
\put(852,212){\makebox(0,0)[lb]{\smash{\SetFigFont{12}{14.4}{\rmdefault}{\mddefault}{\updefault}$-a$}}}
\end{picture}
}
\makebox(10,80){\Bigg)}
\makebox(10,80){\Bigg)}
\makebox(10,80){-}
$$
$$
\makebox(10,80){\Bigg(}
\makebox(10,80){I}
\makebox(10,80){\Bigg(}
\makebox(60,80){\begin{picture}(0,0)%
\includegraphics{ima311.pstex}%
\end{picture}%
\setlength{\unitlength}{4144sp}%
\begingroup\makeatletter\ifx\SetFigFont\undefined%
\gdef\SetFigFont#1#2#3#4#5{%
  \reset@font\fontsize{#1}{#2pt}%
  \fontfamily{#3}\fontseries{#4}\fontshape{#5}%
  \selectfont}%
\fi\endgroup%
\begin{picture}(662,409)(431,-237)
\put(887,-218){\makebox(0,0)[lb]{\smash{\SetFigFont{12}{14.4}{\rmdefault}{\mddefault}{\updefault}$-a$}}}
\put(512,-237){\makebox(0,0)[lb]{\smash{\SetFigFont{12}{14.4}{\rmdefault}{\mddefault}{\updefault}$a$}}}
\end{picture}
}
\makebox(10,80){\Bigg)}
\makebox(10,80){-}
\makebox(10,80){I}
\makebox(10,80){\Bigg(}
\makebox(60,80){\begin{picture}(0,0)%
\includegraphics{ima314.pstex}%
\end{picture}%
\setlength{\unitlength}{4144sp}%
\begingroup\makeatletter\ifx\SetFigFont\undefined%
\gdef\SetFigFont#1#2#3#4#5{%
  \reset@font\fontsize{#1}{#2pt}%
  \fontfamily{#3}\fontseries{#4}\fontshape{#5}%
  \selectfont}%
\fi\endgroup%
\begin{picture}(553,365)(443,232)
\put(499,232){\makebox(0,0)[lb]{\smash{\SetFigFont{12}{14.4}{\rmdefault}{\mddefault}{\updefault}$a$}}}
\put(833,246){\makebox(0,0)[lb]{\smash{\SetFigFont{12}{14.4}{\rmdefault}{\mddefault}{\updefault}$a$}}}
\end{picture}
}
\makebox(10,80){\Bigg)}
\makebox(10,80){\Bigg)}
\makebox(10,80){=}
\makebox(10,80){0}
$$
\begin{center}
{\bf Fig. 31}
\end{center}

\item {\em $4-T$ relation:\/} This relation is obtained by going in

$M_{\mathcal{K}}$ around a diagram with three double points in

a triple point of the projection. For each $a$, $b \in G$, we

have the relation shown in Fig.~32.
$$
\makebox(10,80){I}
\makebox(10,80){\Bigg(}
\makebox(120,80){\begin{picture}(0,0)%
\includegraphics{ima321.pstex}%
\end{picture}%
\setlength{\unitlength}{4144sp}%
\begingroup\makeatletter\ifx\SetFigFont\undefined%
\gdef\SetFigFont#1#2#3#4#5{%
  \reset@font\fontsize{#1}{#2pt}%
  \fontfamily{#3}\fontseries{#4}\fontshape{#5}%
  \selectfont}%
\fi\endgroup%
\begin{picture}(1470,1245)(116,-921)
\put(1051,-598){\makebox(0,0)[b]{\smash{\SetFigFont{12}{14.4}{\rmdefault}{\mddefault}{\updefault}$b$}}}
\put(400,-604){\makebox(0,0)[b]{\smash{\SetFigFont{12}{14.4}{\rmdefault}{\mddefault}{\updefault}$a$}}}
\put(952,-34){\makebox(0,0)[lb]{\smash{\SetFigFont{12}{14.4}{\rmdefault}{\mddefault}{\updefault}$a-b$}}}
\end{picture}
}
\makebox(10,80){\Bigg)}
\makebox(10,80){-}
\makebox(10,80){I}
\makebox(10,80){\Bigg(}
\makebox(120,80){\begin{picture}(0,0)%
\includegraphics{ima322.pstex}%
\end{picture}%
\setlength{\unitlength}{4144sp}%
\begingroup\makeatletter\ifx\SetFigFont\undefined%
\gdef\SetFigFont#1#2#3#4#5{%
  \reset@font\fontsize{#1}{#2pt}%
  \fontfamily{#3}\fontseries{#4}\fontshape{#5}%
  \selectfont}%
\fi\endgroup%
\begin{picture}(1556,1116)(199,-624)
\put(1108,213){\makebox(0,0)[b]{\smash{\SetFigFont{12}{14.4}{\rmdefault}{\mddefault}{\updefault}$a$}}}
\put(1004,-348){\makebox(0,0)[lb]{\smash{\SetFigFont{12}{14.4}{\rmdefault}{\mddefault}{\updefault}$b-a$}}}
\put(818,227){\makebox(0,0)[b]{\smash{\SetFigFont{12}{14.4}{\rmdefault}{\mddefault}{\updefault}$b$}}}
\end{picture}
}
\makebox(10,80){\Bigg)}
$$

$$
\makebox(10,80){+}
\makebox(10,80){I}
\makebox(10,80){\Bigg(}
\makebox(120,80){
\begin{picture}(0,0)%
\includegraphics{ima323.pstex}%
\end{picture}%
\setlength{\unitlength}{4144sp}%
\begingroup\makeatletter\ifx\SetFigFont\undefined%
\gdef\SetFigFont#1#2#3#4#5{%
  \reset@font\fontsize{#1}{#2pt}%
  \fontfamily{#3}\fontseries{#4}\fontshape{#5}%
  \selectfont}%
\fi\endgroup%
\begin{picture}(1510,1304)(1294,-1107)
\put(2084,-179){\makebox(0,0)[lb]{\smash{\SetFigFont{12}{14.4}{\rmdefault}{\mddefault}{\updefault}$b-a$}}}
\put(1498,-698){\makebox(0,0)[b]{\smash{\SetFigFont{12}{14.4}{\rmdefault}{\mddefault}{\updefault}$a$}}}
\put(2255,-805){\makebox(0,0)[b]{\smash{\SetFigFont{12}{14.4}{\rmdefault}{\mddefault}{\updefault}$b$}}}
\end{picture}
}
\makebox(10,80){\Bigg)}
\makebox(10,80){-}
\makebox(10,80){I}
\makebox(10,80){\Bigg(}
\makebox(120,80){
\begin{picture}(0,0)%
\includegraphics{ima324.pstex}%
\end{picture}%
\setlength{\unitlength}{4144sp}%
\begingroup\makeatletter\ifx\SetFigFont\undefined%
\gdef\SetFigFont#1#2#3#4#5{%
  \reset@font\fontsize{#1}{#2pt}%
  \fontfamily{#3}\fontseries{#4}\fontshape{#5}%
  \selectfont}%
\fi\endgroup%
\begin{picture}(1556,1116)(199,-624)
\put(1205,-11){\makebox(0,0)[b]{\smash{\SetFigFont{12}{14.4}{\rmdefault}{\mddefault}{\updefault}$a$}}}
\put(1034,-341){\makebox(0,0)[lb]{\smash{\SetFigFont{12}{14.4}{\rmdefault}{\mddefault}{\updefault}$b-a$}}}
\put(740,216){\makebox(0,0)[b]{\smash{\SetFigFont{12}{14.4}{\rmdefault}{\mddefault}{\updefault}$-b$}}}
\end{picture}
}
\makebox(10,80){\Bigg)}
\makebox(10,80){=}
\makebox(10,80){0}
$$
\begin{center}
{\bf Fig. 32}
\end{center}

\end{itemize}
The proofs are completely analogous to the one for knots in $\mathbb{R}^3$
(see also [Go] where it is done for knots in the solid torus). Each functional $I$ which can be integrated to a knot invariant verifies $1-T$, $2-T$, $4-T$.
\par {\sl Question}: {\it Can each functional which verifies $1-T$, $2-T$, $4-T$ be integrated to a knot invariant?\/}
\par {\bf Remark}.
Goryunov [Go] has shown that the answer is "yes" in the case of the solid
torus. Notice that our chord diagrams are planar and with homological
markings. In [A-M-R] it is shown that the answer to the above question is
"yes" if $\partial F^2 \not= \emptyset$ and if one uses chord diagrams which
are immersed in $F^2$ instead of planar chord diagrams with homological
markings. But then, it seems to be difficult to find such functionals.
\par Let us consider now functionals $I^0: \mathcal{A}^0_G: \to
\mathbb{C}$. Evidently, if $I^0$ can be integrated to a knot invariant under
$G$-pure isotopy, then it has to verify the corresponding relations $1-T^0=
\emptyset$, $2-T^0$, $4-T^0$ similar to the previous ones but where the marking $0 \in G$ is forbidden.
par {\it Example~3.6}
$$
\makebox(50,50){$I^0_1:= T_K \Bigg($}\makebox(50,50){\begin{picture}(0,0)%
\includegraphics{p52b8.pstex}%
\end{picture}%
\setlength{\unitlength}{4144sp}%
\begingroup\makeatletter\ifx\SetFigFont\undefined%
\gdef\SetFigFont#1#2#3#4#5{%
  \reset@font\fontsize{#1}{#2pt}%
  \fontfamily{#3}\fontseries{#4}\fontshape{#5}%
  \selectfont}%
\fi\endgroup%
\begin{picture}(631,833)(564,-298)
\put(1025,400){\makebox(0,0)[lb]{\smash{\SetFigFont{10}{12.0}{\rmdefault}{\mddefault}{\updefault}$a$}}}
\put(727,-298){\makebox(0,0)[lb]{\smash{\SetFigFont{10}{12.0}{\rmdefault}{\mddefault}{\updefault}$-a$}}}
\end{picture}
}\makebox(30,50){; $\emptyset$ }\makebox(10,50){$\Bigg)$}\makebox(10,50){,} \makebox(50,50)
{$a\neq 0$}\makebox(10,50){.}
$$
$I^0_1$ verifies $1-T$, $2-T$, and $4-T^0$ but it does {\em not\/} always verify $4-T$ (see Fig.~33).
$$
\makebox(10,80){I}
\makebox(10,80){\Bigg(}
\makebox(120,80){\begin{picture}(0,0)%
\includegraphics{ima321.pstex}%
\end{picture}%
\setlength{\unitlength}{4144sp}%
\begingroup\makeatletter\ifx\SetFigFont\undefined%
\gdef\SetFigFont#1#2#3#4#5{%
  \reset@font\fontsize{#1}{#2pt}%
  \fontfamily{#3}\fontseries{#4}\fontshape{#5}%
  \selectfont}%
\fi\endgroup%
\begin{picture}(1470,1245)(116,-921)
\put(1051,-598){\makebox(0,0)[b]{\smash{\SetFigFont{12}{14.4}{\rmdefault}{\mddefault}{\updefault}$a$}}}
\put(400,-604){\makebox(0,0)[b]{\smash{\SetFigFont{12}{14.4}{\rmdefault}{\mddefault}{\updefault}$a$}}}
\put(952,-34){\makebox(0,0)[lb]{\smash{\SetFigFont{12}{14.4}{\rmdefault}{\mddefault}{\updefault}$0$}}}
\end{picture}
}
\makebox(10,80){\Bigg)}
\makebox(10,80){-}
\makebox(10,80){I}
\makebox(10,80){\Bigg(}
\makebox(120,80){\begin{picture}(0,0)%
\includegraphics{ima322.pstex}%
\end{picture}%
\setlength{\unitlength}{4144sp}%
\begingroup\makeatletter\ifx\SetFigFont\undefined%
\gdef\SetFigFont#1#2#3#4#5{%
  \reset@font\fontsize{#1}{#2pt}%
  \fontfamily{#3}\fontseries{#4}\fontshape{#5}%
  \selectfont}%
\fi\endgroup%
\begin{picture}(1556,1116)(199,-624)
\put(1108,213){\makebox(0,0)[b]{\smash{\SetFigFont{12}{14.4}{\rmdefault}{\mddefault}{\updefault}$a$}}}
\put(1004,-348){\makebox(0,0)[lb]{\smash{\SetFigFont{12}{14.4}{\rmdefault}{\mddefault}{\updefault}$0$}}}
\put(818,227){\makebox(0,0)[b]{\smash{\SetFigFont{12}{14.4}{\rmdefault}{\mddefault}{\updefault}$a$}}}
\end{picture}
}
\makebox(10,80){\Bigg)}
$$
$$
\makebox(10,80){+}
\makebox(10,80){I}
\makebox(10,80){\Bigg(}
\makebox(120,80){
\begin{picture}(0,0)%
\includegraphics{ima323.pstex}%
\end{picture}%
\setlength{\unitlength}{4144sp}%
\begingroup\makeatletter\ifx\SetFigFont\undefined%
\gdef\SetFigFont#1#2#3#4#5{%
  \reset@font\fontsize{#1}{#2pt}%
  \fontfamily{#3}\fontseries{#4}\fontshape{#5}%
  \selectfont}%
\fi\endgroup%
\begin{picture}(1510,1304)(1294,-1107)
\put(2084,-179){\makebox(0,0)[lb]{\smash{\SetFigFont{12}{14.4}{\rmdefault}{\mddefault}{\updefault}$0$}}}
\put(1498,-698){\makebox(0,0)[b]{\smash{\SetFigFont{12}{14.4}{\rmdefault}{\mddefault}{\updefault}$a$}}}
\put(2255,-805){\makebox(0,0)[b]{\smash{\SetFigFont{12}{14.4}{\rmdefault}{\mddefault}{\updefault}$a$}}}
\end{picture}
}
\makebox(10,80){\Bigg)}
\makebox(10,80){-}
\makebox(10,80){I}
\makebox(10,80){\Bigg(}
\makebox(120,80){
\begin{picture}(0,0)%
\includegraphics{ima324.pstex}%
\end{picture}%
\setlength{\unitlength}{4144sp}%
\begingroup\makeatletter\ifx\SetFigFont\undefined%
\gdef\SetFigFont#1#2#3#4#5{%
  \reset@font\fontsize{#1}{#2pt}%
  \fontfamily{#3}\fontseries{#4}\fontshape{#5}%
  \selectfont}%
\fi\endgroup%
\begin{picture}(1556,1116)(199,-624)
\put(1205,-11){\makebox(0,0)[b]{\smash{\SetFigFont{12}{14.4}{\rmdefault}{\mddefault}{\updefault}$a$}}}
\put(1034,-341){\makebox(0,0)[lb]{\smash{\SetFigFont{12}{14.4}{\rmdefault}{\mddefault}{\updefault}$0$}}}
\put(740,216){\makebox(0,0)[b]{\smash{\SetFigFont{12}{14.4}{\rmdefault}{\mddefault}{\updefault}$-a$}}}
\end{picture}
}
\makebox(10,80){$\Bigg)$}
\makebox(10,80){$=$}
\makebox(10,80){$-1$}
$$
\begin{center}
{\bf Fig. 33}
\end{center}
Indeed, the only non-zero contribution to $I^0_1$ from two of the three
involved crossings comes from the term shown in Fig.~34, and is equal to $w(p)w(q)= -1$.
$$
\makebox(10,80){-}
\makebox(10,80){\Bigg(}
\makebox(10,80){$-I_1^0$}
\makebox(10,80){\Bigg(}
\makebox(120,80){\begin{picture}(0,0)%
\includegraphics{ima34.pstex}%
\end{picture}%
\setlength{\unitlength}{4144sp}%
\begingroup\makeatletter\ifx\SetFigFont\undefined%
\gdef\SetFigFont#1#2#3#4#5{%
  \reset@font\fontsize{#1}{#2pt}%
  \fontfamily{#3}\fontseries{#4}\fontshape{#5}%
  \selectfont}%
\fi\endgroup%
\begin{picture}(1556,1116)(199,-624)
\put(1205,-11){\makebox(0,0)[b]{\smash{\SetFigFont{12}{14.4}{\rmdefault}{\mddefault}{\updefault}$-a$}}}
\put(1034,-341){\makebox(0,0)[lb]{\smash{\SetFigFont{12}{14.4}{\rmdefault}{\mddefault}{\updefault}$0$}}}
\put(740,216){\makebox(0,0)[b]{\smash{\SetFigFont{12}{14.4}{\rmdefault}{\mddefault}{\updefault}$a$}}}
\end{picture}
}
\makebox(10,80){\Bigg)}
\makebox(10,80){\Bigg)}
$$
\begin{center}
{\bf Fig. 34}
\end{center}
{\it Exemple $3.7$}
$$
\makebox(50,50){$I^0_2:= T_K \Bigg($}\makebox(50,50){\begin{picture}(0,0)%
\includegraphics{pag52DH.pstex}%
\end{picture}%
\setlength{\unitlength}{4144sp}%
\begingroup\makeatletter\ifx\SetFigFont\undefined%
\gdef\SetFigFont#1#2#3#4#5{%
  \reset@font\fontsize{#1}{#2pt}%
  \fontfamily{#3}\fontseries{#4}\fontshape{#5}%
  \selectfont}%
\fi\endgroup%
\begin{picture}(631,710)(564,-175)
\put(1025,400){\makebox(0,0)[lb]{\smash{\SetFigFont{10}{12.0}{\rmdefault}{\mddefault}{\updefault}$a$}}}
\put(720,400){\makebox(0,0)[lb]{\smash{\SetFigFont{10}{12.0}{\rmdefault}{\mddefault}{\updefault}$a$}}}
\end{picture}
}\makebox(20,50){; }\makebox(50,50){$\emptyset \Bigg)$}
$$
with $a \not= 0$ and $a \not= -a$ in $G$. $I^0_2$ verifies $1-T$ and $4-T^0$
but it does {\em not\/} always verify $2-T^0$! Indeed, in case 2, we have the relation shown in Fig.~35.
$$
\makebox(10,80){\Bigg(}
\makebox(10,80){$I_2^0$}
\makebox(10,80){\Bigg(}
\makebox(60,80){\begin{picture}(0,0)%
\includegraphics{ima311.pstex}%
\end{picture}%
\setlength{\unitlength}{4144sp}%
\begingroup\makeatletter\ifx\SetFigFont\undefined%
\gdef\SetFigFont#1#2#3#4#5{%
  \reset@font\fontsize{#1}{#2pt}%
  \fontfamily{#3}\fontseries{#4}\fontshape{#5}%
  \selectfont}%
\fi\endgroup%
\begin{picture}(662,409)(431,-237)
\put(887,-218){\makebox(0,0)[lb]{\smash{\SetFigFont{12}{14.4}{\rmdefault}{\mddefault}{\updefault}$-a$}}}
\put(512,-237){\makebox(0,0)[lb]{\smash{\SetFigFont{12}{14.4}{\rmdefault}{\mddefault}{\updefault}$a$}}}
\end{picture}
}
\makebox(10,80){\Bigg)}
\makebox(10,80){-}
\makebox(10,80){$I_2^0$}
\makebox(10,80){\Bigg(}
\makebox(60,80){\begin{picture}(0,0)%
\includegraphics{ima312.pstex}%
\end{picture}%
\setlength{\unitlength}{4144sp}%
\begingroup\makeatletter\ifx\SetFigFont\undefined%
\gdef\SetFigFont#1#2#3#4#5{%
  \reset@font\fontsize{#1}{#2pt}%
  \fontfamily{#3}\fontseries{#4}\fontshape{#5}%
  \selectfont}%
\fi\endgroup%
\begin{picture}(553,338)(435,212)
\put(507,213){\makebox(0,0)[lb]{\smash{\SetFigFont{12}{14.4}{\rmdefault}{\mddefault}{\updefault}$-a$}}}
\put(852,212){\makebox(0,0)[lb]{\smash{\SetFigFont{12}{14.4}{\rmdefault}{\mddefault}{\updefault}$-a$}}}
\end{picture}
}
\makebox(10,80){\Bigg)}
\makebox(10,80){\Bigg)}
\makebox(10,80){-}
$$
$$
\makebox(10,80){\Bigg(}
\makebox(10,80){$I_2^0$}
\makebox(10,80){\Bigg(}
\makebox(60,80){\begin{picture}(0,0)%
\includegraphics{ima311.pstex}%
\end{picture}%
\setlength{\unitlength}{4144sp}%
\begingroup\makeatletter\ifx\SetFigFont\undefined%
\gdef\SetFigFont#1#2#3#4#5{%
  \reset@font\fontsize{#1}{#2pt}%
  \fontfamily{#3}\fontseries{#4}\fontshape{#5}%
  \selectfont}%
\fi\endgroup%
\begin{picture}(662,409)(431,-237)
\put(887,-218){\makebox(0,0)[lb]{\smash{\SetFigFont{12}{14.4}{\rmdefault}{\mddefault}{\updefault}$-a$}}}
\put(512,-237){\makebox(0,0)[lb]{\smash{\SetFigFont{12}{14.4}{\rmdefault}{\mddefault}{\updefault}$a$}}}
\end{picture}
}
\makebox(10,80){\Bigg)}
\makebox(10,80){-}
\makebox(10,80){$I_2^0$}
\makebox(10,80){\Bigg(}
\makebox(60,80){\begin{picture}(0,0)%
\includegraphics{ima314.pstex}%
\end{picture}%
\setlength{\unitlength}{4144sp}%
\begingroup\makeatletter\ifx\SetFigFont\undefined%
\gdef\SetFigFont#1#2#3#4#5{%
  \reset@font\fontsize{#1}{#2pt}%
  \fontfamily{#3}\fontseries{#4}\fontshape{#5}%
  \selectfont}%
\fi\endgroup%
\begin{picture}(553,365)(443,232)
\put(499,232){\makebox(0,0)[lb]{\smash{\SetFigFont{12}{14.4}{\rmdefault}{\mddefault}{\updefault}$a$}}}
\put(833,246){\makebox(0,0)[lb]{\smash{\SetFigFont{12}{14.4}{\rmdefault}{\mddefault}{\updefault}$a$}}}
\end{picture}
}
\makebox(10,80){\Bigg)}
\makebox(10,80){\Bigg)}
\makebox(10,80){=}
\makebox(10,80){$-1$}
$$
\begin{center}
{\bf Fig. 35}
\end{center}
However, $I^0_2$ is an isotopy invariant for $G$-pure global knots because
autotangencies with opposite directions of the branches do not occur in an isotopy through global knots.
par {\sl Question}. {\it If a functional $I^0: \mathcal{A}^0_G \to \mathbb{C}$
  verifies $2-T^0$ and $4-T^0$, can it be integrated to a knot invariant under $G$-pure isotopy?\/}
\par Let us consider now Gauss diagram invariants of degree 2 for global knots
$K$, which are not of finite type. From now on, $G$ will always be a quotient
group of $H_1(F^2; \mathbb{Z})/\langle [K] \rangle$. Hence, a knot is $G$-pure
if and only if there is no marking equal to $0$ in $G$. By the above lemmas
and the definition of the degree of a Gauss diagram invariant, we must have
$m= n= 1$. Consequently, 
$$ 
\makebox(60,50){$D=$}\makebox(60,50){\begin{picture}(0,0)%
\includegraphics{pag52M.pstex}%
\end{picture}%
\setlength{\unitlength}{4144sp}%
\begingroup\makeatletter\ifx\SetFigFont\undefined%
\gdef\SetFigFont#1#2#3#4#5{%
  \reset@font\fontsize{#1}{#2pt}%
  \fontfamily{#3}\fontseries{#4}\fontshape{#5}%
  \selectfont}%
\fi\endgroup%
\begin{picture}(772,947)(494,-258)
\put(931, 39){\makebox(0,0)[lb]{\smash{\SetFigFont{12}{14.4}{\rmdefault}{\mddefault}{\updefault}$p$}}}
\put(893,524){\makebox(0,0)[lb]{\smash{\SetFigFont{12}{14.4}{\rmdefault}{\mddefault}{\updefault}$a$}}}
\end{picture}}
$$
 for some $a \in G$, $a \not= 0$.
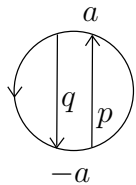
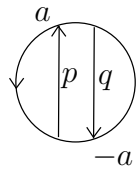
\begin{defin}
\begin{displaymath}
c_{++}(D)= \sum_{\begin{picture}(0,0)%
\includegraphics{page52B.pstex}%
\end{picture}%
\setlength{\unitlength}{4144sp}%
\begingroup\makeatletter\ifx\SetFigFont\undefined%
\gdef\SetFigFont#1#2#3#4#5{%
  \reset@font\fontsize{#1}{#2pt}%
  \fontfamily{#3}\fontseries{#4}\fontshape{#5}%
  \selectfont}%
\fi\endgroup%
\begin{picture}(772,1131)(494,-428)
\put(758,-428){\makebox(0,0)[lb]{\smash{\SetFigFont{12}{14.4}{\rmdefault}{\mddefault}{\updefault}$-a$}}}
\put(826, 29){\makebox(0,0)[lb]{\smash{\SetFigFont{12}{14.4}{\rmdefault}{\mddefault}{\updefault}$q$}}}
\put(1043,-72){\makebox(0,0)[lb]{\smash{\SetFigFont{12}{14.4}{\rmdefault}{\mddefault}{\updefault}$p$}}}
\put(953,538){\makebox(0,0)[lb]{\smash{\SetFigFont{12}{14.4}{\rmdefault}{\mddefault}{\updefault}$a$}}}
\end{picture}
}{w(q)} \qquad
c_{--}(d)= \sum_{\begin{picture}(0,0)%
\includegraphics{page53H.pstex}%
\end{picture}%
\setlength{\unitlength}{4144sp}%
\begingroup\makeatletter\ifx\SetFigFont\undefined%
\gdef\SetFigFont#1#2#3#4#5{%
  \reset@font\fontsize{#1}{#2pt}%
  \fontfamily{#3}\fontseries{#4}\fontshape{#5}%
  \selectfont}%
\fi\endgroup%
\begin{picture}(772,1033)(494,-383)
\put(1002,-383){\makebox(0,0)[lb]{\smash{\SetFigFont{12}{14.4}{\rmdefault}{\mddefault}{\updefault}$-a$}}}
\put(656,485){\makebox(0,0)[lb]{\smash{\SetFigFont{12}{14.4}{\rmdefault}{\mddefault}{\updefault}$a$}}}
\put(1036,104){\makebox(0,0)[lb]{\smash{\SetFigFont{12}{14.4}{\rmdefault}{\mddefault}{\updefault}$q$}}}
\put(818,119){\makebox(0,0)[lb]{\smash{\SetFigFont{12}{14.4}{\rmdefault}{\mddefault}{\updefault}$p$}}}
\end{picture}
}{w(q)} \qquad
\end{displaymath}
\begin{displaymath}
c_{+-}(D)= \sum_{\begin{picture}(0,0)%
\includegraphics{page53M.pstex}%
\end{picture}%
\setlength{\unitlength}{4144sp}%
\begingroup\makeatletter\ifx\SetFigFont\undefined%
\gdef\SetFigFont#1#2#3#4#5{%
  \reset@font\fontsize{#1}{#2pt}%
  \fontfamily{#3}\fontseries{#4}\fontshape{#5}%
  \selectfont}%
\fi\endgroup%
\begin{picture}(772,901)(494,-251)
\put(656,485){\makebox(0,0)[lb]{\smash{\SetFigFont{12}{14.4}{\rmdefault}{\mddefault}{\updefault}$a$}}}
\put(1028,480){\makebox(0,0)[lb]{\smash{\SetFigFont{12}{14.4}{\rmdefault}{\mddefault}{\updefault}$a$}}}
\put(1043,-20){\makebox(0,0)[lb]{\smash{\SetFigFont{12}{14.4}{\rmdefault}{\mddefault}{\updefault}$p$}}}
\put(834,136){\makebox(0,0)[lb]{\smash{\SetFigFont{12}{14.4}{\rmdefault}{\mddefault}{\updefault}$q$}}}
\end{picture}
}{w(q)} \qquad
c_{-+}(D)= \sum_{\begin{picture}(0,0)%
\includegraphics{page53B.pstex}%
\end{picture}%
\setlength{\unitlength}{4144sp}%
\begingroup\makeatletter\ifx\SetFigFont\undefined%
\gdef\SetFigFont#1#2#3#4#5{%
  \reset@font\fontsize{#1}{#2pt}%
  \fontfamily{#3}\fontseries{#4}\fontshape{#5}%
  \selectfont}%
\fi\endgroup%
\begin{picture}(772,901)(494,-251)
\put(656,485){\makebox(0,0)[lb]{\smash{\SetFigFont{12}{14.4}{\rmdefault}{\mddefault}{\updefault}$a$}}}
\put(1028,480){\makebox(0,0)[lb]{\smash{\SetFigFont{12}{14.4}{\rmdefault}{\mddefault}{\updefault}$a$}}}
\put(1056,139){\makebox(0,0)[lb]{\smash{\SetFigFont{12}{14.4}{\rmdefault}{\mddefault}{\updefault}$q$}}}
\put(833,-17){\makebox(0,0)[lb]{\smash{\SetFigFont{12}{14.4}{\rmdefault}{\mddefault}{\updefault}$p$}}}
\end{picture}
}{w(q)}
\end{displaymath}
\end{defin}
\begin{lem}
$c_{++}(D)$, $c_{--}(D)$, $c_{+-}(D)$, $c_{-+}(D)$
{\it are $G$-pure classes of degree 1 of $D$. Moreover, these are the only $G$-pure classes of degree 1.\/}
\end{lem}
\par {\bf Proof}. Evidently, $c_{++}$, $c_{--}$, $c_{+-}$, $c_{-+}$ are
$G$-pure classes of $D$ because $p$ and $q$ cannot get crossed in a $G$-pure
isotopy. Assume now that $p$ and $q$ are crossed. None of the configurations depicted in Fig.~36
$$
\begin{picture}(0,0)%
\includegraphics{ima36.pstex}%
\end{picture}%
\setlength{\unitlength}{4144sp}%
\begingroup\makeatletter\ifx\SetFigFont\undefined%
\gdef\SetFigFont#1#2#3#4#5{%
  \reset@font\fontsize{#1}{#2pt}%
  \fontfamily{#3}\fontseries{#4}\fontshape{#5}%
  \selectfont}%
\fi\endgroup%
\begin{picture}(2688,1241)(281,-615)
\put(1973,476){\makebox(0,0)[lb]{\smash{\SetFigFont{12}{14.4}{\rmdefault}{\mddefault}{\updefault}$a$}}}
\put(2675,464){\makebox(0,0)[lb]{\smash{\SetFigFont{12}{14.4}{\rmdefault}{\mddefault}{\updefault}$a$}}}
\put(447,491){\makebox(0,0)[lb]{\smash{\SetFigFont{12}{14.4}{\rmdefault}{\mddefault}{\updefault}$a$}}}
\put(1171,480){\makebox(0,0)[lb]{\smash{\SetFigFont{12}{14.4}{\rmdefault}{\mddefault}{\updefault}$a$}}}
\put(508, 78){\makebox(0,0)[lb]{\smash{\SetFigFont{12}{14.4}{\rmdefault}{\mddefault}{\updefault}$p$}}}
\put(1039, 38){\makebox(0,0)[lb]{\smash{\SetFigFont{12}{14.4}{\rmdefault}{\mddefault}{\updefault}$q$}}}
\put(2025, 78){\makebox(0,0)[lb]{\smash{\SetFigFont{12}{14.4}{\rmdefault}{\mddefault}{\updefault}$q$}}}
\put(2581, 63){\makebox(0,0)[lb]{\smash{\SetFigFont{12}{14.4}{\rmdefault}{\mddefault}{\updefault}$p$}}}
\put(1523, 22){\makebox(0,0)[lb]{\smash{\SetFigFont{12}{14.4}{\rmdefault}{\mddefault}{\updefault}or}}}
\end{picture}
$$ 
\begin{center}
{\bf Fig. 36}
\end{center}
enters into the class. Indeed, if one of them did, it would be invariant under
Reidemeister moves of type $II$. Assume that the configuration in the left-hand side of Fig.~37
$$
\begin{picture}(0,0)%
\includegraphics{ima37.pstex}%
\end{picture}%
\setlength{\unitlength}{4144sp}%
\begingroup\makeatletter\ifx\SetFigFont\undefined%
\gdef\SetFigFont#1#2#3#4#5{%
  \reset@font\fontsize{#1}{#2pt}%
  \fontfamily{#3}\fontseries{#4}\fontshape{#5}%
  \selectfont}%
\fi\endgroup%
\begin{picture}(2688,1241)(281,-615)
\put(2675,464){\makebox(0,0)[lb]{\smash{\SetFigFont{12}{14.4}{\rmdefault}{\mddefault}{\updefault}$a$}}}
\put(447,491){\makebox(0,0)[lb]{\smash{\SetFigFont{12}{14.4}{\rmdefault}{\mddefault}{\updefault}$a$}}}
\put(508, 78){\makebox(0,0)[lb]{\smash{\SetFigFont{12}{14.4}{\rmdefault}{\mddefault}{\updefault}$p$}}}
\put(1039, 38){\makebox(0,0)[lb]{\smash{\SetFigFont{12}{14.4}{\rmdefault}{\mddefault}{\updefault}$q$}}}
\put(1171,480){\makebox(0,0)[lb]{\smash{\SetFigFont{12}{14.4}{\rmdefault}{\mddefault}{\updefault}$-a$}}}
\put(1973,476){\makebox(0,0)[lb]{\smash{\SetFigFont{12}{14.4}{\rmdefault}{\mddefault}{\updefault}$-a$}}}
\put(2697,-68){\makebox(0,0)[lb]{\smash{\SetFigFont{12}{14.4}{\rmdefault}{\mddefault}{\updefault}$p$}}}
\put(2186, 15){\makebox(0,0)[lb]{\smash{\SetFigFont{12}{14.4}{\rmdefault}{\mddefault}{\updefault}$q$}}}
\end{picture}
$$
\begin{center}
{\bf Fig. 37}
\end{center}
 enters into the class. The stratum $a^{\mbox{ }}_{\dr}(a|-a, 2a)$ of the
 discriminant forces then the configuration in the right-hand side of Fig.~37
 to enter into the class too. But then, the stratum $a^{+}_{\dr}(a|2a,-a)$
 forces the configuration in the left-hand side of Fig.~38 
$$
\begin{picture}(0,0)%
\includegraphics{ima38.pstex}%
\end{picture}%
\setlength{\unitlength}{4144sp}%
\begingroup\makeatletter\ifx\SetFigFont\undefined%
\gdef\SetFigFont#1#2#3#4#5{%
  \reset@font\fontsize{#1}{#2pt}%
  \fontfamily{#3}\fontseries{#4}\fontshape{#5}%
  \selectfont}%
\fi\endgroup%
\begin{picture}(2688,1241)(281,-615)
\put(447,491){\makebox(0,0)[lb]{\smash{\SetFigFont{12}{14.4}{\rmdefault}{\mddefault}{\updefault}$-a$}}}
\put(508, 78){\makebox(0,0)[lb]{\smash{\SetFigFont{12}{14.4}{\rmdefault}{\mddefault}{\updefault}$q$}}}
\put(1039, 38){\makebox(0,0)[lb]{\smash{\SetFigFont{12}{14.4}{\rmdefault}{\mddefault}{\updefault}$p$}}}
\put(1171,480){\makebox(0,0)[lb]{\smash{\SetFigFont{12}{14.4}{\rmdefault}{\mddefault}{\updefault}$a$}}}
\put(2707,-60){\makebox(0,0)[lb]{\smash{\SetFigFont{12}{14.4}{\rmdefault}{\mddefault}{\updefault}$q$}}}
\put(2161,-101){\makebox(0,0)[lb]{\smash{\SetFigFont{12}{14.4}{\rmdefault}{\mddefault}{\updefault}$p$}}}
\put(2644,472){\makebox(0,0)[lb]{\smash{\SetFigFont{12}{14.4}{\rmdefault}{\mddefault}{\updefault}$-a$}}}
\put(2000,486){\makebox(0,0)[lb]{\smash{\SetFigFont{12}{14.4}{\rmdefault}{\mddefault}{\updefault}$a$}}}
\end{picture}
$$
\begin{center}
{\bf Fig. 38}
\end{center}
to be also in the class. Now, the stratum $a^{+}_{\br}(a|2a,-a)$ forces the
configuration in the right-hand side of Fig.~38 to be in the class too. Thus,
{\em all\/} crossings with marking $-a$ enter into the class. This class is
therefore equal to $W_K(-a)$ independently of $p$ (see [F] for the definitions
of the strata and of $W_K(-a)$). Therefore, the corresponding $T$-invariant
would be $W_K(a).W_K(-a)$, which is not new, and is of course of finite type. $\Box$
\begin{propo}
$$
\makebox(70,50){$ T_K:= T_K \Bigg( D=$}\makebox(70,50){\begin{picture}(0,0)%
\includegraphics{pag52M.pstex}%
\end{picture}%
\setlength{\unitlength}{4144sp}%
\begingroup\makeatletter\ifx\SetFigFont\undefined%
\gdef\SetFigFont#1#2#3#4#5{%
  \reset@font\fontsize{#1}{#2pt}%
  \fontfamily{#3}\fontseries{#4}\fontshape{#5}%
  \selectfont}%
\fi\endgroup%
\begin{picture}(772,947)(494,-258)
\put(893,524){\makebox(0,0)[lb]{\smash{\SetFigFont{12}{14.4}{\rmdefault}{\mddefault}{\updefault}$a$}}}
\end{picture}
 } 
$$
$$
\makebox(70,50){ ; $c_{++}(D)= c_1$, $c_{--}(D)= c_2$, $c_{+-}(D)= c_3$, $c_{-+}(D)=c_4\Bigg)$ }
$$
{\it is the universal $T$-invariant of degree 2 which is not of finite type for $G$-pure global knots.\/}
\end{propo}
(Compare with Definition 3.4.) $T_K$ is "universal" means that any other invariant (of degree 2, not of finite type) can be extracted from $T_K$.
\par {\bf Proof}.
The proposition is an immediate consequence of Theorem~1 and Lemma~3.3. $\Box$
In sect.~9, we show an application of the above invariant. Let $K_1, K_2
\hookrightarrow F^2 \times \mathbb{R}$ be two global knots with respect to the
same non-elliptic vector field $v$ on $F^2$. We assume that $K_i, i= 1, 2$ are
not solid torus knots in $F^2 \times \mathbb{R}$. Let $\mathcal{G}$ be the set
of all possible quotient groups $G$ of $H_1(F^2, \mathbb{Z})$ such that $K_1$
and $K_2$ have global representatives which are $G$-pure. Let $\mathcal{T}$ be
the set of all $T$-invariants of $G$-pure knots with respect to some $G \in \mathcal{G}$.
\par {\bf Conjecture}:
{\it If $K_1$ and $K_2$ are not isotopic, then there are $T$-invariants in $\mathcal{T}$ which distinguish them.\/}
\par {\bf Remarks}:
\begin{enumerate}
\item
Remember that the usual invariants of finite type, in particular the free
homotopy class of the knot, are a subset of $T$-invariants for $m= 0$.
\item
For solid torus knots, the $T$-invariants are nothing but the usual invariants
of finite type (extracted from the generalized Kontsevitch integral). Indeed,
$H_1(S^1 \times I; \mathbb{Z})/\langle [K] \rangle$ is a finite group $G$. In
the Gauss diagram of $K$ with markings $a$, $-a$ in $G$, there are no
subdiagrams of the forms depicted in Fig.~39 at all (because $s$ would be a global knot homologous to $0$).
$$
\begin{picture}(0,0)%
\includegraphics{ima39.pstex}%
\end{picture}%
\setlength{\unitlength}{4144sp}%
\begingroup\makeatletter\ifx\SetFigFont\undefined%
\gdef\SetFigFont#1#2#3#4#5{%
  \reset@font\fontsize{#1}{#2pt}%
  \fontfamily{#3}\fontseries{#4}\fontshape{#5}%
  \selectfont}%
\fi\endgroup%
\begin{picture}(3861,1337)(458,-896)
\put(669,299){\makebox(0,0)[lb]{\smash{\SetFigFont{12}{14.4}{\rmdefault}{\mddefault}{\updefault}$a$}}}
\put(1313,285){\makebox(0,0)[lb]{\smash{\SetFigFont{12}{14.4}{\rmdefault}{\mddefault}{\updefault}$a$}}}
\put(2007,306){\makebox(0,0)[lb]{\smash{\SetFigFont{12}{14.4}{\rmdefault}{\mddefault}{\updefault}$a$}}}
\put(2644,-896){\makebox(0,0)[lb]{\smash{\SetFigFont{12}{14.4}{\rmdefault}{\mddefault}{\updefault}$-a$}}}
\put(3994,296){\makebox(0,0)[lb]{\smash{\SetFigFont{12}{14.4}{\rmdefault}{\mddefault}{\updefault}$-a$}}}
\put(3338,-895){\makebox(0,0)[lb]{\smash{\SetFigFont{12}{14.4}{\rmdefault}{\mddefault}{\updefault}$a$}}}
\put(3618,-290){\makebox(0,0)[lb]{\smash{\SetFigFont{12}{14.4}{\rmdefault}{\mddefault}{\updefault}$S$}}}
\put(2272,-320){\makebox(0,0)[lb]{\smash{\SetFigFont{12}{14.4}{\rmdefault}{\mddefault}{\updefault}$S$}}}
\put(967,-324){\makebox(0,0)[lb]{\smash{\SetFigFont{12}{14.4}{\rmdefault}{\mddefault}{\updefault}$S$}}}
\end{picture}
$$
\begin{center}
{\bf Fig. 39}
\end{center}
Hence, there are no specific $T$-invariants with respect to $G$. Therefore, we
have to consider a quotient $G'$ of $G$. This means that some $a \not= 0 \in
G$ becomes $0$ in $G'$. But, if for a closed braid $K$ and a given $a \in G$,
no crossing is of type $[K^+_p]= a$, then it is easily seen that there always
exist a crossing $p$ with $[K^+_p]= -a \in G$. But  this means that $K$ is
never a $G'$-pure global knot and, hence, there are not any specific $T$-invariants.
\item
Global solid torus knots are closed braids. They are classified by Artin's
theorem together with Garside's solution of the conjugacy problem in braid
groups. Unfortunately, this solution has exponential complexity.
\end{enumerate}
\newpage
\section{$T$-invariants separate ${\mathbb Z}/ 2 {\mathbb Z}$-pure global knots in $T^2 \times {\mathbb R}$}
Let $\{ \alpha, \beta \}$ be generators of $H_1(T^2)$ as shown in Fig.~22. It
is more convenient to use the non-generic vector field $v$ which is tangent to
the fibers of $f$ (see Example~3.4). The difference with a generic vector
field (obtained by a small perturbation of $v$) is that for the latter,
positive multiples of $\beta$ can be represented by global knots. But in any
case, these are solid torus knots and we do not consider them.
\begin{defin}
{\rm A global knot $K \hookrightarrow T^2 \times \mathbb{R}$ with respect to
  $v$ is called a {\it $\mathbb{Z}/2\mathbb{Z}$-pure global knot} if}:
\begin{enumerate}
\item
$H_1(T^2; \mathbb{Z})/\langle [K] \rangle \cong \mathbb{Z}$
\item
$K$ is $\mathbb{Z}/2\mathbb{Z}$-pure for $G= \mathbb{Z}/2\mathbb{Z} \cong (H_1(T^2)/\langle [K] \rangle)/2\mathbb{Z}$
\end{enumerate}
\end{defin}
\par {\bf Remark}.
In particular, condition 1 implies that a $\mathbb{Z}/2\mathbb{Z}$-pure global
knot $K$ is a solid torus knot if and only if $K \hookrightarrow T^2$ (i.e. K
is a torus knot). We show a typical example of a $\mathbb{Z}/2\mathbb{Z}$-pure
global knot in Fig.~40. $[K]= 3\alpha+ \beta$ and $\alpha$ is a generator of
$H_1(T^2)/\langle [K] \rangle$. Notice that switching a crossing $p$ does not
change the marking of $p$ for a $\mathbb{Z}/2\mathbb{Z}$-pure global
knot. Consequently, the property of a global knot to be $\mathbb{Z}/
2\mathbb{Z}$-pure or not depends only on the underlying curve $pr(K)
\hookrightarrow T^2$.
$$
\begin{picture}(0,0)%
\includegraphics{im14.pstex}%
\end{picture}%
\setlength{\unitlength}{4144sp}%
\begingroup\makeatletter\ifx\SetFigFont\undefined%
\gdef\SetFigFont#1#2#3#4#5{%
  \reset@font\fontsize{#1}{#2pt}%
  \fontfamily{#3}\fontseries{#4}\fontshape{#5}%
  \selectfont}%
\fi\endgroup%
\begin{picture}(4942,2884)(476,-2500)
\put(1919,-1309){\makebox(0,0)[lb]{\smash{\SetFigFont{10}{12.0}{\rmdefault}{\mddefault}{\updefault}$2\alpha + \beta$}}}
\put(1353,-2029){\makebox(0,0)[lb]{\smash{\SetFigFont{10}{12.0}{\rmdefault}{\mddefault}{\updefault}$2\alpha + \beta$}}}
\put(4116,-2455){\makebox(0,0)[lb]{\smash{\SetFigFont{10}{12.0}{\rmdefault}{\mddefault}{\updefault}$2 \alpha + \beta$}}}
\put(3216,-2129){\makebox(0,0)[lb]{\smash{\SetFigFont{10}{12.0}{\rmdefault}{\mddefault}{\updefault}$2 \alpha + \beta$}}}
\put(2477,-1352){\makebox(0,0)[lb]{\smash{\SetFigFont{10}{12.0}{\rmdefault}{\mddefault}{\updefault}$2\alpha + \beta$}}}
\put(2389,-1659){\makebox(0,0)[lb]{\smash{\SetFigFont{10}{12.0}{\rmdefault}{\mddefault}{\updefault}$\alpha$}}}
\put(2873,-1701){\makebox(0,0)[lb]{\smash{\SetFigFont{10}{12.0}{\rmdefault}{\mddefault}{\updefault}$\alpha$}}}
\put(3052,-1958){\makebox(0,0)[lb]{\smash{\SetFigFont{10}{12.0}{\rmdefault}{\mddefault}{\updefault}$\alpha$}}}
\put(4070,-2169){\makebox(0,0)[lb]{\smash{\SetFigFont{10}{12.0}{\rmdefault}{\mddefault}{\updefault}$\alpha$}}}
\put(1960,-1645){\makebox(0,0)[lb]{\smash{\SetFigFont{10}{12.0}{\rmdefault}{\mddefault}{\updefault}$\alpha$}}}
\end{picture}
$$
\begin{center}
{\bf Fig. 40}
\end{center}
\begin{theor}
{\it Let $K \hookrightarrow T^2 \times \mathbb{R}$ be a $\mathbb{Z}/
  2\mathbb{Z}$- pure global knot with $c$ crossings. Let $d$ be any natural
  number not bigger than $c$. Then, the knot $K$ is uniquely determined by the
  set of all $T$-invariants $T_K$ of finite type, of the degrees $(m= d, n=
  0)$ with respect to $G= \mathbb{Z}/2 \mathbb{Z}$.\/}
\end{theor}
\par {\bf Remark}.
The number of such $T$-invariants is finite. Thus, Theorem~2 proves the
conjecture in sect.~3, in the case of $\mathbb{Z}/2\mathbb{Z}$-pure global
knots in $T^2 \times \mathbb{R}$ and, moreover gives an effective solution to
the problem. Notice that we do not need here the $T$-invariants of infinite type.
\par {\bf Proof}.
The proof consists of two steps.
\par {\sl Step 1:\/} The Gauss diagram with markings in $G \cong \mathbb{Z}/ 2\mathbb{Z}$ determines $K$.
\par {\sl Step 2:\/} The invariants $T_K$ determine the Gauss diagram of $K$.
\par {\sl Step~1:\/}
For local knots, it is well known and rather evident. But it is not at all
obvious for knots in $T^2 \times \mathbb{R}$. By definition, the marking of
each crossing of a $\mathbb{Z}/2 \mathbb{Z}$-pure knot is the non-trivial
element in $\mathbb{Z}/ 2\mathbb{Z}$. Therefore, we do not write it in Gauss
diagrams, configurations, etc \dots The $T$-invariant of degree $(d_1, d_2)=
(0, 0)$ is the free homotopy class of $K$, or (equivalently here), the
homology class represented by $K$ (which is a primitive class because $K$ is
not a solid torus knot). Consequently, we have to show that the Gauss diagram determines $K$ in its given homology class.
\begin{defin}
{\rm A set of arrows in a Gauss diagram is called a {\it bunch of arrows\/}
if their number is even and:}
\end{defin}
\begin{enumerate}
\item
they are near to each other (i.e. there are small arcs on $S^1$ between them
where no other arrow starts or ends)
\item
each two arrows cut in exactly one point
\item
the orientation of the arrows is alternating
\item
all the arrows have the same writhe
\end{enumerate}
 (See an example in Fig.~41.)
$$
\begin{picture}(0,0)%
\includegraphics{ima41.pstex}%
\end{picture}%
\setlength{\unitlength}{4144sp}%
\begingroup\makeatletter\ifx\SetFigFont\undefined%
\gdef\SetFigFont#1#2#3#4#5{%
  \reset@font\fontsize{#1}{#2pt}%
  \fontfamily{#3}\fontseries{#4}\fontshape{#5}%
  \selectfont}%
\fi\endgroup%
\begin{picture}(4983,2792)(446,-2193)
\put(1230,-2045){\makebox(0,0)[lb]{\smash{\SetFigFont{12}{14.4}{\rmdefault}{\mddefault}{\updefault}$\emptyset$}}}
\put(3760,310){\makebox(0,0)[lb]{\smash{\SetFigFont{12}{14.4}{\rmdefault}{\mddefault}{\updefault}$\emptyset$}}}
\put(3760,-2025){\makebox(0,0)[lb]{\smash{\SetFigFont{12}{14.4}{\rmdefault}{\mddefault}{\updefault}$\emptyset$}}}
\put(4235,-2135){\makebox(0,0)[lb]{\smash{\SetFigFont{12}{14.4}{\rmdefault}{\mddefault}{\updefault}$\emptyset$}}}
\put(4710,-1960){\makebox(0,0)[lb]{\smash{\SetFigFont{12}{14.4}{\rmdefault}{\mddefault}{\updefault}$\emptyset$}}}
\put(1577,330){\makebox(0,0)[lb]{\smash{\SetFigFont{12}{14.4}{\rmdefault}{\mddefault}{\updefault}$\emptyset$}}}
\put(4253,443){\makebox(0,0)[lb]{\smash{\SetFigFont{12}{14.4}{\rmdefault}{\mddefault}{\updefault}$\emptyset$}}}
\put(4751,220){\makebox(0,0)[lb]{\smash{\SetFigFont{12}{14.4}{\rmdefault}{\mddefault}{\updefault}$\emptyset$}}}
\put(2441,-75){\makebox(0,0)[lb]{\smash{\SetFigFont{12}{14.4}{\rmdefault}{\mddefault}{\updefault}$s_1$}}}
\end{picture}
$$
\begin{center}
{\bf Fig. 41}
\end{center}
\begin{lem}
{\it After possibly performing Reidemeister moves of type $II$, such that each
 of them decreases the number of crossings (see left-hand part of Fig.~42), 
$$
\begin{picture}(0,0)%
\includegraphics{ima42.pstex}%
\end{picture}%
\setlength{\unitlength}{4144sp}%
\begingroup\makeatletter\ifx\SetFigFont\undefined%
\gdef\SetFigFont#1#2#3#4#5{%
  \reset@font\fontsize{#1}{#2pt}%
  \fontfamily{#3}\fontseries{#4}\fontshape{#5}%
  \selectfont}%
\fi\endgroup%
\begin{picture}(4325,446)(85,-369)
\end{picture}
$$
\begin{center}
{\bf Fig. 42}
\end{center}
the Gauss diagram of a $\mathbb{Z}/2\mathbb{Z}$-pure global knot is of the
following form: some chord diagram where each chord is replaced by some bunch
of arrows. Moreover, there exists a homotopy from $K$ to a torus knot $K'$
which is an isotopy of diagrams besides possibly performing transformations of the type depicted in the right-hand part of Fig.~42.\/}
\end{lem}
(Of course, $K'$ is determined by the homology class represented by $K$.)
\par{\it Example}
The knot in Fig.~40 corresponds to the diagram in the right-hand part of Fig.~43.
$$
\begin{picture}(0,0)%
\includegraphics{ima43.pstex}%
\end{picture}%
\setlength{\unitlength}{4144sp}%
\begingroup\makeatletter\ifx\SetFigFont\undefined%
\gdef\SetFigFont#1#2#3#4#5{%
  \reset@font\fontsize{#1}{#2pt}%
  \fontfamily{#3}\fontseries{#4}\fontshape{#5}%
  \selectfont}%
\fi\endgroup%
\begin{picture}(3836,1578)(310,-787)
\put(2718,245){\makebox(0,0)[lb]{\smash{\SetFigFont{10}{12.0}{\rmdefault}{\mddefault}{\updefault}$+$}}}
\put(3012,437){\makebox(0,0)[lb]{\smash{\SetFigFont{10}{12.0}{\rmdefault}{\mddefault}{\updefault}$+$}}}
\put(3408,377){\makebox(0,0)[lb]{\smash{\SetFigFont{10}{12.0}{\rmdefault}{\mddefault}{\updefault}$+$}}}
\put(3744,362){\makebox(0,0)[lb]{\smash{\SetFigFont{10}{12.0}{\rmdefault}{\mddefault}{\updefault}$+$}}}
\put(2853,-283){\makebox(0,0)[lb]{\smash{\SetFigFont{10}{12.0}{\rmdefault}{\mddefault}{\updefault}$-$}}}
\put(2958,-121){\makebox(0,0)[lb]{\smash{\SetFigFont{10}{12.0}{\rmdefault}{\mddefault}{\updefault}$-$}}}
\put(3456,-166){\makebox(0,0)[lb]{\smash{\SetFigFont{12}{14.4}{\rmdefault}{\mddefault}{\updefault}$+$}}}
\put(3711,-136){\makebox(0,0)[lb]{\smash{\SetFigFont{10}{12.0}{\rmdefault}{\mddefault}{\updefault}$+$}}}
\put(3735,-325){\makebox(0,0)[lb]{\smash{\SetFigFont{10}{12.0}{\rmdefault}{\mddefault}{\updefault}$+$}}}
\put(3540,-517){\makebox(0,0)[lb]{\smash{\SetFigFont{10}{12.0}{\rmdefault}{\mddefault}{\updefault}$+$}}}
\end{picture}
$$
\begin{center}
{\bf Fig. 43}
\end{center}

Notice that e.g. \makebox(40,10){\begin{picture}(0,0)%
\includegraphics{pag63.pstex}%
\end{picture}%
\setlength{\unitlength}{4144sp}%
\begingroup\makeatletter\ifx\SetFigFont\undefined%
\gdef\SetFigFont#1#2#3#4#5{%
  \reset@font\fontsize{#1}{#2pt}%
  \fontfamily{#3}\fontseries{#4}\fontshape{#5}%
  \selectfont}%
\fi\endgroup%
\begin{picture}(290,333)(445,68)
\end{picture}
} can never appear as the underlying chord diagram of a $\mathbb{Z}/2\mathbb{Z}$-pure global knot. (This is an exercise).
\par {\bf Proof of Lemma $4.1$}
We start by eliminating all couples of arrows corresponding to crossings say $q_1$ and $q_2$ as shown in Fig.~44.
$$
\begin{picture}(0,0)%
\includegraphics{ima44.pstex}%
\end{picture}%
\setlength{\unitlength}{4144sp}%
\begingroup\makeatletter\ifx\SetFigFont\undefined%
\gdef\SetFigFont#1#2#3#4#5{%
  \reset@font\fontsize{#1}{#2pt}%
  \fontfamily{#3}\fontseries{#4}\fontshape{#5}%
  \selectfont}%
\fi\endgroup%
\begin{picture}(5257,2606)(458,-2228)
\put(1590,-2228){\makebox(0,0)[lb]{\smash{\SetFigFont{12}{14.4}{\rmdefault}{\mddefault}{\updefault}$I_1$}}}
\put(1960,194){\makebox(0,0)[lb]{\smash{\SetFigFont{12}{14.4}{\rmdefault}{\mddefault}{\updefault}$q_1$}}}
\put(1165,222){\makebox(0,0)[lb]{\smash{\SetFigFont{12}{14.4}{\rmdefault}{\mddefault}{\updefault}$q_2$}}}
\put(1450,-60){\makebox(0,0)[lb]{\smash{\SetFigFont{12}{14.4}{\rmdefault}{\mddefault}{\updefault}$I_2$}}}
\put(1184,-395){\makebox(0,0)[lb]{\smash{\SetFigFont{12}{14.4}{\rmdefault}{\mddefault}{\updefault}$+$}}}
\put(1733,-373){\makebox(0,0)[lb]{\smash{\SetFigFont{12}{14.4}{\rmdefault}{\mddefault}{\updefault}$-$}}}
\end{picture}
$$
\begin{center}
{\bf Fig. 44}
\end{center}
We assume that the arcs $I_1$ and $I_2$ are empty. We have to prove that both
arcs $I_1$ and $I_2$ are small or equivalently that $[K^+_{q_1}]= [K^+_{q_2}]
\in H_1(T^2; \mathbb{Z})$. Then the above operation corresponds to a
Reidemeister move inverse to the one depicted in the left-hand part of
Fig.~42. After possibly performing a (global) isotopy of the diagram, we may
assume that at least one of the two arcs, say $I_1$, is small (see Fig.~45).
$$
\begin{picture}(0,0)%
\includegraphics{im17a.pstex}%
\end{picture}%
\setlength{\unitlength}{4144sp}%
\begingroup\makeatletter\ifx\SetFigFont\undefined%
\gdef\SetFigFont#1#2#3#4#5{%
  \reset@font\fontsize{#1}{#2pt}%
  \fontfamily{#3}\fontseries{#4}\fontshape{#5}%
  \selectfont}%
\fi\endgroup%
\begin{picture}(4942,3698)(476,-3314)
\put(2716,-2011){\makebox(0,0)[rb]{\smash{\SetFigFont{12}{14.4}{\rmdefault}{\mddefault}{\updefault}$q_2$}}}
\put(2446,-1321){\makebox(0,0)[lb]{\smash{\SetFigFont{12}{14.4}{\rmdefault}{\mddefault}{\updefault}$q_1$}}}
\end{picture}
$$
$$
\begin{picture}(0,0)%
\includegraphics{im17.pstex}%
\end{picture}%
\setlength{\unitlength}{4144sp}%
\begingroup\makeatletter\ifx\SetFigFont\undefined%
\gdef\SetFigFont#1#2#3#4#5{%
  \reset@font\fontsize{#1}{#2pt}%
  \fontfamily{#3}\fontseries{#4}\fontshape{#5}%
  \selectfont}%
\fi\endgroup%
\begin{picture}(5376,2711)(210,-6437)
\put(2723,-5625){\makebox(0,0)[lb]{\smash{\SetFigFont{12}{14.4}{\rmdefault}{\mddefault}{\updefault}$I_1$}}}
\put(2483,-5895){\makebox(0,0)[lb]{\smash{\SetFigFont{12}{14.4}{\rmdefault}{\mddefault}{\updefault}$q_1$}}}
\put(3796,-5903){\makebox(0,0)[rb]{\smash{\SetFigFont{12}{14.4}{\rmdefault}{\mddefault}{\updefault}$q_2$}}}
\end{picture}
$$
\begin{center}
{\bf Fig. 45}
\end{center}
If the arc $I_2$ cannot be made small at the same time as $I_1$, then we have
exactly the situation shown in Fig.~45. But then $w(q_1)= w(q_2)$ in contradiction to our assumption on $q_1$ and $q_2$.
We observe now that for a diagram of a $\mathbb{Z}/2\mathbb{Z}$-pure global
knot, there are no possible Reidemeister moves of type $III$ at all. Indeed, in $pr(K)$, there cannot be any triangle 
$$
\begin{picture}(0,0)%
\includegraphics{ima46.pstex}%
\end{picture}%
\setlength{\unitlength}{4144sp}%
\begingroup\makeatletter\ifx\SetFigFont\undefined%
\gdef\SetFigFont#1#2#3#4#5{%
  \reset@font\fontsize{#1}{#2pt}%
  \fontfamily{#3}\fontseries{#4}\fontshape{#5}%
  \selectfont}%
\fi\endgroup%
\begin{picture}(1764,752)(344,-1156)
\end{picture}
$$
\begin{center}
{\bf Fig. 46}
\end{center}
as shown in Fig.~46 because the three crossings $q_i, i=1, 2, 3$ always verify
a relation: $[K^+_{q_3}] =  [K^+_{q_1}] + [K^+_{q_2}]  mod [K]$.
Consequently, they could not be all three non-zero in $G \cong
\mathbb{Z}/2\mathbb{Z}$. After having performed all possible Reidemeister
moves of type $II$ as in the left-hand side of Fig.~42, we obtain a diagram of
$K$ which we call {\em minimal\/}. It is characterized by the fact that it does
not allow any Reidemeister moves, except for those which strictly increase the
number of crossings (i.e. the move inverse to the one in the left-hand side of Fig.~42).
\par {\bf Claim $1$}.
If the minimal diagram is not empty, then it contains always a {\em 2-gon\/} (see the left-hand part of Fig.~47).
$$
\begin{picture}(0,0)%
\includegraphics{ima47.pstex}%
\end{picture}%
\setlength{\unitlength}{4144sp}%
\begingroup\makeatletter\ifx\SetFigFont\undefined%
\gdef\SetFigFont#1#2#3#4#5{%
  \reset@font\fontsize{#1}{#2pt}%
  \fontfamily{#3}\fontseries{#4}\fontshape{#5}%
  \selectfont}%
\fi\endgroup%
\begin{picture}(4946,2685)(446,-1928)
\put(3895,-68){\makebox(0,0)[lb]{\smash{\SetFigFont{12}{14.4}{\rmdefault}{\mddefault}{\updefault}$q_1$}}}
\put(3580,-1023){\makebox(0,0)[lb]{\smash{\SetFigFont{12}{14.4}{\rmdefault}{\mddefault}{\updefault}$q_2$}}}
\put(4290,-1928){\makebox(0,0)[lb]{\smash{\SetFigFont{12}{14.4}{\rmdefault}{\mddefault}{\updefault}$I_2$}}}
\put(4248,562){\makebox(0,0)[lb]{\smash{\SetFigFont{12}{14.4}{\rmdefault}{\mddefault}{\updefault}$I_1$}}}
\end{picture}
$$
\begin{center}
{\bf Fig. 47}
\end{center}
Moreover, the following fragment (Fig.~47, right-hand side) of the minimal
Gauss diagram of a $\mathbb{Z}/2\mathbb{Z}$-pure knot corresponds always to a
2-gon. Here, $w(q_1)= w(q_2)$ and $I_i, i=1, 2$ are both empty.
\par {\bf Proof of the claim}. 
We start with the following observation: Let $D(K)$ be a knot which is
obtained from $K$ after performing a Dehn twist of $T^2$. We do not change the
vector field $v$. If the Dehn twist is along $\beta$, or positive along
$\alpha$, then $D(K)$ is still a global ($\mathbb{Z}/ 2\mathbb{Z}$-pure) knot
with respect to $v$ (see Fig.~22). By definition, the two sides of a 2-gon
form a loop which is homotopic to $0$ in $T^2$. Let us consider first "fake"
2-gons, i.e. 2-gons in $pr(K) \subset T^2$ such that the corresponding loop is
not homotopic to $0$ in $T^2$. Using the above observation, we can restrict our
considerations exactly to the two cases shown in Fig.~48.
$$
\begin{picture}(0,0)%
\includegraphics{im19.pstex}%
\end{picture}%
\setlength{\unitlength}{4144sp}%
\begingroup\makeatletter\ifx\SetFigFont\undefined%
\gdef\SetFigFont#1#2#3#4#5{%
  \reset@font\fontsize{#1}{#2pt}%
  \fontfamily{#3}\fontseries{#4}\fontshape{#5}%
  \selectfont}%
\fi\endgroup%
\begin{picture}(3155,3979)(436,-3491)
\put(2078,-786){\makebox(0,0)[b]{\smash{\SetFigFont{12}{14.4}{\rmdefault}{\mddefault}{\updefault}$I_1$}}}
\put(1772,-192){\makebox(0,0)[b]{\smash{\SetFigFont{12}{14.4}{\rmdefault}{\mddefault}{\updefault}$I_2$}}}
\put(2841,-2638){\makebox(0,0)[lb]{\smash{\SetFigFont{12}{14.4}{\rmdefault}{\mddefault}{\updefault}$I_2$}}}
\put(2426,-3433){\makebox(0,0)[lb]{\smash{\SetFigFont{12}{14.4}{\rmdefault}{\mddefault}{\updefault}$\alpha$}}}
\put(3591,-597){\makebox(0,0)[lb]{\smash{\SetFigFont{12}{14.4}{\rmdefault}{\mddefault}{\updefault}$\beta$}}}
\put(2116,-2558){\makebox(0,0)[lb]{\smash{\SetFigFont{12}{14.4}{\rmdefault}{\mddefault}{\updefault}$I_1$}}}
\put(436,323){\makebox(0,0)[lb]{\smash{\SetFigFont{12}{14.4}{\rmdefault}{\mddefault}{\updefault}Case $1$}}}
\put(466,-1274){\makebox(0,0)[lb]{\smash{\SetFigFont{12}{14.4}{\rmdefault}{\mddefault}{\updefault}Case $2$}}}
\end{picture}
$$
\begin{center}
{\bf Fig.48}
\end{center}
$I_1 \cup I_2$ does not cut the rest of the knot $K$. Consequently, in the
second case, $K$ is a solid torus knot which is not a torus knot (the minimal
diagram is not empty). We do not consider these knots. In the first case, $K$
represents $2\alpha + x\beta, x \in \mathbb{Z}$ and $x$ odd. One easily checks
that in the minimal diagram of $K$ there is always some 2-gon. An example is shown in Fig.~49.\newpage
$$
\begin{picture}(0,0)%
\includegraphics{im20.pstex}%
\end{picture}%
\setlength{\unitlength}{4144sp}%
\begingroup\makeatletter\ifx\SetFigFont\undefined%
\gdef\SetFigFont#1#2#3#4#5{%
  \reset@font\fontsize{#1}{#2pt}%
  \fontfamily{#3}\fontseries{#4}\fontshape{#5}%
  \selectfont}%
\fi\endgroup%
\begin{picture}(3231,1567)(585,-1423)
\put(2624,-1187){\makebox(0,0)[lb]{\smash{\SetFigFont{12}{14.4}{\rmdefault}{\mddefault}{\updefault}$q_2$}}}
\put(2165,-1191){\makebox(0,0)[lb]{\smash{\SetFigFont{12}{14.4}{\rmdefault}{\mddefault}{\updefault}$q_1$}}}
\put(1304,-703){\makebox(0,0)[lb]{\smash{\SetFigFont{12}{14.4}{\rmdefault}{\mddefault}{\updefault}$I_2$}}}
\put(1026,-1104){\makebox(0,0)[lb]{\smash{\SetFigFont{12}{14.4}{\rmdefault}{\mddefault}{\updefault}$I_1$}}}
\end{picture}
$$\nopagebreak\nolinebreak 
\hspace{-30mm}\begin{center}
{\bf Fig. 49}
\end{center}
If we change in the diagram in Fig.~47 exactly one of the crossings $q_1$ or
$q_2$ to its inverse, then we obtain a fragment as shown in Fig.~44.
We have already proven that the fragment in Fig.~44 corresponds to 
$$
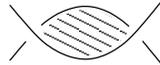
\begin{picture}(0,0)%
\includegraphics{page69.pstex}%
\end{picture}%
\setlength{\unitlength}{4144sp}%
\begingroup\makeatletter\ifx\SetFigFont\undefined%
\gdef\SetFigFont#1#2#3#4#5{%
  \reset@font\fontsize{#1}{#2pt}%
  \fontfamily{#3}\fontseries{#4}\fontshape{#5}%
  \selectfont}%
\fi\endgroup%
\begin{picture}(927,389)(85,-312)
\end{picture}
$$
Consequently, the fragment in Fig.~47 corresponds to 
$$
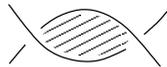
\begin{picture}(0,0)%
\includegraphics{page69bis.pstex}%
\end{picture}%
\setlength{\unitlength}{4144sp}%
\begingroup\makeatletter\ifx\SetFigFont\undefined%
\gdef\SetFigFont#1#2#3#4#5{%
  \reset@font\fontsize{#1}{#2pt}%
  \fontfamily{#3}\fontseries{#4}\fontshape{#5}%
  \selectfont}%
\fi\endgroup%
\begin{picture}(990,395)(418,46)
\end{picture}
$$
 which is a 2-gon.
\par Assume now that the diagram of $K$ (in general position) contains no
2-gons at all. We have already proven that the diagram of $K$ contains no triangles whose three sides form a loop which is contractible in $T^2$.
\par {\bf Sub-claim}. {\it The sides of each $n$-gon in $pr(K) \subset T^2$
  form a contractible loop if $n \geq 3$.}
\par Indeed, either we are in the situation analogue to case 2 in Fig.~48, and
hence, $K$ is a solid torus knot, or we are in the situation analogue to case
1in Fig.~48. But this is not possible if $n \geq 3$ as shown in Fig.~50: 
$$
\begin{picture}(0,0)%
\includegraphics{im21.pstex}%
\end{picture}%
\setlength{\unitlength}{4144sp}%
\begingroup\makeatletter\ifx\SetFigFont\undefined%
\gdef\SetFigFont#1#2#3#4#5{%
  \reset@font\fontsize{#1}{#2pt}%
  \fontfamily{#3}\fontseries{#4}\fontshape{#5}%
  \selectfont}%
\fi\endgroup%
\begin{picture}(2762,1011)(699,-1244)
\put(1694,-797){\makebox(0,0)[lb]{\smash{\SetFigFont{12}{14.4}{\rmdefault}{\mddefault}{\updefault}$q_2$}}}
\put(1544,-1119){\makebox(0,0)[lb]{\smash{\SetFigFont{12}{14.4}{\rmdefault}{\mddefault}{\updefault}$q_1$}}}
\put(2732,-883){\makebox(0,0)[lb]{\smash{\SetFigFont{12}{14.4}{\rmdefault}{\mddefault}{\updefault}$q_3$}}}
\put(3445,-429){\makebox(0,0)[lb]{\smash{\SetFigFont{12}{14.4}{\rmdefault}{\mddefault}{\updefault}$\beta$}}}
\put(1821,-455){\makebox(0,0)[lb]{\smash{\SetFigFont{14}{16.8}{\rmdefault}{\mddefault}{\updefault}?}}}
\end{picture}
$$
\begin{center}
{\bf Fig. 50}
\end{center}
The branch of $K$ through $q_1$ and

$q_2$ is blocked. This proves the subclaim.

The assuption (no 2-gons) together with the subclaim imply that the

4-valent graph $pr(K) \subset T^2$ splits $T^2$ into contractible

4-gons, 5-gons \dots Let $v_0$ be the number of vertices and $v_1$

be the number of edges of $pr(K)$. Let $v_2$ be the number of

components of $T^2 \setminus pr(K)$. Evidently, $v_1= 2v_0$.

One has: $v_0- v_1+ v_2= \chi(T^2)= 0$ and hence, $v_0= v_2$.

We denote by $\sharp(.)$ the number of (.).

$\sharp$(angles)$= 4v_0= 4v_2$. On the other hand, $\sharp$(angles)$= 4\sharp$(4-gons)$+ 5\sharp$(5-gons)$+ \dots$

This implies that $0= \sharp$(5-gons)= $\sharp$(6-gons)$= \dots$

Consequently, $T^2 \setminus pr(K)$ consists only of contractible

4-gons.
We take one of them. It has at least two opposite sides which

have the same orientation (induced by the orientation of $K$). We add

to the 4-gon the two neighbouring 4-gons corresponding to the remaining

two sides (see Fig.51).
$$
\begin{picture}(0,0)%
\includegraphics{im22.pstex}%
\end{picture}%
\setlength{\unitlength}{4144sp}%
\begingroup\makeatletter\ifx\SetFigFont\undefined%
\gdef\SetFigFont#1#2#3#4#5{%
  \reset@font\fontsize{#1}{#2pt}%
  \fontfamily{#3}\fontseries{#4}\fontshape{#5}%
  \selectfont}%
\fi\endgroup%
\begin{picture}(3177,1003)(438,-1281)
\put(1777,-804){\makebox(0,0)[lb]{\smash{\SetFigFont{12}{14.4}{\rmdefault}{\mddefault}{\updefault}$4-gon$}}}
\end{picture}
$$ 
\begin{center}
{\bf Fig. 51}
\end{center}
We continue the process and at the end, we obtain an orientable immersed band
in $T^2$. But the boundary of the band has two components contradicting the
fact that $K$ is a knot. This proves Claim~1.
\par We take now the existing 2-gon and make a homotopy as indicated in
Fig.~52. $$
\begin{picture}(0,0)%
\includegraphics{ima52.pstex}%
\end{picture}%
\setlength{\unitlength}{4144sp}%
\begingroup\makeatletter\ifx\SetFigFont\undefined%
\gdef\SetFigFont#1#2#3#4#5{%
  \reset@font\fontsize{#1}{#2pt}%
  \fontfamily{#3}\fontseries{#4}\fontshape{#5}%
  \selectfont}%
\fi\endgroup%
\begin{picture}(4779,1392)(446,-1493)
\put(1246,-1298){\makebox(0,0)[lb]{\smash{\SetFigFont{12}{14.4}{\rmdefault}{\mddefault}{\updefault}$K$}}}
\put(3916,-1493){\makebox(0,0)[lb]{\smash{\SetFigFont{12}{14.4}{\rmdefault}{\mddefault}{\updefault}$K'$}}}
\end{picture}
 $$\begin{center}
 {\bf Fig. 52}
 \end{center}
$K'$ is again a $\mathbb{Z}/2\mathbb{Z}$-pure global knot.
\par {\bf Claim~$2$}
{\it If $K$ was already minimal, then $K'$ is minimal too}.
\par {\bf Proof of Claim~2}
If $K'$ is not minimal, then it contains a fragment as shown in Fig.~53
$$
 \begin{picture}(0,0)%
 \includegraphics{im23.pstex}%
 \end{picture}%
 \setlength{\unitlength}{4144sp}%
 \begingroup\makeatletter\ifx\SetFigFont\undefined%
 \gdef\SetFigFont#1#2#3#4#5{%
   \reset@font\fontsize{#1}{#2pt}%
   \fontfamily{#3}\fontseries{#4}\fontshape{#5}%
   \selectfont}%
 \fi\endgroup%
 \begin{picture}(3215,926)(621,-1186)
 \put(2016,-695){\makebox(0,0)[lb]{\smash{\SetFigFont{17}{20.4}{\rmdefault}{\mddefault}{\updefault}or}}}
 \end{picture}
 $$\begin{center}
 {\bf Fig. 53}
 \end{center}
 We have already proven that we can transform $K$ into a torus knot by

 performing only the operations 
$$
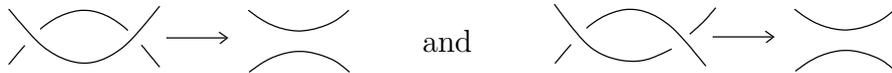
\begin{picture}(0,0)%
\includegraphics{page74.pstex}%
\end{picture}%
\setlength{\unitlength}{4144sp}%
\begingroup\makeatletter\ifx\SetFigFont\undefined%
\gdef\SetFigFont#1#2#3#4#5{%
  \reset@font\fontsize{#1}{#2pt}%
  \fontfamily{#3}\fontseries{#4}\fontshape{#5}%
  \selectfont}%
\fi\endgroup%
\begin{picture}(5331,431)(85,-342)
\put(2574,-195){\makebox(0,0)[lb]{\smash{\SetFigFont{12}{14.4}{\rmdefault}{\mddefault}{\updefault}and}}}
\end{picture}
$$
on the diagram of $K$. Therefore, we

may assume that we have eliminated all crossings of $K$ outside of the

above fragment. Again, by using appropriate Dehn twists, we can reduce

our considerations to the two cases shown in Fig.~54. (We need only

$pr(K) \subset T^2$.)

$$
\begin{picture}(0,0)%
\includegraphics{im24.pstex}%
\end{picture}%
\setlength{\unitlength}{4144sp}%
\begingroup\makeatletter\ifx\SetFigFont\undefined%
\gdef\SetFigFont#1#2#3#4#5{%
  \reset@font\fontsize{#1}{#2pt}%
  \fontfamily{#3}\fontseries{#4}\fontshape{#5}%
  \selectfont}%
\fi\endgroup%
\begin{picture}(3691,2176)(466,-1861)
\put(2141,-1164){\makebox(0,0)[lb]{\smash{\SetFigFont{12}{14.4}{\rmdefault}{\mddefault}{\updefault}$q$}}}
\put(466,180){\makebox(0,0)[lb]{\smash{\SetFigFont{12}{14.4}{\rmdefault}{\mddefault}{\updefault}Case $1$}}}
\end{picture}
$$
$$
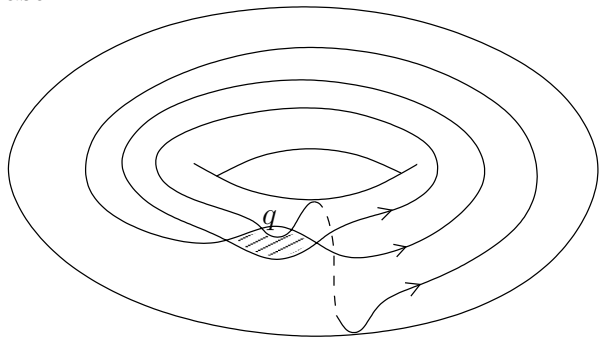
\begin{picture}(0,0)%
\includegraphics{im24bis.pstex}%
\end{picture}%
\setlength{\unitlength}{4144sp}%
\begingroup\makeatletter\ifx\SetFigFont\undefined%
\gdef\SetFigFont#1#2#3#4#5{%
  \reset@font\fontsize{#1}{#2pt}%
  \fontfamily{#3}\fontseries{#4}\fontshape{#5}%
  \selectfont}%
\fi\endgroup%
\begin{picture}(3651,2266)(443,-1936)
\put(1848,-1209){\makebox(0,0)[lb]{\smash{\SetFigFont{12}{14.4}{\rmdefault}{\mddefault}{\updefault}$q$}}}
\put(443,195){\makebox(0,0)[lb]{\smash{\SetFigFont{12}{14.4}{\rmdefault}{\mddefault}{\updefault}Case $2$}}}
\end{picture}
$$
\begin{center}
{\bf Fig. 54}
\end{center}

In case 1, we have $[K^+_q]= 2\alpha$ mod $[K]$ and in case 2, we have

$[K^+_q]= 2\alpha+ 2\beta$ mod $[K]$. Consequently, $[K^+_q]= 0 \in G$

and $K$ was not $\mathbb{Z}/2\mathbb{Z}$-pure. This proves Claim~2.

By Claim~2, we can reduce the minimal diagram of $K$ to a torus knot

by using {\it only\/} the operation 
$$
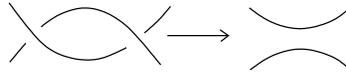
\begin{picture}(0,0)%
\includegraphics{page75.pstex}%
\end{picture}%
\setlength{\unitlength}{4144sp}%
\begingroup\makeatletter\ifx\SetFigFont\undefined%
\gdef\SetFigFont#1#2#3#4#5{%
  \reset@font\fontsize{#1}{#2pt}%
  \fontfamily{#3}\fontseries{#4}\fontshape{#5}%
  \selectfont}%
\fi\endgroup%
\begin{picture}(2067,425)(431,-133)
\end{picture}
$$
Moreover, we have proven that

this operation corresponds exactly to the operation in Fig.~55.
$$
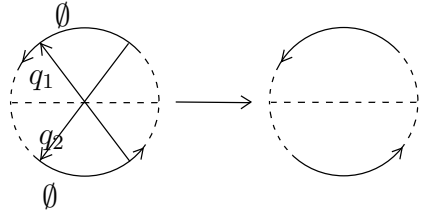
\begin{picture}(0,0)%
\includegraphics{im25.pstex}%
\end{picture}%
\setlength{\unitlength}{4144sp}%
\begingroup\makeatletter\ifx\SetFigFont\undefined%
\gdef\SetFigFont#1#2#3#4#5{%
  \reset@font\fontsize{#1}{#2pt}%
  \fontfamily{#3}\fontseries{#4}\fontshape{#5}%
  \selectfont}%
\fi\endgroup%
\begin{picture}(2543,1673)(383,-2737)
\put(567,-1564){\makebox(0,0)[lb]{\smash{\SetFigFont{12}{14.4}{\rmdefault}{\mddefault}{\updefault}$q_1$}}}
\put(627,-1925){\makebox(0,0)[lb]{\smash{\SetFigFont{12}{14.4}{\rmdefault}{\mddefault}{\updefault}$q_2$}}}
\put(721,-1220){\makebox(0,0)[lb]{\smash{\SetFigFont{12}{14.4}{\rmdefault}{\mddefault}{\updefault}$\emptyset$}}}
\put(646,-2303){\makebox(0,0)[lb]{\smash{\SetFigFont{12}{14.4}{\rmdefault}{\mddefault}{\updefault}$\emptyset$}}}
\put(383,-2679){\makebox(0,0)[lb]{\smash{\SetFigFont{12}{14.4}{\rmdefault}{\mddefault}{\updefault}$w(q_1)=w(q_2)$}}}
\end{picture}
$$
\begin{center}
{\bf Fig. 55}
\end{center}
Suppose, that we have a fragment as shown in Fig.~56.

$$
\begin{picture}(0,0)%
\includegraphics{im26.pstex}%
\end{picture}%
\setlength{\unitlength}{4144sp}%
\begingroup\makeatletter\ifx\SetFigFont\undefined%
\gdef\SetFigFont#1#2#3#4#5{%
  \reset@font\fontsize{#1}{#2pt}%
  \fontfamily{#3}\fontseries{#4}\fontshape{#5}%
  \selectfont}%
\fi\endgroup%
\begin{picture}(3216,3873)(643,-5158)
\put(1433,-2379){\makebox(0,0)[lb]{\smash{\SetFigFont{12}{14.4}{\rmdefault}{\mddefault}{\updefault}$q_1$}}}
\put(2048,-1921){\makebox(0,0)[lb]{\smash{\SetFigFont{12}{14.4}{\rmdefault}{\mddefault}{\updefault}$q_2$}}}
\put(2461,-2116){\makebox(0,0)[lb]{\smash{\SetFigFont{12}{14.4}{\rmdefault}{\mddefault}{\updefault}$q_3$}}}
\put(2776,-2536){\makebox(0,0)[lb]{\smash{\SetFigFont{12}{14.4}{\rmdefault}{\mddefault}{\updefault}$q_4$}}}
\put(2109,-5100){\makebox(0,0)[lb]{\smash{\SetFigFont{12}{14.4}{\rmdefault}{\mddefault}{\updefault}$\emptyset$}}}
\put(2806,-4890){\makebox(0,0)[lb]{\smash{\SetFigFont{12}{14.4}{\rmdefault}{\mddefault}{\updefault}$\emptyset$}}}
\put(1336,-4920){\makebox(0,0)[lb]{\smash{\SetFigFont{12}{14.4}{\rmdefault}{\mddefault}{\updefault}$\emptyset$}}}
\put(1358,-1568){\makebox(0,0)[lb]{\smash{\SetFigFont{12}{14.4}{\rmdefault}{\mddefault}{\updefault}$\emptyset$}}}
\put(1995,-1441){\makebox(0,0)[lb]{\smash{\SetFigFont{12}{14.4}{\rmdefault}{\mddefault}{\updefault}$\emptyset$}}}
\put(2685,-1651){\makebox(0,0)[lb]{\smash{\SetFigFont{12}{14.4}{\rmdefault}{\mddefault}{\updefault}$\emptyset$}}}
\end{picture}
$$
\begin{center}
{\bf Fig. 56}
\end{center}

In fact, we have already shown that this implies automtically $w(q_1)=
w(q_2)$, $w(q_3)= w(q_4)$, $w(q_2)= -w(q_3)$, and that $q_2$, $q_3$ can be eliminated by a move 
$$
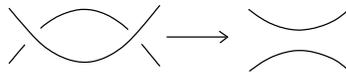
\begin{picture}(0,0)%
\includegraphics{page76H.pstex}%
\end{picture}%
\setlength{\unitlength}{4144sp}%
\begingroup\makeatletter\ifx\SetFigFont\undefined%
\gdef\SetFigFont#1#2#3#4#5{%
  \reset@font\fontsize{#1}{#2pt}%
  \fontfamily{#3}\fontseries{#4}\fontshape{#5}%
  \selectfont}%
\fi\endgroup%
\begin{picture}(2063,419)(85,-342)
\end{picture}
$$
Consequently, repeating the operation  
$$
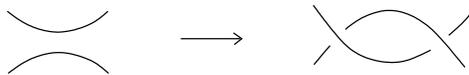
\begin{picture}(0,0)%
\includegraphics{page76B.pstex}%
\end{picture}%
\setlength{\unitlength}{4144sp}%
\begingroup\makeatletter\ifx\SetFigFont\undefined%
\gdef\SetFigFont#1#2#3#4#5{%
  \reset@font\fontsize{#1}{#2pt}%
  \fontfamily{#3}\fontseries{#4}\fontshape{#5}%
  \selectfont}%
\fi\endgroup%
\begin{picture}(2821,431)(1518,-342)
\end{picture}
$$
creates just bunches of arrows and Lemma~4.1 is proven.

\begin{lem}

{\em Let $K$ be a $\mathbb{Z}/2\mathbb{Z}$-pure global knot.

Then, $K$ is determinated by its Gauss diagram with markings in

$G \cong \mathbb{Z}/2\mathbb{Z}$ together with the homology class

$[K] \in H_1(T^2; \mathbb{Z})$.\/}
\end{lem}
\par {\bf Proof}.

As we have seen in the proof of Lemma~4.1, we can detect all possible

moves 
$$
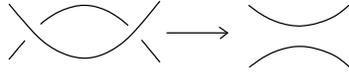
\begin{picture}(0,0)%
\includegraphics{page76H.pstex}%
\end{picture}%
\setlength{\unitlength}{4144sp}%
\begingroup\makeatletter\ifx\SetFigFont\undefined%
\gdef\SetFigFont#1#2#3#4#5{%
  \reset@font\fontsize{#1}{#2pt}%
  \fontfamily{#3}\fontseries{#4}\fontshape{#5}%
  \selectfont}%
\fi\endgroup%
\begin{picture}(2063,419)(85,-342)
\end{picture}
$$
with the Gauss diagram. Performing these moves, we obtain the minimal diagram
of $K$. By Lemma~4.1, the minimal diagram of $K$ is obtained from the torus
knot $K'$ (which is determinated by its homology class), by performing only operations 
$$
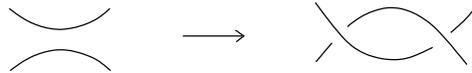
\begin{picture}(0,0)%
\includegraphics{page76B.pstex}%
\end{picture}%
\setlength{\unitlength}{4144sp}%
\begingroup\makeatletter\ifx\SetFigFont\undefined%
\gdef\SetFigFont#1#2#3#4#5{%
  \reset@font\fontsize{#1}{#2pt}%
  \fontfamily{#3}\fontseries{#4}\fontshape{#5}%
  \selectfont}%
\fi\endgroup%
\begin{picture}(2821,431)(1518,-342)
\end{picture}
$$
on the diagram of $K'$. Each such operation corresponds to a bunch of two
arrows. Thus, we only need to show that the resulting knot is completely
determined by the place of the bunch in the Gauss diagram, the directions of the arrows and their writhe. Indeed, the operation shown in Fig.~57
$$
\begin{picture}(0,0)%
\includegraphics{ima57.pstex}%
\end{picture}%
\setlength{\unitlength}{4144sp}%
\begingroup\makeatletter\ifx\SetFigFont\undefined%
\gdef\SetFigFont#1#2#3#4#5{%
  \reset@font\fontsize{#1}{#2pt}%
  \fontfamily{#3}\fontseries{#4}\fontshape{#5}%
  \selectfont}%
\fi\endgroup%
\begin{picture}(5160,2767)(548,-2225)
\put(4576,-2143){\makebox(0,0)[lb]{\smash{\SetFigFont{12}{14.4}{\rmdefault}{\mddefault}{\updefault}$I_2$}}}
\put(1808,197){\makebox(0,0)[lb]{\smash{\SetFigFont{12}{14.4}{\rmdefault}{\mddefault}{\updefault}$I_1$}}}
\put(1441,-2225){\makebox(0,0)[lb]{\smash{\SetFigFont{12}{14.4}{\rmdefault}{\mddefault}{\updefault}$I_2$}}}
\put(3981,-194){\makebox(0,0)[lb]{\smash{\SetFigFont{12}{14.4}{\rmdefault}{\mddefault}{\updefault}$q_1$}}}
\put(4011,-1148){\makebox(0,0)[lb]{\smash{\SetFigFont{12}{14.4}{\rmdefault}{\mddefault}{\updefault}$q_2$}}}
\put(4386,347){\makebox(0,0)[lb]{\smash{\SetFigFont{12}{14.4}{\rmdefault}{\mddefault}{\updefault}$I_1$}}}
\end{picture}
$$
\begin{center}
{\bf Fig. 57}
\end{center}
corresponds to the change in Fig.~58 
$$
\begin{picture}(0,0)%
\includegraphics{im27.pstex}%
\end{picture}%
\setlength{\unitlength}{4144sp}%
\begingroup\makeatletter\ifx\SetFigFont\undefined%
\gdef\SetFigFont#1#2#3#4#5{%
  \reset@font\fontsize{#1}{#2pt}%
  \fontfamily{#3}\fontseries{#4}\fontshape{#5}%
  \selectfont}%
\fi\endgroup%
\begin{picture}(5020,1126)(458,-2802)
\put(4383,-2029){\makebox(0,0)[lb]{\smash{\SetFigFont{12}{14.4}{\rmdefault}{\mddefault}{\updefault}$I_1$}}}
\put(4278,-2438){\makebox(0,0)[lb]{\smash{\SetFigFont{12}{14.4}{\rmdefault}{\mddefault}{\updefault}$I_2$}}}
\put(1038,-2469){\makebox(0,0)[lb]{\smash{\SetFigFont{12}{14.4}{\rmdefault}{\mddefault}{\updefault}$I_1$}}}
\put(1028,-1953){\makebox(0,0)[lb]{\smash{\SetFigFont{12}{14.4}{\rmdefault}{\mddefault}{\updefault}$I_2$}}}
\put(5478,-2302){\makebox(0,0)[lb]{\smash{\SetFigFont{12}{14.4}{\rmdefault}{\mddefault}{\updefault}$\beta$}}}
\end{picture}
$$
\begin{center}
{\bf Fig. 58}
\end{center}
if $w(q_1)= w(q_2)= +1$, or to the

change in Fig.~59 if $w(q_1)= w(q_2)= -1$. 
$$
\begin{picture}(0,0)%
\includegraphics{im28.pstex}%
\end{picture}%
\setlength{\unitlength}{4144sp}%
\begingroup\makeatletter\ifx\SetFigFont\undefined%
\gdef\SetFigFont#1#2#3#4#5{%
  \reset@font\fontsize{#1}{#2pt}%
  \fontfamily{#3}\fontseries{#4}\fontshape{#5}%
  \selectfont}%
\fi\endgroup%
\begin{picture}(5020,1119)(458,-2794)
\put(1038,-2469){\makebox(0,0)[lb]{\smash{\SetFigFont{12}{14.4}{\rmdefault}{\mddefault}{\updefault}$I_1$}}}
\put(1028,-1953){\makebox(0,0)[lb]{\smash{\SetFigFont{12}{14.4}{\rmdefault}{\mddefault}{\updefault}$I_2$}}}
\put(5478,-2302){\makebox(0,0)[lb]{\smash{\SetFigFont{12}{14.4}{\rmdefault}{\mddefault}{\updefault}$\beta$}}}
\put(4268,-1958){\makebox(0,0)[lb]{\smash{\SetFigFont{12}{14.4}{\rmdefault}{\mddefault}{\updefault}$I_2$}}}
\put(4358,-2466){\makebox(0,0)[lb]{\smash{\SetFigFont{12}{14.4}{\rmdefault}{\mddefault}{\updefault}$I_1$}}}
\end{picture}
$$
\begin{center}
{\bf Fig. 59}
\end{center}
Lemma~4.2 is proven.

{\sl Step 2}
\begin{lem}
{\it Let $K$ be a $\mathbb{Z}/2\mathbb{Z}$-pure global knot.}\\
{\it A) Let }
$$
\begin{picture}(0,0)%
\includegraphics{ima60.pstex}%
\end{picture}%
\setlength{\unitlength}{4144sp}%
\begingroup\makeatletter\ifx\SetFigFont\undefined%
\gdef\SetFigFont#1#2#3#4#5{%
  \reset@font\fontsize{#1}{#2pt}%
  \fontfamily{#3}\fontseries{#4}\fontshape{#5}%
  \selectfont}%
\fi\endgroup%
\begin{picture}(1350,1484)(653,-892)
\put(1647,-278){\makebox(0,0)[lb]{\smash{\SetFigFont{12}{14.4}{\rmdefault}{\mddefault}{\updefault}$q_2$}}}
\put(1141,-233){\makebox(0,0)[lb]{\smash{\SetFigFont{12}{14.4}{\rmdefault}{\mddefault}{\updefault}$q_1$}}}
\put(1253,427){\makebox(0,0)[lb]{\smash{\SetFigFont{12}{14.4}{\rmdefault}{\mddefault}{\updefault}$I$}}}
\end{picture}
$$
\begin{center}
{\bf Fig. 60}
\end{center}
{\it the diagram in Fig.~60 occur as subdiagram of the minimal diagram
  of $K$ in such a way that $I= \emptyset$ (i.e. there do not start or end any arrows in $I$). Then $w(q_1)= w(q_2)$.}\\
 {\it Let}
$$
\begin{picture}(0,0)%
\includegraphics{ima61.pstex}%
\end{picture}%
\setlength{\unitlength}{4144sp}%
\begingroup\makeatletter\ifx\SetFigFont\undefined%
\gdef\SetFigFont#1#2#3#4#5{%
  \reset@font\fontsize{#1}{#2pt}%
  \fontfamily{#3}\fontseries{#4}\fontshape{#5}%
  \selectfont}%
\fi\endgroup%
\begin{picture}(1377,1484)(629,-892)
\put(1647,-278){\makebox(0,0)[lb]{\smash{\SetFigFont{12}{14.4}{\rmdefault}{\mddefault}{\updefault}$q_2$}}}
\put(1141,-233){\makebox(0,0)[lb]{\smash{\SetFigFont{12}{14.4}{\rmdefault}{\mddefault}{\updefault}$q_1$}}}
\put(1253,427){\makebox(0,0)[lb]{\smash{\SetFigFont{12}{14.4}{\rmdefault}{\mddefault}{\updefault}$I$}}}
\put(629,-795){\makebox(0,0)[lb]{\smash{\SetFigFont{12}{14.4}{\rmdefault}{\mddefault}{\updefault}$J$}}}
\end{picture}
$$ 
\begin{center}
{\bf Fig. 61}
\end{center}
 {\it the diagram in Fig.~61 occur as subdiagram of the minimal Gauss
diagram of $K$ in such a way that $I= \emptyset$, $J= \emptyset$.
Then $w(q_3)= -w(q_1)= -w(q_2)$.\/}
\end{lem}
\par {\bf Remark}.
\par Lemma~4.1 implies that all three crossings belong to different %
bunches. Evidently, Lemma~4.3 allows to calculate all the writhes of a minimal
diagram if one knows the writhe of {\em one\/} bunch of arrows.
\par {\bf Proof of Lemma~4.3}.
We will prove only A). The proof of B) is similar, and is therefore
omitted. Using Lemma~4.1, we can eliminate all crossings of $K$, except of the four crossings shown in Fig.~62.
$$
\begin{picture}(0,0)%
\includegraphics{ima62.pstex}%
\end{picture}%
\setlength{\unitlength}{4144sp}%
\begingroup\makeatletter\ifx\SetFigFont\undefined%
\gdef\SetFigFont#1#2#3#4#5{%
  \reset@font\fontsize{#1}{#2pt}%
  \fontfamily{#3}\fontseries{#4}\fontshape{#5}%
  \selectfont}%
\fi\endgroup%
\begin{picture}(1350,1484)(653,-892)
\put(1253,427){\makebox(0,0)[lb]{\smash{\SetFigFont{12}{14.4}{\rmdefault}{\mddefault}{\updefault}$I$}}}
\put(1364,142){\makebox(0,0)[lb]{\smash{\SetFigFont{12}{14.4}{\rmdefault}{\mddefault}{\updefault}$q_2$}}}
\put(1176, 15){\makebox(0,0)[lb]{\smash{\SetFigFont{12}{14.4}{\rmdefault}{\mddefault}{\updefault}$q_1$}}}
\put(1142,-459){\makebox(0,0)[lb]{\smash{\SetFigFont{12}{14.4}{\rmdefault}{\mddefault}{\updefault}$p_1$}}}
\put(1414,-646){\makebox(0,0)[lb]{\smash{\SetFigFont{12}{14.4}{\rmdefault}{\mddefault}{\updefault}$p_2$}}}
\end{picture}
$$
\begin{center}
{\bf Fig. 62}
\end{center}
We know already that $w(q_1)= w(p_1)$, $w(q_2)= w(p_2)$. After suitable Dehn
twists, and after making $I$ small, $K$ is transformed into a knot $K'$ so
that one of the possibilities depicted in Fig.~63 is realized. In both cases,
$w(q_1)= w(q_2)$. $\Box$

$$
\begin{picture}(0,0)%
\includegraphics{im29.pstex}%
\end{picture}%
\setlength{\unitlength}{4144sp}%
\begingroup\makeatletter\ifx\SetFigFont\undefined%
\gdef\SetFigFont#1#2#3#4#5{%
  \reset@font\fontsize{#1}{#2pt}%
  \fontfamily{#3}\fontseries{#4}\fontshape{#5}%
  \selectfont}%
\fi\endgroup%
\begin{picture}(4519,3013)(451,-2613)
\put(451,235){\makebox(0,0)[lb]{\smash{\SetFigFont{12}{14.4}{\rmdefault}{\mddefault}{\updefault}Case $1$}}}
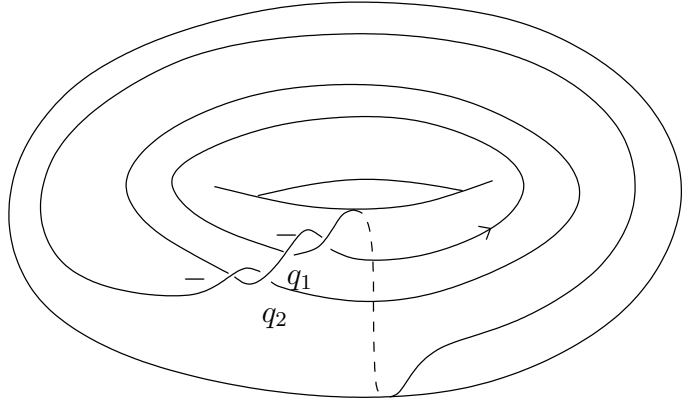
\put(2454,-1723){\makebox(0,0)[lb]{\smash{\SetFigFont{12}{14.4}{\rmdefault}{\mddefault}{\updefault}$q_2$}}}
\put(2604,-1499){\makebox(0,0)[lb]{\smash{\SetFigFont{12}{14.4}{\rmdefault}{\mddefault}{\updefault}$q_1$}}}
\put(2528,-1250){\makebox(0,0)[lb]{\smash{\SetFigFont{12}{14.4}{\rmdefault}{\mddefault}{\updefault}$-$}}}
\put(1980,-1512){\makebox(0,0)[lb]{\smash{\SetFigFont{12}{14.4}{\rmdefault}{\mddefault}{\updefault}$-$}}}
\put(2566,-2555){\makebox(0,0)[lb]{\smash{\SetFigFont{12}{14.4}{\rmdefault}{\mddefault}{\updefault}$[K']=3\alpha+ \beta$}}}
\end{picture}
$$
$$
\begin{picture}(0,0)%
\includegraphics{im29bis.pstex}%
\end{picture}%
\setlength{\unitlength}{4144sp}%
\begingroup\makeatletter\ifx\SetFigFont\undefined%
\gdef\SetFigFont#1#2#3#4#5{%
  \reset@font\fontsize{#1}{#2pt}%
  \fontfamily{#3}\fontseries{#4}\fontshape{#5}%
  \selectfont}%
\fi\endgroup%
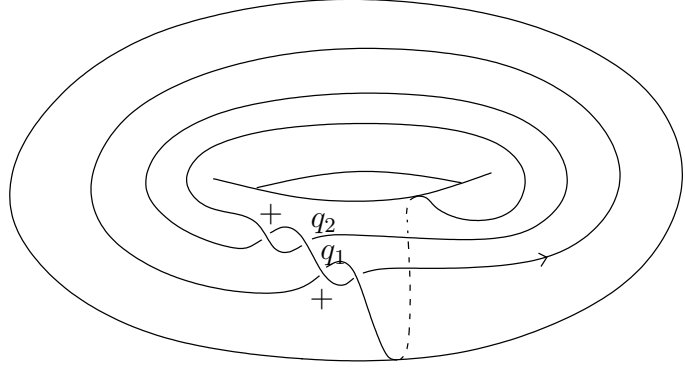
\begin{picture}(4549,2678)(496,-2349)
\put(2502,-1173){\makebox(0,0)[lb]{\smash{\SetFigFont{12}{14.4}{\rmdefault}{\mddefault}{\updefault}$+$}}}
\put(2812,-1203){\makebox(0,0)[lb]{\smash{\SetFigFont{12}{14.4}{\rmdefault}{\mddefault}{\updefault}$q_2$}}}
\put(2872,-1393){\makebox(0,0)[lb]{\smash{\SetFigFont{12}{14.4}{\rmdefault}{\mddefault}{\updefault}$q_1$}}}
\put(2807,-1683){\makebox(0,0)[lb]{\smash{\SetFigFont{12}{14.4}{\rmdefault}{\mddefault}{\updefault}$+$}}}
\put(2506,-2291){\makebox(0,0)[lb]{\smash{\SetFigFont{12}{14.4}{\rmdefault}{\mddefault}{\updefault}$[K']=3\alpha-\beta$}}}
\put(496,164){\makebox(0,0)[lb]{\smash{\SetFigFont{12}{14.4}{\rmdefault}{\mddefault}{\updefault}Case $2$}}}
\end{picture}
$$
\begin{center}
{\bf Fig. 63}
\end{center}
Let $K$ be the diagram of a $\mathbb{Z}/2\mathbb{Z}$-pure global knot with $c$
crossings and let $K'$ be the corresponding minimal diagram with $c'$
crossings. (We remind that a "diagram" is a knot together with his regular
projection into $T^2$.) Each $T$-invariant $T_K$ of degree $(d, 0)$ is $0$ for
$d> c'$. Indeed, $K$ is isotopic to $K'$ and in the Gauss diagram of $K'$,
there are not any configurations of $d$ arrows. Let $D$ be the Gauss diagram
of $K'$, and let $\bar D$ be the Gauss diagram of $K'$, without the
writhes. Thus, $\bar D$ is a $\mathbb{Z}/2\mathbb{Z}$-pure configuration of degree $c'$ (see Def.~3.2).
By Theorem~1,
$$
T_K(\bar D; \emptyset):= \sum_{\bar D}w(p_1) \cdots w(p_{c'})
$$
is an isotopy invariant of $K$. In each bunch, there is an  even number of
arrows, and, consequently, $T_K(\bar D; \emptyset)= +1$. By Lemma~4.3, there
are only two possibilities for the writhes of $D$. Therefore, the
$T$-invariant (of finite type) $T_K(\bar D; \emptyset)$ of degree $(c', 0)$
almost determines $K$: there are at most two knots with the same
invariant. Their Gauss diagrams are obtained one from the other by a
simultaneous switch of the writhes. Let $K_1$ and $K_2$ be the corresponding
knots (we know already that they are determined by their Gauss diagrams). We
have to distinguish them by $T$-invariants of smaller degree. 
\par Let $p$ be an arrow of $\bar D$ and let $\bar D_p$ be the configuration
$\bar D \setminus p$ of degree $c'- 1$. Of course, different $p$ could
determine the same configuration $\bar D_p$. Evidently, for each configuration $\bar D_p$, we have
$$
T_{K_1}(\bar D_p; \emptyset):= \sum_{\bar D_p \subset D(K_1)}\prod_{p_i \in \bar D_p}{w(p_i)}= -T_{K_2}(\bar D_p; \emptyset)
$$
Consequently, if $T_K(\bar D_p; \emptyset)\not= 0$ for some $\bar D_p$, then
$T_K(\bar D; \emptyset)$ together with $T_K(\bar D_p; \emptyset)$ determine
$K$. Assume that for all $p$, $T_K(\bar D_p; \emptyset)= 0$. Evidently, for
all couples $(p_1, p_2)$ of arrows in $\bar D$ and the corresponding configurations $\bar D_{(p_1, p_2)}:= \bar D \setminus \{ p_1, p_2 \}$, we have
$$
T_{K_1}(\bar D_{(p_1, p_2)}; \emptyset)= 
$$
$$
T_{K_2}(\bar D_{(p_1, p_2)}; \emptyset)
$$
Therefore, we go on with considering all triples $(p_1, p_2, p_3)$ of arrows in $\bar D$. For the configurations $\bar D_{(p_1, p_2, p_3)}:= \bar D \setminus \{p_1, p_2, p_3 \}$, we have 
$$
T_{K_1}(\bar D_{(p_1, p_2, p_3)}; \emptyset)= 
$$
$$
-T_{K_2}(\bar D_{(p_1, p_2, p_3)}; \emptyset)
$$
If again for all triples $(p_1, p_2, p_3)$, one has $T_K(\bar D_{(p_1, p_2,
  p_3)}; \emptyset)= 0$, then we continue with 5-tuples $(p_1, p_2, p_3, p_4, p_5)$ and so on \dots
At the end, we have either distinguished $K_1$ from $K_2$ or proven that
$T_{K_1}(c; \emptyset)= T_{K_2}(c; \emptyset)$ for any $\mathbb{Z}/
2\mathbb{Z}$-pure configuration $c$ (of course, $T_K(c; \emptyset)= 0$ for all configurations $c$ which are not subconfigurations of $\bar D$).
One easily sees that in the latter case, the Gauss diagrams of $K_1$ and $K_2$
are isotopic, and hence, by Lemma~4.2, $K_1$ and $K_2$ are isotopic too. Theorem~2 is proven. $\Box$
\newpage
\section{Non-invertibility of knots in $T^2 \times \mathbb{R}$}
Let $flip: T^2 \to T^2$ be the hyper-elliptic involution shown in
Fig.~64. 
$$
\begin{picture}(0,0)%
\includegraphics{im30.pstex}%
\end{picture}%
\setlength{\unitlength}{4144sp}%
\begingroup\makeatletter\ifx\SetFigFont\undefined%
\gdef\SetFigFont#1#2#3#4#5{%
  \reset@font\fontsize{#1}{#2pt}%
  \fontfamily{#3}\fontseries{#4}\fontshape{#5}%
  \selectfont}%
\fi\endgroup%
\begin{picture}(4435,1738)(100,-1649)
\put(4370,-645){\makebox(0,0)[lb]{\smash{\SetFigFont{12}{14.4}{\rmdefault}{\mddefault}{\updefault}$\pi$}}}
\put(2550,-1360){\makebox(0,0)[lb]{\smash{\SetFigFont{12}{14.4}{\rmdefault}{\mddefault}{\updefault}$\beta$}}}
\put(3505,-90){\makebox(0,0)[lb]{\smash{\SetFigFont{12}{14.4}{\rmdefault}{\mddefault}{\updefault}$T^2$}}}
\end{picture}
$$
\begin{center}
{\bf Fig. 64}
\end{center}
The orientation preserving involution $flip \times id$ on $T^2 \times
\mathbb{R}$ will also be called $flip$ for simplicity. $flip$ acts as $-1$ on $H_1(T^2; \mathbb{Z})$.
\par Let $K \hookrightarrow T^2 \times \mathbb{R}$ be any oriented knot, and
let $-flip(K)= flip(-K)$ be the knot obtained from $flip(K)$ by reversing its o\begin{defin}
{\rm $K$ is called {\it invertible\/} if $K$ is ambient isotopic to $-flip(K)$ in
$T^2 \times \mathbb{R}$. Otherwise, $K$ is called {\it non-invertible\/}}.
\end{defin}
\par {\bf Remarks}.
\begin{enumerate}
\item
We show in the next section that quantum invariants do not detect non-invertibility.
\item
We show in sect.~7 that our notion of invertibility for knots in $T^2
\times{R}$ coincides with the usual notion of invertibility for certain links in $S^3$.
\item
The knot $-flip(K)$ is always homotopic to $K$ in $T^2 \times \mathbb{R}$
\end{enumerate}
\par Let $v$ be our standard vector field on $T^2$ (see sect.~4). Let $K
\hookrightarrow T^2 \times \mathbb{R}$ be a (canonically oriented) global knot
with respect to $v$. Let $G$ be a quotient group of $H_1(T^2;
\mathbb{Z})/\langle [K] \rangle$. We assume that $K$ is a $G$-pure global
knot. Let $p$ be a crossing of $K$, and let $p'$ be the corresponding crossing of $-flip(K)$.
\begin{lem}
{\it $-flip(K)$ is a $G$-pure global knot with respect to $v$ too. Moreover,\/}
\begin{enumerate}
\item
$w(p)= w(p')$
\item
$[K^+_{p'}]= -[K^+_p]$ {\it in\/} $G$
\item
{\it Let $D \subset \mathbb{R}^2$ be the Gauss diagram of $K$ without writhes
  and homological markings. Then, $-flip(D) \subset \mathbb{R}^2$ is obtained
  from $D$ by a reflection with respect to any line in $\mathbb{R}^2$,
  followed by the reversion of the orientation of the circle.\/}
\end{enumerate}
\end{lem}
{\bf Proof}.
$flip(K)$ is a knot transversal to $v$, but with the wrong orientation.
Hence, $-flip(K)$ is a global knot. $flip: T^2 \times \mathbb{R} \to T^2
\times \mathbb{R}$ preserves the orientation and, consequently, $w(p)= w(p')$
for each crossing $p$. The involution $flip$ maps $K^+_p$ to
$K^+_{p'}$. Reversing the orientation of $flip(K)$, the knot $-K^-_{p'}$ for
$flip(K)$ is mapped to the knot $K^+_{p'}$ for $-flip(K)$. Thus, $[K^+_{p'}]=
[K^-_p]=[K]- [K^+_p]$ in $H_1(T^2; \mathbb{Z})$, and hence,
$[K^+_{p'}]=-[K^+_p]$ in $G$. If $[K^+_p] \not= 0$ in $G$, then $K$ is
$G$-pure. Therefore, if $K$ is $G$-pure, $-flip(K)$ is $G$-pure too. Let $D_K$
be the Gauss diagram of $K$ without writhes and homological markings (in
$G$). The circle of $D_K$ is always supposed to be embedded in the standard
way in $\mathbb{R}^2: \br$ $D_{flip(K)}$ is exactly the same Gauss diagram
(but the knots $K$ and $flip(K)$ are embedded in different ways), because
$flip$ preserves the orientation of the lines $\mathbb{R}$ and, hence,
preserves undercrosses and overcrosses.
\par Changing the orientation of $flip(K)$ changes only the orientation of the
Gauss diagram $D_{flip(K)}= D_K$. To obtain the standard embedding of the
circle in the plane, we only need to perform a reflection with respect to a line in the plane. $\Box$
\begin{defin}
{\rm Let $D$ be any $G$-pure configuration with markings in $G$. The {\it
    inverse configuration\/} $\bar D$ is obtained from $D$ by the successive operations:}
\end{defin}
\begin{enumerate}
\item
{\em a reflection with respect to a line in the plane\/}
\item
{\em reversing the orientation of the circle\/}
\item
{\em replacing each marking $a \in G$ by $-a \in G$.\/}
\end{enumerate}
{\bf Remark}
The inverse configuration is also a $G$-pure configuration.
\begin{lem}
{\it Let $K \hookrightarrow T^2 \times \mathbb{R}$ be a $G$-pure global knot
  and let $D$ be any $G$-pure configuration. If $K$ is invertible, then for the $T$-invariants (of finite type), the following holds}:
$$
T_K(D; \emptyset)= T_K(\bar D; \emptyset)
$$
\end{lem}
{\bf Proof}.
$K$ and $-flip(K)= flip(-K)$ are $G$-pure global knots (with respect to the
same $v$), and they are homotopic. Lemma~5.2 follows then immediately from
Theorem~1, Lemma~5.1, and the definition of the inverse configuration $\bar D$. $\Box$
\par {\bf Remark}
Lemma~5.2 can be generalized in a straightforward way to the case of general
$T$-invariants $T_K(D; c_1(D)= c_1, \dots, c_k(D)= c_k)$ for $G$-pure global
knots (see Def.~3.4). In particular, the {\em inverse class\/} $\bar c(\bar
D)$ is defined exactly as $c(D)$, replacing $D$ by $\bar D$ and each
configuration $\mathcal{D}_i$ by its inverse configuration $\bar
\mathcal{D}_i$ (see Def.~3.3 and 5.2). For example, if $K$ is invertible, then $T_K(\emptyset; c(\emptyset))=T_K(\emptyset; \bar c(\emptyset))$ for each $G$-pure class $c(\emptyset)$.
\begin{lem}
{\it Let $K \hookrightarrow T^2 \times \mathbb{R}$ be a $\mathbb{Z}/
  2\mathbb{Z}$-pure global knot, let $D$ be the corresponding minimal
  configuration, and let $D_s \subset D$ be a subconfiguration of highest odd
  degree such that $T_K(D_s; \emptyset) \not= 0$ (see sect.~4). then, $K$ is invertible if and only if\/}
$$
T_K(D; \emptyset)= T_K(\bar D; \emptyset)
$$
{\em and\/}
$$
T_K(D_s; \emptyset)= T_K(\bar D_s; \emptyset)
$$
\end{lem}
{\bf Proof}.
As shown in the proof of Theorem~2, $T_K(D; \emptyset)$ and $T_K(D_s;
\emptyset)$ determine the knot $K$. Lemma~5.3 follows then immediately from Lemma~5.2. $\Box$
{\it Example~5.1}
$$
\begin{picture}(0,0)%
\includegraphics{im31.pstex}%
\end{picture}%
\setlength{\unitlength}{4144sp}%
\begingroup\makeatletter\ifx\SetFigFont\undefined%
\gdef\SetFigFont#1#2#3#4#5{%
  \reset@font\fontsize{#1}{#2pt}%
  \fontfamily{#3}\fontseries{#4}\fontshape{#5}%
  \selectfont}%
\fi\endgroup%
\begin{picture}(4036,2470)(310,-2157)
\put(1001,-1924){\makebox(0,0)[lb]{\smash{\SetFigFont{12}{14.4}{\rmdefault}{\mddefault}{\updefault}$+$}}}
\put(1201,-1929){\makebox(0,0)[lb]{\smash{\SetFigFont{12}{14.4}{\rmdefault}{\mddefault}{\updefault}$+$}}}
\put(1361,-1654){\makebox(0,0)[lb]{\smash{\SetFigFont{12}{14.4}{\rmdefault}{\mddefault}{\updefault}$-$}}}
\put(1551,-1634){\makebox(0,0)[lb]{\smash{\SetFigFont{12}{14.4}{\rmdefault}{\mddefault}{\updefault}$-$}}}
\put(1786,-1634){\makebox(0,0)[lb]{\smash{\SetFigFont{12}{14.4}{\rmdefault}{\mddefault}{\updefault}$-$}}}
\put(2016,-1649){\makebox(0,0)[lb]{\smash{\SetFigFont{12}{14.4}{\rmdefault}{\mddefault}{\updefault}$-$}}}
\put(2181,-1174){\makebox(0,0)[lb]{\smash{\SetFigFont{12}{14.4}{\rmdefault}{\mddefault}{\updefault}$-$}}}
\put(2406,-1154){\makebox(0,0)[lb]{\smash{\SetFigFont{12}{14.4}{\rmdefault}{\mddefault}{\updefault}$-$}}}
\put(4081,-909){\makebox(0,0)[lb]{\smash{\SetFigFont{12}{14.4}{\rmdefault}{\mddefault}{\updefault}$K$}}}
\end{picture}
$$
\begin{center}
{\bf Fig. 65}
\end{center}
The knot shown in Fig.~65 represents $4\alpha+ \beta$ in $H_1(T^2; \mathbb{Z})$ and is a $\mathbb{Z}/2\mathbb{Z}$-pure global knot. Its Gauss diagram is shown in Fig.~66.
$$
\begin{picture}(0,0)%
\includegraphics{im32.pstex}%
\end{picture}%
\setlength{\unitlength}{4144sp}%
\begingroup\makeatletter\ifx\SetFigFont\undefined%
\gdef\SetFigFont#1#2#3#4#5{%
  \reset@font\fontsize{#1}{#2pt}%
  \fontfamily{#3}\fontseries{#4}\fontshape{#5}%
  \selectfont}%
\fi\endgroup%
\begin{picture}(2390,2415)(1061,-2161)
\put(2244,-23){\makebox(0,0)[lb]{\smash{\SetFigFont{12}{14.4}{\rmdefault}{\mddefault}{\updefault}$-$}}}
\put(2744,-438){\makebox(0,0)[lb]{\smash{\SetFigFont{12}{14.4}{\rmdefault}{\mddefault}{\updefault}$-$}}}
\put(1681,-347){\makebox(0,0)[lb]{\smash{\SetFigFont{12}{14.4}{\rmdefault}{\mddefault}{\updefault}$+$}}}
\put(1264,-1130){\makebox(0,0)[lb]{\smash{\SetFigFont{12}{14.4}{\rmdefault}{\mddefault}{\updefault}$+$}}}
\put(1592,-1863){\makebox(0,0)[lb]{\smash{\SetFigFont{12}{14.4}{\rmdefault}{\mddefault}{\updefault}$-$}}}
\put(2060,-1842){\makebox(0,0)[lb]{\smash{\SetFigFont{12}{14.4}{\rmdefault}{\mddefault}{\updefault}$-$}}}
\put(2531,-1722){\makebox(0,0)[lb]{\smash{\SetFigFont{12}{14.4}{\rmdefault}{\mddefault}{\updefault}$-$}}}
\put(2918,-1410){\makebox(0,0)[lb]{\smash{\SetFigFont{12}{14.4}{\rmdefault}{\mddefault}{\updefault}$-$}}}
\end{picture}
$$
\begin{center}
{\bf Fig. 66}
\end{center}
Let $D$ be the configuration of degree 6 shown in Fig.~67.
$$
\begin{picture}(0,0)%
\includegraphics{im33.pstex}%
\end{picture}%
\setlength{\unitlength}{4144sp}%
\begingroup\makeatletter\ifx\SetFigFont\undefined%
\gdef\SetFigFont#1#2#3#4#5{%
  \reset@font\fontsize{#1}{#2pt}%
  \fontfamily{#3}\fontseries{#4}\fontshape{#5}%
  \selectfont}%
\fi\endgroup%
\begin{picture}(2390,2415)(1061,-2161)
\end{picture}
$$
\begin{center}
{\bf Fig. 67}
\end{center}
Evidently, each configuration which does not contain a subconfiguration as depicted in Fig.~68 
$$
\begin{picture}(0,0)%
\includegraphics{ima68.pstex}%
\end{picture}%
\setlength{\unitlength}{4144sp}%
\begingroup\makeatletter\ifx\SetFigFont\undefined%
\gdef\SetFigFont#1#2#3#4#5{%
  \reset@font\fontsize{#1}{#2pt}%
  \fontfamily{#3}\fontseries{#4}\fontshape{#5}%
  \selectfont}%
\fi\endgroup%
\begin{picture}(2234,2558)(3474,-2098)
\put(4547,328){\makebox(0,0)[lb]{\smash{\SetFigFont{12}{14.4}{\rmdefault}{\mddefault}{\updefault}$\emptyset$}}}
\put(4420,-2049){\makebox(0,0)[lb]{\smash{\SetFigFont{12}{14.4}{\rmdefault}{\mddefault}{\updefault}$\emptyset$}}}
\end{picture}
$$
\begin{center}
{\bf Fig. 68}
\end{center}
is a $\mathbb{Z}/2\mathbb{Z}$-pure configuration (see Def.~3.2). Consequently,
$D$ is a $\mathbb{Z}/2\mathbb{Z}$-pure configuration. Using Fig.~66 it takes
some seconds to calculate $T_K(D; \emptyset)= -1$. The inverse configuration
$\bar D$ is shown in Fig.~69.
$$
\begin{picture}(0,0)%
\includegraphics{im34.pstex}%
\end{picture}%
\setlength{\unitlength}{4144sp}%
\begingroup\makeatletter\ifx\SetFigFont\undefined%
\gdef\SetFigFont#1#2#3#4#5{%
  \reset@font\fontsize{#1}{#2pt}%
  \fontfamily{#3}\fontseries{#4}\fontshape{#5}%
  \selectfont}%
\fi\endgroup%
\begin{picture}(2390,2412)(1073,-2140)
\end{picture}
$$
\begin{center}
{\bf Fig. 69}
\end{center}
We see immediately that $\bar D$ does not appear at all in the Gauss diagram
of $K$ (the cyclic ordering of the bunches has changed). Therefore, $T_K(\bar
D; \emptyset)= 0$, and the knot $K$ is not invertible according to Lemma~5.3.
{\it Example~5.2}
This is a more complicated example. 
$$
\begin{picture}(0,0)%
\includegraphics{im35.pstex}%
\end{picture}%
\setlength{\unitlength}{4144sp}%
\begingroup\makeatletter\ifx\SetFigFont\undefined%
\gdef\SetFigFont#1#2#3#4#5{%
  \reset@font\fontsize{#1}{#2pt}%
  \fontfamily{#3}\fontseries{#4}\fontshape{#5}%
  \selectfont}%
\fi\endgroup%
\begin{picture}(4229,2616)(929,-2322)
\put(3001,-1223){\makebox(0,0)[lb]{\smash{\SetFigFont{12}{14.4}{\rmdefault}{\mddefault}{\updefault}$+$}}}
\put(3501,-1278){\makebox(0,0)[lb]{\smash{\SetFigFont{12}{14.4}{\rmdefault}{\mddefault}{\updefault}$+$}}}
\put(2686,-1623){\makebox(0,0)[lb]{\smash{\SetFigFont{12}{14.4}{\rmdefault}{\mddefault}{\updefault}$+$}}}
\put(2491,-1613){\makebox(0,0)[lb]{\smash{\SetFigFont{12}{14.4}{\rmdefault}{\mddefault}{\updefault}$+$}}}
\put(2211,-1798){\makebox(0,0)[lb]{\smash{\SetFigFont{12}{14.4}{\rmdefault}{\mddefault}{\updefault}$+$}}}
\put(1991,-1798){\makebox(0,0)[lb]{\smash{\SetFigFont{12}{14.4}{\rmdefault}{\mddefault}{\updefault}$+$}}}
\put(3161,-1628){\makebox(0,0)[lb]{\smash{\SetFigFont{12}{14.4}{\rmdefault}{\mddefault}{\updefault}$-$}}}
\put(3431,-1613){\makebox(0,0)[lb]{\smash{\SetFigFont{12}{14.4}{\rmdefault}{\mddefault}{\updefault}$-$}}}
\put(1776,-2003){\makebox(0,0)[lb]{\smash{\SetFigFont{12}{14.4}{\rmdefault}{\mddefault}{\updefault}$-$}}}
\put(1516,-1968){\makebox(0,0)[lb]{\smash{\SetFigFont{12}{14.4}{\rmdefault}{\mddefault}{\updefault}$-$}}}
\put(5158,-1483){\makebox(0,0)[lb]{\smash{\SetFigFont{12}{14.4}{\rmdefault}{\mddefault}{\updefault}$K$}}}
\end{picture}
$$
\begin{center}
{\bf Fig. 70}
\end{center}
The knot $K$ drawn in Fig.~70 is a global knot which represents $5\alpha+ \beta$ in $H_1(T^2; \mathbb{Z})$.
\par Let $G:= (H_1(T^2)/\langle [k] \rangle)/3\mathbb{Z} \cong \mathbb{Z}/
3\mathbb{Z}= \{ 0, a, -a \}$, where the class $a$ is represented by $\alpha$. The Gauss diagram of $K$ is shown in Fig.~71.
$$
\begin{picture}(0,0)%
\includegraphics{im36.pstex}%
\end{picture}%
\setlength{\unitlength}{4144sp}%
\begingroup\makeatletter\ifx\SetFigFont\undefined%
\gdef\SetFigFont#1#2#3#4#5{%
  \reset@font\fontsize{#1}{#2pt}%
  \fontfamily{#3}\fontseries{#4}\fontshape{#5}%
  \selectfont}%
\fi\endgroup%
\begin{picture}(2685,2816)(811,-2377)
\put(1671,-331){\makebox(0,0)[lb]{\smash{\SetFigFont{12}{14.4}{\rmdefault}{\mddefault}{\updefault}$-$}}}
\put(2156,-11){\makebox(0,0)[lb]{\smash{\SetFigFont{12}{14.4}{\rmdefault}{\mddefault}{\updefault}$-$}}}
\put(2086,274){\makebox(0,0)[lb]{\smash{\SetFigFont{12}{14.4}{\rmdefault}{\mddefault}{\updefault}$a$}}}
\put(816,-551){\makebox(0,0)[lb]{\smash{\SetFigFont{12}{14.4}{\rmdefault}{\mddefault}{\updefault}$-a$}}}
\put(1221,-926){\makebox(0,0)[lb]{\smash{\SetFigFont{12}{14.4}{\rmdefault}{\mddefault}{\updefault}$+$}}}
\put(1391,-1576){\makebox(0,0)[lb]{\smash{\SetFigFont{12}{14.4}{\rmdefault}{\mddefault}{\updefault}$+$}}}
\put(3066,-1191){\makebox(0,0)[lb]{\smash{\SetFigFont{12}{14.4}{\rmdefault}{\mddefault}{\updefault}$+$}}}
\put(2941,-1526){\makebox(0,0)[lb]{\smash{\SetFigFont{12}{14.4}{\rmdefault}{\mddefault}{\updefault}$+$}}}
\put(2056,-1901){\makebox(0,0)[lb]{\smash{\SetFigFont{12}{14.4}{\rmdefault}{\mddefault}{\updefault}$+$}}}
\put(3016,-791){\makebox(0,0)[lb]{\smash{\SetFigFont{12}{14.4}{\rmdefault}{\mddefault}{\updefault}$-$}}}
\put(2521,-756){\makebox(0,0)[lb]{\smash{\SetFigFont{12}{14.4}{\rmdefault}{\mddefault}{\updefault}$-$}}}
\put(2826,129){\makebox(0,0)[lb]{\smash{\SetFigFont{12}{14.4}{\rmdefault}{\mddefault}{\updefault}$a$}}}
\put(3326,-321){\makebox(0,0)[lb]{\smash{\SetFigFont{12}{14.4}{\rmdefault}{\mddefault}{\updefault}$-a$}}}
\put(3496,-751){\makebox(0,0)[lb]{\smash{\SetFigFont{12}{14.4}{\rmdefault}{\mddefault}{\updefault}$a$}}}
\put(3216,-1876){\makebox(0,0)[lb]{\smash{\SetFigFont{12}{14.4}{\rmdefault}{\mddefault}{\updefault}$-a$}}}
\put(811,-1001){\makebox(0,0)[lb]{\smash{\SetFigFont{12}{14.4}{\rmdefault}{\mddefault}{\updefault}$a$}}}
\put(2146,-2377){\makebox(0,0)[lb]{\smash{\SetFigFont{12}{14.4}{\rmdefault}{\mddefault}{\updefault}$a$}}}
\put(2616,-2307){\makebox(0,0)[lb]{\smash{\SetFigFont{12}{14.4}{\rmdefault}{\mddefault}{\updefault}$-a$}}}
\put(811,-1596){\makebox(0,0)[lb]{\smash{\SetFigFont{12}{14.4}{\rmdefault}{\mddefault}{\updefault}$-a$}}}
\put(2586,-161){\makebox(0,0)[lb]{\smash{\SetFigFont{12}{14.4}{\rmdefault}{\mddefault}{\updefault}$+$}}}
\end{picture}
$$
\begin{center}
{\bf Fig. 71}
\end{center}
Hence, $K$ is $\mathbb{Z}/3\mathbb{Z}$-pure. Let $c(\emptyset)$ be the class of degree 5 shown in Fig.~72 (the weight functions are always the products of the writhes of the 5 crossings).
$$
\begin{picture}(0,0)%
\includegraphics{ima72.pstex}%
\end{picture}%
\setlength{\unitlength}{4144sp}%
\begingroup\makeatletter\ifx\SetFigFont\undefined%
\gdef\SetFigFont#1#2#3#4#5{%
  \reset@font\fontsize{#1}{#2pt}%
  \fontfamily{#3}\fontseries{#4}\fontshape{#5}%
  \selectfont}%
\fi\endgroup%
\begin{picture}(5654,4367)(76,-3661)
\put(1426,-2521){\makebox(0,0)[lb]{\smash{\SetFigFont{12}{14.4}{\rmdefault}{\mddefault}{\updefault}$a$}}}
\put(1264,-3628){\makebox(0,0)[lb]{\smash{\SetFigFont{12}{14.4}{\rmdefault}{\mddefault}{\updefault}$-a$}}}
\put(811,-3394){\makebox(0,0)[rb]{\smash{\SetFigFont{12}{14.4}{\rmdefault}{\mddefault}{\updefault}$a$}}}
\put(1600,-2662){\makebox(0,0)[lb]{\smash{\SetFigFont{12}{14.4}{\rmdefault}{\mddefault}{\updefault}$a$}}}
\put(931,-3598){\makebox(0,0)[lb]{\smash{\SetFigFont{12}{14.4}{\rmdefault}{\mddefault}{\updefault}$a$}}}
\put(1168,571){\makebox(0,0)[lb]{\smash{\SetFigFont{12}{14.4}{\rmdefault}{\mddefault}{\updefault}$a$}}}
\put(1435,502){\makebox(0,0)[lb]{\smash{\SetFigFont{12}{14.4}{\rmdefault}{\mddefault}{\updefault}$a$}}}
\put(973,-479){\makebox(0,0)[rb]{\smash{\SetFigFont{12}{14.4}{\rmdefault}{\mddefault}{\updefault}$a$}}}
\put(1339,-533){\makebox(0,0)[lb]{\smash{\SetFigFont{12}{14.4}{\rmdefault}{\mddefault}{\updefault}$-a$}}}
\put(778,-245){\makebox(0,0)[rb]{\smash{\SetFigFont{12}{14.4}{\rmdefault}{\mddefault}{\updefault}$-a$}}}
\put(3379,493){\makebox(0,0)[lb]{\smash{\SetFigFont{12}{14.4}{\rmdefault}{\mddefault}{\updefault}$a$}}}
\put(3538,361){\makebox(0,0)[lb]{\smash{\SetFigFont{12}{14.4}{\rmdefault}{\mddefault}{\updefault}$a$}}}
\put(3190,-569){\makebox(0,0)[lb]{\smash{\SetFigFont{12}{14.4}{\rmdefault}{\mddefault}{\updefault}$-a$}}}
\put(2887,-548){\makebox(0,0)[lb]{\smash{\SetFigFont{12}{14.4}{\rmdefault}{\mddefault}{\updefault}$a$}}}
\put(2671,-236){\makebox(0,0)[rb]{\smash{\SetFigFont{12}{14.4}{\rmdefault}{\mddefault}{\updefault}$-a$}}}
\put(5574,388){\makebox(0,0)[lb]{\smash{\SetFigFont{12}{14.4}{\rmdefault}{\mddefault}{\updefault}$a$}}}
\put(5334,556){\makebox(0,0)[lb]{\smash{\SetFigFont{12}{14.4}{\rmdefault}{\mddefault}{\updefault}$a$}}}
\put(5163,-542){\makebox(0,0)[lb]{\smash{\SetFigFont{12}{14.4}{\rmdefault}{\mddefault}{\updefault}$a$}}}
\put(5109,-533){\makebox(0,0)[rb]{\smash{\SetFigFont{12}{14.4}{\rmdefault}{\mddefault}{\updefault}$-a$}}}
\put(4659,-95){\makebox(0,0)[rb]{\smash{\SetFigFont{12}{14.4}{\rmdefault}{\mddefault}{\updefault}$-a$}}}
\put(1396,-1012){\makebox(0,0)[lb]{\smash{\SetFigFont{12}{14.4}{\rmdefault}{\mddefault}{\updefault}$a$}}}
\put(1570,-1135){\makebox(0,0)[lb]{\smash{\SetFigFont{12}{14.4}{\rmdefault}{\mddefault}{\updefault}$a$}}}
\put(1222,-2071){\makebox(0,0)[lb]{\smash{\SetFigFont{12}{14.4}{\rmdefault}{\mddefault}{\updefault}$-a$}}}
\put(811,-1882){\makebox(0,0)[rb]{\smash{\SetFigFont{12}{14.4}{\rmdefault}{\mddefault}{\updefault}$-a$}}}
\put(913,-2029){\makebox(0,0)[b]{\smash{\SetFigFont{12}{14.4}{\rmdefault}{\mddefault}{\updefault}$a$}}}
\put(5538,-880){\makebox(0,0)[lb]{\smash{\SetFigFont{12}{14.4}{\rmdefault}{\mddefault}{\updefault}$a$}}}
\put(5253,-775){\makebox(0,0)[lb]{\smash{\SetFigFont{12}{14.4}{\rmdefault}{\mddefault}{\updefault}$a$}}}
\put(4764,-1582){\makebox(0,0)[rb]{\smash{\SetFigFont{12}{14.4}{\rmdefault}{\mddefault}{\updefault}$-a$}}}
\put(5127,-1891){\makebox(0,0)[rb]{\smash{\SetFigFont{12}{14.4}{\rmdefault}{\mddefault}{\updefault}$-a$}}}
\put(5190,-1894){\makebox(0,0)[lb]{\smash{\SetFigFont{12}{14.4}{\rmdefault}{\mddefault}{\updefault}$a$}}}
\put(3442,-2554){\makebox(0,0)[lb]{\smash{\SetFigFont{12}{14.4}{\rmdefault}{\mddefault}{\updefault}$a$}}}
\put(3280,-3661){\makebox(0,0)[lb]{\smash{\SetFigFont{12}{14.4}{\rmdefault}{\mddefault}{\updefault}$-a$}}}
\put(3616,-2695){\makebox(0,0)[lb]{\smash{\SetFigFont{12}{14.4}{\rmdefault}{\mddefault}{\updefault}$a$}}}
\put(2947,-3631){\makebox(0,0)[lb]{\smash{\SetFigFont{12}{14.4}{\rmdefault}{\mddefault}{\updefault}$-a$}}}
\put(2879,-3444){\makebox(0,0)[rb]{\smash{\SetFigFont{12}{14.4}{\rmdefault}{\mddefault}{\updefault}$a$}}}
\put( 76, 22){\makebox(0,0)[lb]{\smash{\SetFigFont{12}{14.4}{\rmdefault}{\mddefault}{\updefault}$c(\emptyset):=$}}}
\put(2108, 22){\makebox(0,0)[lb]{\smash{\SetFigFont{12}{14.4}{\rmdefault}{\mddefault}{\updefault}$+$}}}
\put(1966,-1523){\makebox(0,0)[lb]{\smash{\SetFigFont{12}{14.4}{\rmdefault}{\mddefault}{\updefault}$+$}}}
\put(1996,-3136){\makebox(0,0)[lb]{\smash{\SetFigFont{12}{14.4}{\rmdefault}{\mddefault}{\updefault}$+$}}}
\put(4073,-1373){\makebox(0,0)[lb]{\smash{\SetFigFont{12}{14.4}{\rmdefault}{\mddefault}{\updefault}$+$}}}
\put(4005,104){\makebox(0,0)[lb]{\smash{\SetFigFont{12}{14.4}{\rmdefault}{\mddefault}{\updefault}$+$}}}
\put(240,-1478){\makebox(0,0)[lb]{\smash{\SetFigFont{12}{14.4}{\rmdefault}{\mddefault}{\updefault}$+$}}}
\put(188,-3316){\makebox(0,0)[lb]{\smash{\SetFigFont{12}{14.4}{\rmdefault}{\mddefault}{\updefault}$+$}}}
\put(4149,-3132){\makebox(0,0)[lb]{\smash{\SetFigFont{12}{14.4}{\rmdefault}{\mddefault}{\updefault}$+$ $\ldots$}}}
\put(3326,-889){\makebox(0,0)[lb]{\smash{\SetFigFont{12}{14.4}{\rmdefault}{\mddefault}{\updefault}$a$}}}
\put(3164,-1996){\makebox(0,0)[lb]{\smash{\SetFigFont{12}{14.4}{\rmdefault}{\mddefault}{\updefault}$-a$}}}
\put(2711,-1762){\makebox(0,0)[rb]{\smash{\SetFigFont{12}{14.4}{\rmdefault}{\mddefault}{\updefault}$a$}}}
\put(3500,-1030){\makebox(0,0)[lb]{\smash{\SetFigFont{12}{14.4}{\rmdefault}{\mddefault}{\updefault}$a$}}}
\put(2831,-1966){\makebox(0,0)[lb]{\smash{\SetFigFont{12}{14.4}{\rmdefault}{\mddefault}{\updefault}$-a$}}}
\end{picture}
$$
\begin{center}
{\bf Fig. 72}
\end{center}
The only possible strata of triple points in the discriminant for a
$\mathbb{Z}/3\mathbb{Z}$-pure isotopy are $a^{\pm}_{\dr(\br)}(a|-a,|-a)$ and $a^{\pm}_{\br(\br)}(-a|a,|a)$ (see \cite{F}, sect.~1). Therefore, the changings depicted in Fig.~73 are the only possible ones for a couple of crossings in a $\mathbb{Z}/3\mathbb{Z}$-pure isotopy.
$$
\begin{picture}(0,0)%
\includegraphics{ima73.pstex}%
\end{picture}%
\setlength{\unitlength}{4144sp}%
\begingroup\makeatletter\ifx\SetFigFont\undefined%
\gdef\SetFigFont#1#2#3#4#5{%
  \reset@font\fontsize{#1}{#2pt}%
  \fontfamily{#3}\fontseries{#4}\fontshape{#5}%
  \selectfont}%
\fi\endgroup%
\begin{picture}(5461,3082)(353,-2378)
\put(1606,-1771){\makebox(0,0)[lb]{\smash{\SetFigFont{12}{14.4}{\rmdefault}{\mddefault}{\updefault}or}}}
\put(3856,-16){\makebox(0,0)[lb]{\smash{\SetFigFont{12}{14.4}{\rmdefault}{\mddefault}{\updefault}or}}}
\put(638,569){\makebox(0,0)[lb]{\smash{\SetFigFont{12}{14.4}{\rmdefault}{\mddefault}{\updefault}$x$}}}
\put(1126,547){\makebox(0,0)[lb]{\smash{\SetFigFont{12}{14.4}{\rmdefault}{\mddefault}{\updefault}$y$}}}
\put(2760,554){\makebox(0,0)[lb]{\smash{\SetFigFont{12}{14.4}{\rmdefault}{\mddefault}{\updefault}$x$}}}
\put(5235,562){\makebox(0,0)[lb]{\smash{\SetFigFont{12}{14.4}{\rmdefault}{\mddefault}{\updefault}$x$}}}
\put(540,-1231){\makebox(0,0)[lb]{\smash{\SetFigFont{12}{14.4}{\rmdefault}{\mddefault}{\updefault}$x$}}}
\put(1147,-2378){\makebox(0,0)[lb]{\smash{\SetFigFont{12}{14.4}{\rmdefault}{\mddefault}{\updefault}$x$}}}
\put(3270,-1201){\makebox(0,0)[lb]{\smash{\SetFigFont{12}{14.4}{\rmdefault}{\mddefault}{\updefault}$x$}}}
\put(2512,-2356){\makebox(0,0)[lb]{\smash{\SetFigFont{12}{14.4}{\rmdefault}{\mddefault}{\updefault}$x$}}}
\put(4654,-1162){\makebox(0,0)[lb]{\smash{\SetFigFont{12}{14.4}{\rmdefault}{\mddefault}{\updefault}$x$}}}
\put(5577,-1136){\makebox(0,0)[lb]{\smash{\SetFigFont{12}{14.4}{\rmdefault}{\mddefault}{\updefault}$x$}}}
\put(3473,554){\makebox(0,0)[lb]{\smash{\SetFigFont{12}{14.4}{\rmdefault}{\mddefault}{\updefault}$y$}}}
\put(4560,584){\makebox(0,0)[lb]{\smash{\SetFigFont{12}{14.4}{\rmdefault}{\mddefault}{\updefault}$y$}}}
\put(631,-983){\makebox(0,0)[lb]{\smash{\SetFigFont{12}{14.4}{\rmdefault}{\mddefault}{\updefault}if $x \neq y$}}}
\end{picture}
$$
\begin{center}
{\bf Fig. 73}
\end{center}
Here, $x, y \in \{ a, -a \}$. $c(\emptyset)$ is obtained from the
configuration shown in the left-hand part of Fig.~74 by applying {\em all\/}
possible changings to it. We have shown some of these changings
hereabove. Notice that no chord can ever get crossed with the isolated chord,
because the part of the configuration shown in the right-hand part of Fig.~74 cannot change at all.
$$
\begin{picture}(0,0)%
\includegraphics{ima74.pstex}%
\end{picture}%
\setlength{\unitlength}{4144sp}%
\begingroup\makeatletter\ifx\SetFigFont\undefined%
\gdef\SetFigFont#1#2#3#4#5{%
  \reset@font\fontsize{#1}{#2pt}%
  \fontfamily{#3}\fontseries{#4}\fontshape{#5}%
  \selectfont}%
\fi\endgroup%
\begin{picture}(3428,1427)(454,-721)
\put(3316,526){\makebox(0,0)[lb]{\smash{\SetFigFont{12}{14.4}{\rmdefault}{\mddefault}{\updefault}$a$}}}
\put(3583,457){\makebox(0,0)[lb]{\smash{\SetFigFont{12}{14.4}{\rmdefault}{\mddefault}{\updefault}$a$}}}
\put(1168,571){\makebox(0,0)[lb]{\smash{\SetFigFont{12}{14.4}{\rmdefault}{\mddefault}{\updefault}$a$}}}
\put(1435,502){\makebox(0,0)[lb]{\smash{\SetFigFont{12}{14.4}{\rmdefault}{\mddefault}{\updefault}$a$}}}
\put(973,-479){\makebox(0,0)[rb]{\smash{\SetFigFont{12}{14.4}{\rmdefault}{\mddefault}{\updefault}$a$}}}
\put(1339,-533){\makebox(0,0)[lb]{\smash{\SetFigFont{12}{14.4}{\rmdefault}{\mddefault}{\updefault}$-a$}}}
\put(778,-245){\makebox(0,0)[rb]{\smash{\SetFigFont{12}{14.4}{\rmdefault}{\mddefault}{\updefault}$-a$}}}
\put(3426,-585){\makebox(0,0)[lb]{\smash{\SetFigFont{12}{14.4}{\rmdefault}{\mddefault}{\updefault}$a$}}}
\end{picture}
$$
\begin{center}
{\bf Fig. 74}
\end{center}
Consequently, $c(\emptyset)$ is a $\mathbb{Z}/3\mathbb{Z}$-pure class of
degree 5 and $T_K(\emptyset; c(\emptyset))$ is an isotopy invariant of
$K$. Notice that $p$ is the only arrow in the Gauss diagram of $K$ with
marking $a$ and such that there are arrows in $K^+_p$. Using this fact, we
easily calculate $T_K(\emptyset; c(\emptyset))= -1$. For the convenience of
the reader, we give the Gauss diagram of $flip(-K)$ in Fig.~75.
$$
\begin{picture}(0,0)%
\includegraphics{im37.pstex}%
\end{picture}%
\setlength{\unitlength}{4144sp}%
\begingroup\makeatletter\ifx\SetFigFont\undefined%
\gdef\SetFigFont#1#2#3#4#5{%
  \reset@font\fontsize{#1}{#2pt}%
  \fontfamily{#3}\fontseries{#4}\fontshape{#5}%
  \selectfont}%
\fi\endgroup%
\begin{picture}(2790,2706)(381,-2660)
\put(741,-2179){\makebox(0,0)[lb]{\smash{\SetFigFont{12}{14.4}{\rmdefault}{\mddefault}{\updefault}$a$}}}
\put(1011,-1904){\makebox(0,0)[lb]{\smash{\SetFigFont{12}{14.4}{\rmdefault}{\mddefault}{\updefault}$+$}}}
\put(916,-1404){\makebox(0,0)[lb]{\smash{\SetFigFont{12}{14.4}{\rmdefault}{\mddefault}{\updefault}$+$}}}
\put(381,-1039){\makebox(0,0)[lb]{\smash{\SetFigFont{12}{14.4}{\rmdefault}{\mddefault}{\updefault}$-a$}}}
\put(686,-679){\makebox(0,0)[lb]{\smash{\SetFigFont{12}{14.4}{\rmdefault}{\mddefault}{\updefault}$a$}}}
\put(1406,-964){\makebox(0,0)[lb]{\smash{\SetFigFont{12}{14.4}{\rmdefault}{\mddefault}{\updefault}$-$}}}
\put(1011,-1054){\makebox(0,0)[lb]{\smash{\SetFigFont{12}{14.4}{\rmdefault}{\mddefault}{\updefault}$-$}}}
\put(1986,-2229){\makebox(0,0)[lb]{\smash{\SetFigFont{12}{14.4}{\rmdefault}{\mddefault}{\updefault}$+$}}}
\put(1461,-524){\makebox(0,0)[lb]{\smash{\SetFigFont{12}{14.4}{\rmdefault}{\mddefault}{\updefault}$+$}}}
\put(1966,-2655){\makebox(0,0)[lb]{\smash{\SetFigFont{12}{14.4}{\rmdefault}{\mddefault}{\updefault}$-a$}}}
\put(1511,-2660){\makebox(0,0)[lb]{\smash{\SetFigFont{12}{14.4}{\rmdefault}{\mddefault}{\updefault}$a$}}}
\put(2756,-1379){\makebox(0,0)[lb]{\smash{\SetFigFont{12}{14.4}{\rmdefault}{\mddefault}{\updefault}$+$}}}
\put(2626,-1839){\makebox(0,0)[lb]{\smash{\SetFigFont{12}{14.4}{\rmdefault}{\mddefault}{\updefault}$+$}}}
\put(2981,-2094){\makebox(0,0)[lb]{\smash{\SetFigFont{12}{14.4}{\rmdefault}{\mddefault}{\updefault}$a$}}}
\put(3171,-1409){\makebox(0,0)[lb]{\smash{\SetFigFont{12}{14.4}{\rmdefault}{\mddefault}{\updefault}$-a$}}}
\put(2161,-119){\makebox(0,0)[lb]{\smash{\SetFigFont{12}{14.4}{\rmdefault}{\mddefault}{\updefault}$-a$}}}
\put(3071,-899){\makebox(0,0)[lb]{\smash{\SetFigFont{12}{14.4}{\rmdefault}{\mddefault}{\updefault}$a$}}}
\put(2471,-869){\makebox(0,0)[lb]{\smash{\SetFigFont{12}{14.4}{\rmdefault}{\mddefault}{\updefault}$-$}}}
\put(2066,-569){\makebox(0,0)[lb]{\smash{\SetFigFont{12}{14.4}{\rmdefault}{\mddefault}{\updefault}$-$}}}
\put(1396,-219){\makebox(0,0)[rb]{\smash{\SetFigFont{12}{14.4}{\rmdefault}{\mddefault}{\updefault}$-a$}}}
\end{picture}
$$
\begin{center}
{\bf Fig. 75}
\end{center}
Hence, $T_{flip(-K)}(\emptyset; c(\emptyset))= T_K(\emptyset; \bar c(
\emptyset))= 0$. (Again, the cyclic order of the two couples of crossed arrows
with respect to the isolated arrow $p'$ has changed.) We have proven that $K$
i\newpage
\section{A remark on quantum invariants for knots in $T^2 \times \mathbb{R}$}
\par Let $L \hookrightarrow T^2 \times \mathbb{R}$ be any oriented link. There
are generalized HOMFLY-PT and Kauffman polynomials for $L$ (see e.g.\cite{H-P}).
\begin{lem}
{\it The generalized HOMFLY-PT polynomials of $L$ and of $flip(-L)$ coincide. The generalized Kauffman polynomials of $L$ and of $flip(-L)$
coincide.}
\end{lem}
\par {\bf Proof}.$L$ can be reduced to a linear combination of "initial knots"
by using skein relations. These combinations are the same for $flip(-L)$
besides the fact that each "initial knot" $K$ has to be replaced by
$flip(-K)$. Consequently, if we find a set of "initial knots" such that for
each of them, $K$ is isotopic to $flip(-K)$, then the Lemma follows. As well
known, we have to choose a knot $K$ in each free homotopy class of oriented
loops in $T^2 \times \mathbb{R}$. Notice that  "$flip \circ -$" acts as the identity on $\pi_1(T^2)$.
\par Any primitive class in $H_1(T^2) \cong \pi_1(T^2)$ can be represented by a torus knot which is invariant under "$flip \circ -$".
\par Let $K \hookrightarrow T^2$. Each class $n[K], n \not= 0$ can be
represented as the closure $\hat \beta$ of the braid $\beta= \sigma_1\sigma_2
\cdots \sigma_{n-1}$ in a tubular neighbourhood $V$ of $K \hookrightarrow T^2
\times \mathbb{R}$ which is a solid torus. (Remember that $\sigma_i$ are the
standard generators of $B_n$.) We easily see that $flip(- \hat \beta)= \hat
\gamma$, where $\gamma= \sigma_{n-1}\sigma_{n-2} \cdots \sigma_2\sigma_1$ in
the same (invariant under $flip$) solid torus $V$. But, as well known, $\beta$
is conjugate to $\gamma$ in $V$ and hence, $\hat \beta$ and $\hat \gamma$ are the same knot.
$\Box$
\par 
{\bf Remarks}.
\begin{enumerate}
\item
Evidently, Lemma~6.1 is still true if one replaces $L$ by any cable of $L$.
\item
Lemma~6.1 implies that the above quantum invariants (and possibly all quantum
invariants) can never detect non-invertibility of knots in $T^2 \times \mathbb{R}$.
\item
It was already well known that quantum invariants never detect the
non-invertibility of links in $S^3$ (see e.g. \cite{K}).
\end{enumerate}
We have shown in sect.~5 that $T$-invariants detect the non-invertibility of
knots in $T^2 \times \mathbb{R}$. Thus, these $T$-invariants (of degrees 5 and
6 in the examples) cannot be extracted from the HOMFLY-PT or Kauffman polynomials of the knot or any of its cables.
\newpage
\section{Non-invertibility of links in $S^3$}
Our results about knots in $T^2 \times \mathbb{R}$ can be interpreted as
results about certain links in $S^3$. Let $T^2 \times \mathbb{R}$ be the tubular neighbourhood of the standardly embedded torus in $S^3$.
\par Let $T_1$ and $T_2$ be the cores of the corresponding solid tori $S^3
\setminus T^2$. To each knot $K \hookrightarrow T^2 \times \mathbb{R} \hookrightarrow S^3$, we associate the link $K \cup T_1 \cup T_2 \hookrightarrow S^3$.
\begin{lem}
{\it Two knots $K, K' \hookrightarrow T^2 \times \mathbb{R}$ are isotopic if and only if the corresponding ordered links $K \cup T_1 \cup T_2$, $K' \cup T_1 \cup T_2 \hookrightarrow S^3$ are isotopic.\/}
\end{lem}
{\bf Proof.}
Lemma~1.7 of \cite{F} implies that the ordered links $K \cup T_1 \cup T_2$ and
$K' \cup T_1 \cup T_2$ are isotopic if and only if the ordered links $K \cup
T_1$ and $K \cup T_2$ are isotopic in the solid torus $S^3 \setminus T_2$. It
is also well known that each isotopy of the solid torus, which is the identity
near the boundary and which maps the core of the solid torus to itself, can be
isotopically deformed to an isotopy which leaves the core pointwise fixed. $\Box$
\par We will use Lemma~7.1 in order to study the invertibility of the link $K
\cup T_1 \cup T_2 \hookrightarrow S^3$. Instead of Lemma~7.1, we could use the
fact that, there is only one isotopy which inverts the Hopf link $H= T_1 \cup
T_2$, up to isotopy of isotopies. Indeed, an isotopy which inverts $H$ inverts
also the meridians and longitudes for $T_1$ and $T_2$. Therefore, such an
isotopy induces an orientation preserving homeomorphism of the incompressible
torus $T^2$ in $S^3 \setminus H$. This homeomorphism acts as $-1$ on $H_1(T^2;
\mathbb{Z})$. As the mapping class group of $T^2$ is $SL(2; \mathbb{Z})$, this
homeomorphism is isotopic to $flip$.
\par {\it Thus, the non-invertibility of the link $L$ in Fig.~1 follows from
  the non-invertibility of the knot $K$ in Fig.~65.\/} Indeed, to the knot
$K$, we have to add the Hopf link $T_1 \cup T_2$. Notice that the resulting
link $L$ is naturally ordered: $K \hookrightarrow S^3$ is not the trivial
knot, $lk(K, T_2)= 4$, $lk(K, T_1)= 1$. Hence, if $L$ is invertible, then $L$
is invertible as an ordered link (i.e. respecting the ordering), and one can
apply Lemma~7.1. $flip$, seen as an involution on $S^3$, maps simultaneously
$T_1$ to $-T_1$ and $T_2$ to $-T_2$. Thus, $L$ is isotopic to $-L$ if and only
if $K$ is isotopic to $flip(-K)$ in $T^2 \times \mathbb{R}$. But we have shown
that this is not the case, using the $T$-invariant in Example~5.1.
\newpage
\section{$T$-invariants which are not of finite type are usefull too}
Let $h= (id, -id): T^2 \times \mathbb{R} \to T^2 \times \mathbb{R}$, and let
$K \hookrightarrow T^2 \times \mathbb{R}$ be a global knot. Then $h(K)$ is
called the {\em mirror image\/} of $K$ and is denoted as usually by $K!$
Clearly, $K!$ is a global knot which is always homotopic to $K$. We give an
example of a $\mathbb{Z}/2\mathbb{Z}$-pure global knot $K$ which we
distinguish from $K!$. We do this in two ways: first with a $T$-invariant of
degree $2$ but which is not of finite type, and then with a $T$-invariant of
degree $8$ which is of finite type. We prove moreover that $K$ and $K!$ cannot
be distinguished by any Gauss diagram invariant (see \cite{F}), or by a
$T$-invariant of finite type of degree not bigger than 2. {\em Hence,
  $T$-invariants which are not of finite type are sometimes more effective than $T$-invariants of finite type.\/}\\
\enlargethispage{-\baselineskip}
\enlargethispage{-\baselineskip}
\enlargethispage{-\baselineskip}
{\it Example~8.1}
\enlargethispage{-\baselineskip}
\enlargethispage{-\baselineskip}
\enlargethispage{-\baselineskip}
$$
\begin{picture}(0,0)%
\includegraphics{im38.pstex}%
\end{picture}%
\setlength{\unitlength}{4144sp}%
\begingroup\makeatletter\ifx\SetFigFont\undefined%
\gdef\SetFigFont#1#2#3#4#5{%
  \reset@font\fontsize{#1}{#2pt}%
  \fontfamily{#3}\fontseries{#4}\fontshape{#5}%
  \selectfont}%
\fi\endgroup%
\begin{picture}(5388,2963)(204,-2834)
\put(5318,-1290){\makebox(0,0)[lb]{\smash{\SetFigFont{12}{14.4}{\rmdefault}{\mddefault}{\updefault}$K$}}}
\put(1732,-2479){\makebox(0,0)[lb]{\smash{\SetFigFont{12}{14.4}{\rmdefault}{\mddefault}{\updefault}$+$}}}
\put(1956,-2460){\makebox(0,0)[lb]{\smash{\SetFigFont{12}{14.4}{\rmdefault}{\mddefault}{\updefault}$+$}}}
\put(3712,-2435){\makebox(0,0)[lb]{\smash{\SetFigFont{12}{14.4}{\rmdefault}{\mddefault}{\updefault}$+$}}}
\put(4110,-2348){\makebox(0,0)[lb]{\smash{\SetFigFont{12}{14.4}{\rmdefault}{\mddefault}{\updefault}$+$}}}
\put(3083,-2024){\makebox(0,0)[lb]{\smash{\SetFigFont{12}{14.4}{\rmdefault}{\mddefault}{\updefault}$-$}}}
\put(3431,-1900){\makebox(0,0)[lb]{\smash{\SetFigFont{12}{14.4}{\rmdefault}{\mddefault}{\updefault}$-$}}}
\put(2236,-1663){\makebox(0,0)[lb]{\smash{\SetFigFont{12}{14.4}{\rmdefault}{\mddefault}{\updefault}$-$}}}
\put(2435,-1788){\makebox(0,0)[lb]{\smash{\SetFigFont{12}{14.4}{\rmdefault}{\mddefault}{\updefault}$-$}}}
\end{picture}
$$
\enlargethispage{-\baselineskip}
\begin{center}
{\bf Fig. 76}
\end{center}
The Gauss diagram of the knot in Fig.~76 is shown in Fig.~77.
$$
\begin{picture}(0,0)%
\includegraphics{im39.pstex}%
\end{picture}%
\setlength{\unitlength}{4144sp}%
\begingroup\makeatletter\ifx\SetFigFont\undefined%
\gdef\SetFigFont#1#2#3#4#5{%
  \reset@font\fontsize{#1}{#2pt}%
  \fontfamily{#3}\fontseries{#4}\fontshape{#5}%
  \selectfont}%
\fi\endgroup%
\begin{picture}(3620,3079)(285,-2684)
\put(3698,-2054){\makebox(0,0)[lb]{\smash{\SetFigFont{12}{14.4}{\rmdefault}{\mddefault}{\updefault}$2\alpha+\beta$}}}
\put(2002,109){\makebox(0,0)[rb]{\smash{\SetFigFont{12}{14.4}{\rmdefault}{\mddefault}{\updefault}$\alpha$}}}
\put(2157,239){\makebox(0,0)[rb]{\smash{\SetFigFont{12}{14.4}{\rmdefault}{\mddefault}{\updefault}$\alpha$}}}
\put(3422, -1){\makebox(0,0)[lb]{\smash{\SetFigFont{12}{14.4}{\rmdefault}{\mddefault}{\updefault}$\alpha$}}}
\put(3808,-479){\makebox(0,0)[lb]{\smash{\SetFigFont{12}{14.4}{\rmdefault}{\mddefault}{\updefault}$2\alpha+\beta$}}}
\put(1441,-645){\makebox(0,0)[lb]{\smash{\SetFigFont{12}{14.4}{\rmdefault}{\mddefault}{\updefault}$+$}}}
\put(2213,-15){\makebox(0,0)[lb]{\smash{\SetFigFont{12}{14.4}{\rmdefault}{\mddefault}{\updefault}$+$}}}
\put(2468,  7){\makebox(0,0)[lb]{\smash{\SetFigFont{12}{14.4}{\rmdefault}{\mddefault}{\updefault}$+$}}}
\put(1441,-1216){\makebox(0,0)[lb]{\smash{\SetFigFont{12}{14.4}{\rmdefault}{\mddefault}{\updefault}$+$}}}
\put(2671,-923){\makebox(0,0)[lb]{\smash{\SetFigFont{12}{14.4}{\rmdefault}{\mddefault}{\updefault}$-$}}}
\put(3661,-1223){\makebox(0,0)[lb]{\smash{\SetFigFont{12}{14.4}{\rmdefault}{\mddefault}{\updefault}$-$}}}
\put(2567,-2626){\makebox(0,0)[lb]{\smash{\SetFigFont{12}{14.4}{\rmdefault}{\mddefault}{\updefault}$\alpha$}}}
\put(3105,-106){\makebox(0,0)[rb]{\smash{\SetFigFont{12}{14.4}{\rmdefault}{\mddefault}{\updefault}$3$}}}
\put(3226,-376){\makebox(0,0)[lb]{\smash{\SetFigFont{12}{14.4}{\rmdefault}{\mddefault}{\updefault}$4$}}}
\put(3601,-729){\makebox(0,0)[rb]{\smash{\SetFigFont{12}{14.4}{\rmdefault}{\mddefault}{\updefault}$-$}}}
\put(2910,-1179){\makebox(0,0)[lb]{\smash{\SetFigFont{12}{14.4}{\rmdefault}{\mddefault}{\updefault}$8$}}}
\put(3458,-1014){\makebox(0,0)[rb]{\smash{\SetFigFont{12}{14.4}{\rmdefault}{\mddefault}{\updefault}$1$}}}
\put(2678,-1695){\makebox(0,0)[rb]{\smash{\SetFigFont{12}{14.4}{\rmdefault}{\mddefault}{\updefault}$7$}}}
\put(2588,-1441){\makebox(0,0)[rb]{\smash{\SetFigFont{12}{14.4}{\rmdefault}{\mddefault}{\updefault}$-$}}}
\put(2724,-2161){\makebox(0,0)[rb]{\smash{\SetFigFont{12}{14.4}{\rmdefault}{\mddefault}{\updefault}$2$}}}
\put(2401,-1966){\makebox(0,0)[lb]{\smash{\SetFigFont{12}{14.4}{\rmdefault}{\mddefault}{\updefault}$6$}}}
\put(1943,-1749){\makebox(0,0)[rb]{\smash{\SetFigFont{12}{14.4}{\rmdefault}{\mddefault}{\updefault}$5$}}}
\put(1756,-2491){\makebox(0,0)[b]{\smash{\SetFigFont{12}{14.4}{\rmdefault}{\mddefault}{\updefault}$2\alpha+\beta$}}}
\put(901,-826){\makebox(0,0)[b]{\smash{\SetFigFont{12}{14.4}{\rmdefault}{\mddefault}{\updefault}$2\alpha+\beta$}}}
\end{picture}
$$
\begin{center}
{\bf Fig. 77}
\end{center}
\enlargethispage{-2cm}
For the convenience of the reader, we have affected numbers to the
crossings. We see that there appear only two different homology classes as
markings. In particular, $K$ is $\mathbb{Z}/ 2\mathbb{Z}$-pure for $G:=
(H_1(T^2)/\langle [K] \rangle)/2 \mathbb{Z}= \{ 0, a \}$. The Gauss diagram of
$K!$ is obtained from the one of $K$ by replacing all arrows, writhes and
markings by their opposites (but remember that $a= -a$). We start by comparing
the invariants of degree 1 (see also \cite{F}, sect.~2.2).
$$
W_K(\alpha)= W_{K!}(\alpha)= W_K(2\alpha+ \beta)= W_{K!}(2\alpha+ \beta)=0
$$
If we see $K$ and $K!$ as knots in $S^3$ using the embedding $T^2 \times
\mathbb{R} \hookrightarrow S^3$, then we easily calculate $v_2(K)= v_2(K!)$
for the only Vassiliev invariant of degree 2 (notice that $K!$ is not the
mirror image of $K$ in $S^3$ because of the two additional crossings seen in
Fig.~76). All the Gauss diagram invariants and $T$-invariants of degree 2
which are of finite type are linear combinations of all possible
configurations of degree 2 (see also \cite{F}, sect.~2.4). The weight function
is always the product of the two writhes (because of the invariance under
Reidemeister moves of type $II$). Therefore, this function is invariant under
taking the mirror image. In Fig.~78, we indicate how the configurations change
by taking the mirror image. 
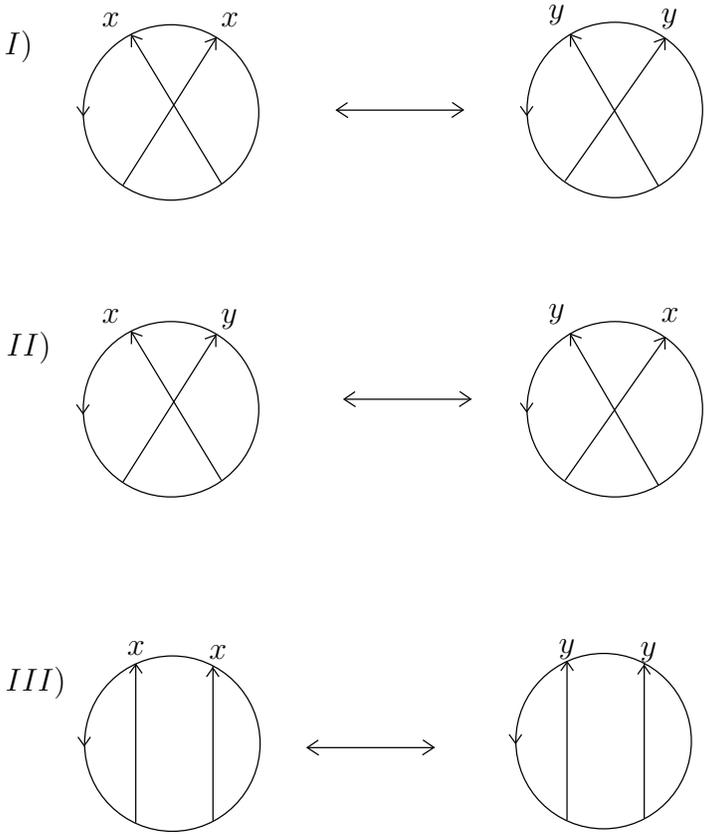
\enlargethispage{-2cm}
\vspace{2cm}
$$
\begin{picture}(0,0)%
\includegraphics{ima78.pstex}%
\end{picture}%
\setlength{\unitlength}{4144sp}%
\begingroup\makeatletter\ifx\SetFigFont\undefined%
\gdef\SetFigFont#1#2#3#4#5{%
  \reset@font\fontsize{#1}{#2pt}%
  \fontfamily{#3}\fontseries{#4}\fontshape{#5}%
  \selectfont}%
\fi\endgroup%
\begin{picture}(4180,4994)(466,-4605)
\put(1052,224){\makebox(0,0)[lb]{\smash{\SetFigFont{12}{14.4}{\rmdefault}{\mddefault}{\updefault}$x$}}}
\put(1765,224){\makebox(0,0)[lb]{\smash{\SetFigFont{12}{14.4}{\rmdefault}{\mddefault}{\updefault}$x$}}}
\put(4397,247){\makebox(0,0)[lb]{\smash{\SetFigFont{12}{14.4}{\rmdefault}{\mddefault}{\updefault}$y$}}}
\put(3722,269){\makebox(0,0)[lb]{\smash{\SetFigFont{12}{14.4}{\rmdefault}{\mddefault}{\updefault}$y$}}}
\put(1052,-1553){\makebox(0,0)[lb]{\smash{\SetFigFont{12}{14.4}{\rmdefault}{\mddefault}{\updefault}$x$}}}
\put(1765,-1553){\makebox(0,0)[lb]{\smash{\SetFigFont{12}{14.4}{\rmdefault}{\mddefault}{\updefault}$y$}}}
\put(4397,-1545){\makebox(0,0)[lb]{\smash{\SetFigFont{12}{14.4}{\rmdefault}{\mddefault}{\updefault}$x$}}}
\put(3722,-1523){\makebox(0,0)[lb]{\smash{\SetFigFont{12}{14.4}{\rmdefault}{\mddefault}{\updefault}$y$}}}
\put(1203,-3540){\makebox(0,0)[lb]{\smash{\SetFigFont{12}{14.4}{\rmdefault}{\mddefault}{\updefault}$x$}}}
\put(1691,-3562){\makebox(0,0)[lb]{\smash{\SetFigFont{12}{14.4}{\rmdefault}{\mddefault}{\updefault}$x$}}}
\put(3783,-3525){\makebox(0,0)[lb]{\smash{\SetFigFont{12}{14.4}{\rmdefault}{\mddefault}{\updefault}$y$}}}
\put(4271,-3547){\makebox(0,0)[lb]{\smash{\SetFigFont{12}{14.4}{\rmdefault}{\mddefault}{\updefault}$y$}}}
\put(466, 54){\makebox(0,0)[lb]{\smash{\SetFigFont{12}{14.4}{\rmdefault}{\mddefault}{\updefault}$I)$}}}
\put(481,-1745){\makebox(0,0)[lb]{\smash{\SetFigFont{12}{14.4}{\rmdefault}{\mddefault}{\updefault}$II)$}}}
\put(473,-3755){\makebox(0,0)[lb]{\smash{\SetFigFont{12}{14.4}{\rmdefault}{\mddefault}{\updefault}$III)$}}}
\end{picture}
$$
$$
\begin{picture}(0,0)%
\includegraphics{ima78bis.pstex}%
\end{picture}%
\setlength{\unitlength}{4144sp}%
\begingroup\makeatletter\ifx\SetFigFont\undefined%
\gdef\SetFigFont#1#2#3#4#5{%
  \reset@font\fontsize{#1}{#2pt}%
  \fontfamily{#3}\fontseries{#4}\fontshape{#5}%
  \selectfont}%
\fi\endgroup%
\begin{picture}(4112,7400)(443,-7255)
\put(458,-163){\makebox(0,0)[lb]{\smash{\SetFigFont{12}{14.4}{\rmdefault}{\mddefault}{\updefault}$IV)$}}}
\put(1185, 10){\makebox(0,0)[lb]{\smash{\SetFigFont{12}{14.4}{\rmdefault}{\mddefault}{\updefault}$x$}}}
\put(1673,-12){\makebox(0,0)[lb]{\smash{\SetFigFont{12}{14.4}{\rmdefault}{\mddefault}{\updefault}$y$}}}
\put(3735, -5){\makebox(0,0)[lb]{\smash{\SetFigFont{12}{14.4}{\rmdefault}{\mddefault}{\updefault}$x$}}}
\put(4223,-27){\makebox(0,0)[lb]{\smash{\SetFigFont{12}{14.4}{\rmdefault}{\mddefault}{\updefault}$y$}}}
\put(4342,-1948){\makebox(0,0)[lb]{\smash{\SetFigFont{12}{14.4}{\rmdefault}{\mddefault}{\updefault}$y$}}}
\put(3449,-3125){\makebox(0,0)[lb]{\smash{\SetFigFont{12}{14.4}{\rmdefault}{\mddefault}{\updefault}$y$}}}
\put(1177,-1866){\makebox(0,0)[lb]{\smash{\SetFigFont{12}{14.4}{\rmdefault}{\mddefault}{\updefault}$x$}}}
\put(1784,-3013){\makebox(0,0)[lb]{\smash{\SetFigFont{12}{14.4}{\rmdefault}{\mddefault}{\updefault}$x$}}}
\put(458,-2233){\makebox(0,0)[lb]{\smash{\SetFigFont{12}{14.4}{\rmdefault}{\mddefault}{\updefault}$V)$}}}
\put(451,-4026){\makebox(0,0)[lb]{\smash{\SetFigFont{12}{14.4}{\rmdefault}{\mddefault}{\updefault}$VI)$}}}
\put(4365,-3771){\makebox(0,0)[lb]{\smash{\SetFigFont{12}{14.4}{\rmdefault}{\mddefault}{\updefault}$x$}}}
\put(3472,-4948){\makebox(0,0)[lb]{\smash{\SetFigFont{12}{14.4}{\rmdefault}{\mddefault}{\updefault}$x$}}}
\put(1200,-3689){\makebox(0,0)[lb]{\smash{\SetFigFont{12}{14.4}{\rmdefault}{\mddefault}{\updefault}$y$}}}
\put(1807,-4836){\makebox(0,0)[lb]{\smash{\SetFigFont{12}{14.4}{\rmdefault}{\mddefault}{\updefault}$y$}}}
\put(443,-6288){\makebox(0,0)[lb]{\smash{\SetFigFont{12}{14.4}{\rmdefault}{\mddefault}{\updefault}$VII)$}}}
\put(4366,-6033){\makebox(0,0)[lb]{\smash{\SetFigFont{12}{14.4}{\rmdefault}{\mddefault}{\updefault}$x$}}}
\put(3473,-7210){\makebox(0,0)[lb]{\smash{\SetFigFont{12}{14.4}{\rmdefault}{\mddefault}{\updefault}$y$}}}
\put(1201,-5951){\makebox(0,0)[lb]{\smash{\SetFigFont{12}{14.4}{\rmdefault}{\mddefault}{\updefault}$x$}}}
\put(1808,-7098){\makebox(0,0)[lb]{\smash{\SetFigFont{12}{14.4}{\rmdefault}{\mddefault}{\updefault}$y$}}}
\end{picture}
$$
\begin{center}
{\bf Fig. 78}
\end{center}
In this Figure, $x, y \in \{ \alpha, 2\alpha+ \beta \}$ and $x \not= y$.
\par {\bf Remarks}.
\begin{enumerate}
\item
$I)$ can not enter in any invariant because of the invariance under Reidemeister moves of type $II$ (see also Lemma~3.3).
\item
$IV)$ is invariant.
\end{enumerate}
Thus, if in our example, the left-hand side is equal to the right-hand side
for $II)$, $III)$, $V)$, $VI)$, $VII)$, then {\em all\/} invariants of finite
type of degree 2 coincide for $K$ and $K!$. We easily calculate the values on both sides and it turns out that they
coincide: $II)= +2$, $III)= 0$, $V)= 0$, $VI)= -2$, $VII)= 0$.
\par Let us consider $T$-invariants of infinite type (see Prop.~3.1).
Let 
$$
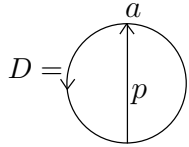
\makebox(10,60){ $D=$ } \makebox(60,60){\begin{picture}(0,0)%
\includegraphics{pag52M.pstex}%
\end{picture}%
\setlength{\unitlength}{4144sp}%
\begingroup\makeatletter\ifx\SetFigFont\undefined%
\gdef\SetFigFont#1#2#3#4#5{%
  \reset@font\fontsize{#1}{#2pt}%
  \fontfamily{#3}\fontseries{#4}\fontshape{#5}%
  \selectfont}%
\fi\endgroup%
\begin{picture}(772,947)(494,-258)
\put(931, 39){\makebox(0,0)[lb]{\smash{\SetFigFont{12}{14.4}{\rmdefault}{\mddefault}{\updefault}$p$}}}
\put(893,524){\makebox(0,0)[lb]{\smash{\SetFigFont{12}{14.4}{\rmdefault}{\mddefault}{\updefault}$a$}}}
\end{picture}
}
$$
for $G \cong \mathbb{Z}/2\mathbb{Z}= \{ 0, a \}$. We easily calculate
$$
T_K(D; c_{++}(D)= +2)= -1
$$
($D$ is the crossing number 1 with the crossing $q$ which is either the crossing number 4 or 5.)
$$
T_K(D; c_{++}(D)= -1)= +2
$$
($D$ is 4 and 5, $q$ is 1.)
$$
T_K(D; c_{++}(D)= 0)= -1
$$
($D$ is 2, 3, 6, 7, 8. There are no $q$'s.) For all other $c \in \mathbb{Z}$,
$T_K(D; c_{++}(D)= c)= 0$ We keep the numbers for the crossings of $K!$.
$$
T_{K!}(D; c_{++}(D)= -2)= +1
$$
($D$ is 2, the $q$'s are 3 and 6.)
$$
T_{K!}(D; c_{++}(D)= +1)= -2
$$
($D$ is 3 and 6, $q$ is 2.)
$$
T_{K!}(D; c_{++}(D)= 0)= +1
$$
and all other $T_{K!}(D; c_{++}(D)= c)= 0$. Consequently, $K$ and $K!$ are not
isotopic and we have proven it with an invariant of quadratic complexity. The
Gauss diagram of $K$ without the writhes is a $\mathbb{Z}/2\mathbb{Z}$-pure
configuration $D$ of degree 8. Clearly, it is different from the corresponding
configuration for $K!$. Thus, $K$ and $K!$ are also distinguished by the
$T$-invariant of degree 8 of finite type $T_K(D; \emptyset)$. We do not know
wether or not there are $G$-pure global knots which can be distinguished {\em
  only\/} by $T$-invariants which are not of finite type. But in any case, our
example shows that these invariants do it sometimes in a more effective way than the invariants of finite type.
\newpage
\section{$T$-invariants are not well defined for general knots}
Let $K \hookrightarrow S^3= (\mathbb{R}^2 \times \mathbb{R}) \cup \{ \infty \}$ be a knot and let $m$ be a meridian of $K$. The meridian
$m$ is isotopic to $(0 \times \mathbb{R}) \cup \{ \infty \}$ and hence, we can
consider $K$ as a knot in $(\mathbb{R}^2 \setminus 0) \times \mathbb{R}$. If
two knots $K$, $K'$ are isotopic in $S^3$, then in fact, they are already
isotopic in $(\mathbb{R}^2 \setminus 0) \times \mathbb{R}$, where $(0 \times
\mathbb{R}) \cup \{ \infty \}$ is a meridian for both knots. We consider the
projection $(\mathbb{R}^2 \setminus 0) \times \mathbb{R} \to \mathbb{R}^2
\setminus 0$. Assume that $K$ and $K'$ are isotopic. If they have the same
writhe and the same Whitney index, then they are regularly isotopic in $(\mathbb{R}^2 \setminus 0) \times \mathbb{R}$ (see  \cite{F}, sect.~2).
If we take now the same cable or satellite for two regularly isotopic knots,
then the resulting knots are again (regularly) isotopic in $(\mathbb{R}^2 \setminus 0) \times \mathbb{R}$.
\par Let $K$ be the figure-eight knot. As well known, $K$ is isotopic to its
mirror image $K!$. As satellite, we take the positive (untwisted) Whitehead
double. Consequently, the two knots $W$ and $W'$ shown in Fig.~79 are isotopic in the solid torus $S^3 \setminus m= (\mathbb{R}^2 \setminus 0) \times \mathbb{R}$.
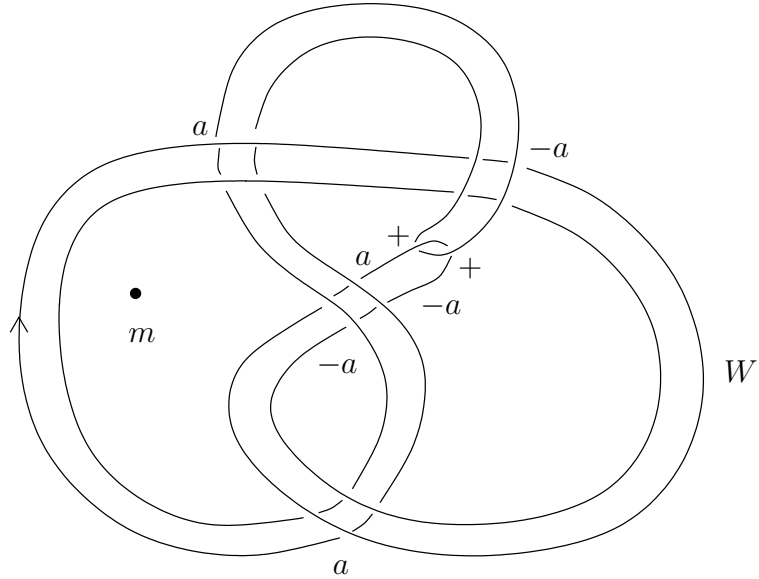
$$
\begin{picture}(0,0)%
\includegraphics{ima79.pstex}%
\end{picture}%
\setlength{\unitlength}{4144sp}%
\begingroup\makeatletter\ifx\SetFigFont\undefined%
\gdef\SetFigFont#1#2#3#4#5{%
  \reset@font\fontsize{#1}{#2pt}%
  \fontfamily{#3}\fontseries{#4}\fontshape{#5}%
  \selectfont}%
\fi\endgroup%
\begin{picture}(4287,3420)(769,-3076)
\put(1973,-450){\makebox(0,0)[rb]{\smash{\SetFigFont{12}{14.4}{\rmdefault}{\mddefault}{\updefault}$a$}}}
\put(3885,-579){\makebox(0,0)[lb]{\smash{\SetFigFont{12}{14.4}{\rmdefault}{\mddefault}{\updefault}$-a$}}}
\put(3466,-1290){\makebox(0,0)[lb]{\smash{\SetFigFont{12}{14.4}{\rmdefault}{\mddefault}{\updefault}$+$}}}
\put(3241,-1516){\makebox(0,0)[lb]{\smash{\SetFigFont{12}{14.4}{\rmdefault}{\mddefault}{\updefault}$-a$}}}
\put(3182,-1097){\makebox(0,0)[rb]{\smash{\SetFigFont{12}{14.4}{\rmdefault}{\mddefault}{\updefault}$+$}}}
\put(2948,-1223){\makebox(0,0)[rb]{\smash{\SetFigFont{12}{14.4}{\rmdefault}{\mddefault}{\updefault}$a$}}}
\put(2716,-3076){\makebox(0,0)[lb]{\smash{\SetFigFont{12}{14.4}{\rmdefault}{\mddefault}{\updefault}$a$}}}
\put(5056,-1928){\makebox(0,0)[lb]{\smash{\SetFigFont{12}{14.4}{\rmdefault}{\mddefault}{\updefault}$W$}}}
\put(1493,-1688){\makebox(0,0)[lb]{\smash{\SetFigFont{12}{14.4}{\rmdefault}{\mddefault}{\updefault}$m$}}}
\put(2619,-1876){\makebox(0,0)[lb]{\smash{\SetFigFont{12}{14.4}{\rmdefault}{\mddefault}{\updefault}$-a$}}}
\end{picture}
$$
$$
\begin{picture}(0,0)%
\includegraphics{ima79bis.pstex}%
\end{picture}%
\setlength{\unitlength}{4144sp}%
\begingroup\makeatletter\ifx\SetFigFont\undefined%
\gdef\SetFigFont#1#2#3#4#5{%
  \reset@font\fontsize{#1}{#2pt}%
  \fontfamily{#3}\fontseries{#4}\fontshape{#5}%
  \selectfont}%
\fi\endgroup%
\begin{picture}(4310,3424)(746,-3076)
\put(1973,-450){\makebox(0,0)[rb]{\smash{\SetFigFont{12}{14.4}{\rmdefault}{\mddefault}{\updefault}$-a$}}}
\put(3885,-579){\makebox(0,0)[lb]{\smash{\SetFigFont{12}{14.4}{\rmdefault}{\mddefault}{\updefault}$a$}}}
\put(3466,-1290){\makebox(0,0)[lb]{\smash{\SetFigFont{12}{14.4}{\rmdefault}{\mddefault}{\updefault}$+$}}}
\put(3241,-1516){\makebox(0,0)[lb]{\smash{\SetFigFont{12}{14.4}{\rmdefault}{\mddefault}{\updefault}$-a$}}}
\put(3182,-1097){\makebox(0,0)[rb]{\smash{\SetFigFont{12}{14.4}{\rmdefault}{\mddefault}{\updefault}$+$}}}
\put(2716,-3076){\makebox(0,0)[lb]{\smash{\SetFigFont{12}{14.4}{\rmdefault}{\mddefault}{\updefault}$-a$}}}
\put(5056,-1928){\makebox(0,0)[lb]{\smash{\SetFigFont{12}{14.4}{\rmdefault}{\mddefault}{\updefault}$W'$}}}
\put(1493,-1688){\makebox(0,0)[lb]{\smash{\SetFigFont{12}{14.4}{\rmdefault}{\mddefault}{\updefault}$m$}}}
\put(2948,-1223){\makebox(0,0)[rb]{\smash{\SetFigFont{12}{14.4}{\rmdefault}{\mddefault}{\updefault}$a$}}}
\put(2701,-1853){\makebox(0,0)[lb]{\smash{\SetFigFont{12}{14.4}{\rmdefault}{\mddefault}{\updefault}$a$}}}
\end{picture}
$$
\begin{center}
{\bf Fig. 79}
\end{center}

\par Let $v$ be a Morse-Smale vector field on $\mathbb{R}^2$, which has a
critical point of index 1 in $0= m \cap \mathbb{R}^2$, and such that $v$ is
transversal to $pr(W)$. But $pr(W)= pr(W')$ and hence, $v$ is transversal to
$pr(W')$ too. Of course, $W$ and $W'$ are {\em not\/} global knots, because
$v$ has critical points of index 1 different from 0. Let "$a$" be the generator$$
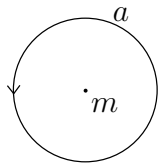
\begin{picture}(0,0)%
\includegraphics{gener.pstex}%
\end{picture}%
\setlength{\unitlength}{4144sp}%
\begingroup\makeatletter\ifx\SetFigFont\undefined%
\gdef\SetFigFont#1#2#3#4#5{%
  \reset@font\fontsize{#1}{#2pt}%
  \fontfamily{#3}\fontseries{#4}\fontshape{#5}%
  \selectfont}%
\fi\endgroup%
\begin{picture}(916,988)(1803,-1285)
\put(2453,-432){\makebox(0,0)[lb]{\smash{\SetFigFont{12}{14.4}{\rmdefault}{\mddefault}{\updefault}$a$}}}
\put(2325,-977){\makebox(0,0)[lb]{\smash{\SetFigFont{12}{14.4}{\rmdefault}{\mddefault}{\updefault}$m$}}}
\end{picture}
$$
of $H_1(S^3 \setminus m; \mathbb{Z})=H_1(\mathbb{R}^2 \setminus 0;
\mathbb{Z})$. One has $[W]= [W']= 0$ in $H_1(\mathbb{R}^2 \setminus 0;
\mathbb{Z})$ and we easily see that $W$ and $W'$ are $\mathbb{Z}$-pure knots
(the markings are shown in Fig.~79 too). It follows from the proof of Theorem~1 that each $T$-invariant $T_W$
\begin{enumerate}
\item
is invariant for each isotopy transversal to $v$ and which is not necessarily $\mathbb{Z}$-pure
\item
is invariant for each $\mathbb{Z}$-pure isotopy which is not necessarily transversal to $v$.
\end{enumerate}
Let 
$$
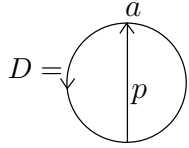
\makebox(10,60){$D=$} \makebox(60,60){\begin{picture}(0,0)%
\includegraphics{pag52M.pstex}%
\end{picture}%
\setlength{\unitlength}{4144sp}%
\begingroup\makeatletter\ifx\SetFigFont\undefined%
\gdef\SetFigFont#1#2#3#4#5{%
  \reset@font\fontsize{#1}{#2pt}%
  \fontfamily{#3}\fontseries{#4}\fontshape{#5}%
  \selectfont}%
\fi\endgroup%
\begin{picture}(772,947)(494,-258)
\put(931, 39){\makebox(0,0)[lb]{\smash{\SetFigFont{12}{14.4}{\rmdefault}{\mddefault}{\updefault}$p$}}}
\put(893,524){\makebox(0,0)[lb]{\smash{\SetFigFont{12}{14.4}{\rmdefault}{\mddefault}{\updefault}$a$}}}
\end{picture}
}
$$
 and let $c(D)$ be any of the classes of Def.~3.7.
We easily calculate $T_W(D; c(D)=c)= T_{W'}(D; c(D)=c)$ for any $c \in
\mathbb{Z}$. But e.g. $T_W(D; c_{++}(D)=0, c_{+-}(D)=-1)= +2$ and $T_{W'}(D; c_{++}(D)=0, c_{+-}(D)=-1)= 0$. We have shown above that $W$ and $W'$ are
actually isotopic in $(\mathbb{R}^2 \setminus 0) \times \mathbb{R}$.
This example has three important consequences:
\begin{enumerate}
\item
$W$ and $W'$ are transversal to $v$ and they are isotopic. But there is no isotopy transversal to $v$ joining them.
\item
$W$ and $W'$ are $\mathbb{Z}$-pure and they are isotopic. But there is
no $\mathbb{Z}$-pure isotopy joining them (i.e. there are cycles with marking 0 which cannot be eliminated).
\item
Multi-classes contain more information than the classes taken individually.
\end{enumerate}

\newpage
Unfortunately, my preprint "New invariants in knot theory" (November 1999)
contains some serious errors. I apologize by the reader for this. The present
preprint is the result of the correction of these errors. It will be added to
the final version of the monography "Gauss diagram invariants for knots and links".

I am very grateful to S\'everine for her constant support!


\newpage
\begin{thebibliography}{99}

\bibitem[A-M-R]{A-M-R} J.~Andersen, J.~Mattes, N.~Reshetikhlin: {\em
    Quantization of the algebra of chord diagrams\/} Math. Proc. Cambridge Philos. Soc. 124 (1998)

\bibitem[B-L]{B-L} J.~Birman, X.S.~Lin: {\em Knot polynomials and Vassiliev
    invariants\/} Inv. Math. 111 (1993)

\bibitem[BN]{BN} D.~Bar Natan: {\em On the Vassiliev knot invariants\/}
  Topology 34 (1995)

\bibitem[F]{F} T.~Fiedler: {\em Gauss diagram invariants for knots and links\/}
Monography, Kluwer academic publishers (to appear)

\bibitem[Go]{Go} V.~Goryunov: {\em Finite order invariants of framed knots in
    a solid torus and in Arnold's $\mathcal{J}^+$-theory of plane curves\/} "Geometry and Physics", Lect. Notes in Pure and Appl. Math. (1996)

\bibitem[G-P-V]{G-P-V} M.~Gussarov, M.~Polyak, O.~Viro: {\em Finite type
    invariants of classical and virtual knots\/} Topology 39 (2000)

\bibitem[H-P]{H-P} J.~Hoste, J.~Przytycki: {\em A survey of skein modules of
    3-manifolds\/} Knots 90 (Osaka, 1990) de Gruyter (1992)

\bibitem[K]{K}M.~Kontsevitch: {\em Vassiliev's knot invariants\/} Adv. in
  Sov. Math. 16 (2) (1993)

\bibitem[L]{L} X.S~Lin: {\em Finite type link invariants and the
    non-invertibility of links\/} Math. Res. Lett. 3 (1996)

\bibitem[M]{M} H.~Morton: {\em Infinitely many fibred knots having the same
    Alexander polynomial\/} Topology 17 (1978)

\bibitem[P-S]{P-S} V.~Praslov, A.~Sossinsky: {\em Knots, links, braids and
    3-manifolds\/} Transl. of Math. Monographs 154 AMS (1997)

\bibitem[P-V]{P-V}M.~Polyak, O.~Viro: {\em Gauss diagram formulas for
    Vassiliev invariants\/} Internat. Math. Res. Notices 11 (1994)

\bibitem[V]{V} V.~Vassiliev: {\em Cohomology of knot spaces\/} Theory of
  Singularities and its Applications Amer. Math. Soc., Providence (1990)

\end{thebibliography}
\end{document}